\documentclass{amsart}
\usepackage{amsmath,amssymb}
\newcommand     {\comment}[1]   {}
\newcommand{\mute}[2] {}
\newcommand     {\printname}[1] {}
\newcommand{\labell}[1] {\label{#1}\printname{#1}}

\def    \calM   {{\mathcal M}}

\newcommand\ft {{\mathfrak t}}

\def    \oursetminus  {{\smallsetminus}}

\def    \sect    {{\operatorname{sec}}}

\def    \inv    {^{-1}}
\def    \to     {\longrightarrow}
\def    \ssminus        {\smallsetminus}

\newcommand     {\ol}[1]   {\overline{#1}}

\renewcommand{\Tilde}{\widetilde}
\usepackage{amssymb,latexsym}
\newcommand{\1}{{{\mathchoice {\rm 1\mskip-4mu l} {\rm 1\mskip-4mu l}
{\rm 1\mskip-4.5mu l} {\rm 1\mskip-5mu l}}}}

\newcommand{\codim}{{\rm codim\ }}
\newcommand{\Fmax}{F_{\rm max}}
\newcommand{\Fmin}{F_{\rm min}}
\newcommand{\coker}{{\rm coker}}

\newcommand{\TW}{{\Tilde W}}
\newcommand{\TY}{{\Tilde Y}}
\newcommand{\TX}{{\Tilde X}}
\newcommand{\TGa}{{\Tilde \Gamma}}

\newcommand{\End}{{\rm End}}

\newcommand{\PSL}{{\rm PSL}}

\newcommand{\pJ}{{{\ov{\p}_{J}}}}
\newcommand{\R}{{\mathbb R}}

\newcommand{\Q}{{\mathbb Q}}
\newcommand{\QH}{{\rm QH}}
\newcommand{\Hom}{{\rm Hom}}

\newcommand{\Flux}{{\rm Flux}}
\newcommand{\Symp}{{\rm Symp}}

\newcommand{\GW}{{\rm GW}}
\newcommand{\TJ}{{\widetilde{J}}}
\newcommand{\TV}{{\widetilde{V}}}
\newcommand{\TZ}{{\widetilde{Z}}}
\newcommand{\TLl}{{\widetilde{\Ll}}}
\newcommand{\TWw}{{\widetilde{\Ww}}}
\newcommand{\TU}{{\widetilde{U}}}
\newcommand{\ttau}{{\tilde\tau}}
\newcommand{\tu}{{\tilde{u}}}
\newcommand{\ts}{{\tilde{s}}}

\newcommand{\pbar}{{\ov {\p}}}
\newcommand{\ev}{{\rm ev}}
\newcommand{\odd}{{\rm odd}}
\newcommand{\Hor}{{\rm Hor}}
\newcommand{\ver}{{\rm vert}}
\newcommand{\PD}{{\rm PD}}

\newcommand{\Z}{{\mathbb Z}}
\newcommand{\C}{{\mathbb C}}

\newcommand{\ov}{\overline }
\newcommand{\oMm}{{\overline {\Mm}}}
\newcommand{\oMmn}{{\overline {\Mm}\,\!^{\nu}}}
\newcommand{\oMmnk}{{\overline {\Mm}\,\!^{\nu}_{0,k}}}

\newcommand{\Ww}{{\mathcal {W}}}

\newcommand{\bw}{{\bf w}}

\newcommand{\bu}{{\bf u}}
\newcommand{\bz}{{\bf z}}

\newcommand{\p}{{\partial}}
\newcommand{\al}{{\alpha}}

\newcommand{\be}{{\beta}}

\newcommand{\Om}{{\Omega}}
\newcommand{\om}{{\omega}}
\newcommand{\io}{{\iota}}
\newcommand{\eps}{{\varepsilon}}

\newcommand{\ga}{{\gamma}}
\newcommand{\Ga}{{\Gamma}}
\newcommand{\ka}{{\kappa}}

\newcommand{\si}{{\sigma}}

\newcommand{\Aa}{{\mathcal A}}
\newcommand{\Bb}{{\mathcal B}}
\newcommand{\Cc}{{\mathcal C}}
\newcommand{\Dd}{{\mathcal D}}
\newcommand{\Ll}{{\mathcal L}}
\newcommand{\Jj}{{\mathcal J}}
\newcommand{\Ee}{{\mathcal E}}
\newcommand{\Ff}{{\mathcal F}}

\newcommand{\Nn}{{\mathcal N}}
\newcommand{\Mm}{{\mathcal M}}

\newcommand{\Ss}{{\mathcal S}}

\newcommand{\La}{{\Lambda}}
\newcommand{\Aut}{{\rm Aut}}
\newcommand{\grad}{{\rm grad\,}}
\newcommand{\pr}{{\rm pr}}

\newcommand{\Si}{{\Sigma}}
\newcommand{\Ham}{{\rm Ham}}
\newcommand{\Stab}{{\rm Stab}}

\newcommand{\im}{{\rm im }}

\newcommand{\MS}{{\medskip}}

\newcommand{\NI}{{\noindent}}

\newcommand{\QED}{\hfill$\Box$\medskip}

\newcommand{\CP}{{\mathbb{CP }}}

\newtheorem{theorem}{Theorem}[section]
\newtheorem{cor}[theorem]{Corollary}
\newtheorem{corollary}[theorem]{Corollary}
\newtheorem{defn}[theorem]{Definition}
\newtheorem{thm}[theorem]{Theorem}
\newtheorem{definition}[theorem]{Definition}
\newtheorem{example}[theorem]{Example}

\newtheorem{remark}[theorem]{Remark}
\newtheorem{rmk}[theorem]{Remark}
\newtheorem{lemma}[theorem]{Lemma}

\newtheorem{prop}[theorem]{Proposition}
\newtheorem{proposition}[theorem]{Proposition}

\begin{document}

\title{Topological properties of Hamiltonian circle actions}
\author{Dusa McDuff}
\thanks{First author partially supported by NSF grant DMS 0305939, 
and second by NSF grant DMS 0204448.}
\address{Department of Mathematics,
 Stony Brook University, Stony Brook, 
NY 11794-3651, USA}
\email{dusa@math.sunysb.edu}
\author{Susan Tolman}
\address{Department of Mathematics,
 University of Illinois at Urbana--Champaign, 
IL, USA}
\email{tolman@math.uiuc.edu}
\date{April 14 2004, revised May 9 2005}

\maketitle

\begin{abstract}  
This paper studies  Hamiltonian circle actions, i.e. circle subgroups of the group $\Ham(M,\om)$ of Hamiltonian symplectomorphisms of a closed symplectic manifold $(M,\om)$.  Our main tool is the  Seidel 
representation of $\pi_1(\Ham(M,\om))$ in the units of the quantum
homology ring.  We show that if the weights of the action at the points at which the moment 
map is a maximum are sufficiently small then the circle 
represents a nonzero element of $\pi_1(\Ham(M,\om))$.  Further, if  the isotropy has order at most two and the circle contracts in 
$\Ham(M,\om)$ then the homology of $M$ is invariant under an involution. For example, the image of the normalized moment map is a symmetric interval $[-a,a]$.
If the action is semifree (i.e. the isotropy weights are $0$ or $\pm 1$) then we calculate the leading order term in
the Seidel representation, an important technical tool in understanding the quantum cohomology of manifolds that admit semifree Hamiltonian circle actions.
If the manifold is toric, we use our results about 
this representation to
describe the basic multiplicative structure of the 
quantum cohomology ring of an arbitrary toric manifold.
There are two important technical ingredients;
one relates the equivariant cohomology  of $M$ to the Morse flow of 
the moment map, and the other is
a version of the localization principle
for calculating Gromov--Witten invariants on 
symplectic manifolds with $S^1$-actions.
\end{abstract}


\section{Introduction}

This paper grew out of an attempt to 
understand when a circle action on a 
symplectic manifold $(M, \om)$ gives rise to an {\bf essential}
(i.e. noncontractible)
loop in the symplectomorphism 
group $\Symp(M, \om)$.  Since nonHamiltonian loops have nonzero image under the Flux homomorphism
$$
\Flux: \pi_1(\Symp(M,\om)) \to H^1(M,\R),
$$
it suffices to restrict attention to circle subgroups of the Hamiltonian group $\Ham: = \Ham(M,\om)$.  
These actions are generated by functions $K:M\to \R$, the {\bf moment map}.  We shall always assume that $K$ is 
normalized, i.e. that $\int_M K\om^n = 0$, and shall denote the
 corresponding circle action by $\La_K$.

The main tool that we shall use is the Seidel representation of 
$\pi_1(\Ham)$ in the group of automorphisms of the quantum 
cohomology ring $QH^*(M)$ of the manifold $(M,\om)$.  
This is very difficult to calculate in general;
some of the reasons for this are explained in
the examples in \S\ref{ss:ex}.
However, we make some progress in the case that
the fixed point set $M^{S^1}$ has a {\bf semifree} component $F$, 
i.e. a component $F$ so that the action is semifree on some neighborhood
of $F$.
(Recall that a circle action is semifree if the stabilizer of every
point is trivial or the whole circle.)
If the action is semifree on the whole manifold,
we calculate in Theorem~\ref{thm:semifree} the leading order 
term of the Seidel automorphism on quantum homology.  This provides 
the technical basis for Gonzalez's proof~\cite{Gonz} that, if in addition all the fixed points are isolated, then the manifold $(M,\om)$ has the same quantum cohomology as a product of $2$-spheres.

Another interesting special case is when the manifold is toric
 and the maximal fixed component of $\La_K$ corresponds to one of the facets of the moment polytope.   In this case, our results throw light on the multiplicative relations in the small quantum cohomology ring of $M$ for arbitrary toric manifolds, though we can calculate them  only in the Fano case: see \S\ref{ss:toric}.

We shall first state our results on the Hamiltonian group and then 
discuss properties of the Seidel representation.

\subsection{Results on the Hamiltonian group}\labell{ss:Ham}

 We get most information when one of the components 
on which $K$ is a maximum or minimum is semifree.

\begin{thm}\label{thm:simplemax}  Suppose that the Hamiltonian circle $\La_K$   
has a semifree maximal or minimal fixed point component.  Then 
$\La_K$ is essential in $\Ham(M,\om)$.
\end{thm}

This generalizes the result of McDuff--Slimowitz~\cite{MSlim}
stating that $\La_K$ is essential 
if the action is semifree.
It follows from 
Theorem~\ref{thm:max} (i) below, which 
calculates the leading order term of the Seidel element in the case
when the  maximal fixed point component is semifree. Because we are dealing with the maximal component, 
the proof is elementary; and the result itself  is
well known to experts, though as far as we 
know it is not formally stated in the literature.  
In contrast our later results are new.
Moreover their proofs are considerably harder.  Because most $S^1$-manifolds do not admit invariant, $\om$-compatible and semipositive almost complex structures $J$ (i.e. they are not symplectically NEF), it is important to work with general symplectic manifolds and hence with the virtual moduli cycle.  In 
\S\ref{ss:local} we prove new results about 
localization for Gromov--Witten invariants, that are familiar 
in the algebraic context but not in  
the symplectic case.

The proof of Theorem~\ref{thm:simplemax}  extends to cases where the 
isotropy of an extremal fixed component $F$ is small compared 
to the \lq\lq size" of $F$.
We can also deal with nonextremal fixed components provided that
the points above $F$ have sufficiently small isotropy.

To be precise,  we need a few definitions.
We say that a subset $N \subset M$ has {\bf at most $k$-fold isotropy}
if the stablizer of every point in $N$ has at most $k$ components. 
For each fixed point $x$ we 
denote by $m(x)$ 
the sum of the weights at $x$.
Finally,
given a fixed component $F$,
after
choosing an $\om$-compatible $S^1$-invariant almost complex structure 
$J$, decompose the negative normal bundle of $F$
as a sum of complex vector
bundles $E_1\oplus\dots\oplus E_\ell$
with weights $-k_1,\ldots,-k_\ell$, where $k_i\ge 1$.
The  associated {\bf obstruction bundle} is the bundle  
\begin{equation}\label{eq:ee}
\Ee: = (E_1\otimes \C^{k_1 - 1})\;\oplus\;\dots\;\oplus\; 
(E_\ell\otimes \C^{k_\ell - 1}).  
\end{equation}
Note that summands $E_i$ with $k_i=1$ do not appear in $\Ee$
since $E_i\otimes \C^{k_i - 1} = \{0\}$ in this case.
We say that $F$ is {\bf homologically visible} if 
the positive weights along $F$ are all $+1$, and
the associated obstruction bundle has nonzero Euler class.
We denote by
$e(\Ee)\in H_*(F)$ the Poincar\'e dual of this Euler class.
In particular, 
if all the negative weights are $-1$,
$e(\Ee) = [F]\in H_{\dim F}(F)$.

\begin{thm}\labell{thm:asimplemax}
Consider a Hamiltonian circle action $\Lambda_K$ on
a compact symplectic manifold $(M,\omega)$  with moment map $K$.
Assume that  $\Lambda_K$ is inessential in $\Ham(M, \om)$. 
Let $F$ be a homologically visible fixed component.
\smallskip

\NI{\rm (i)}
$F$ cannot be the maximal fixed component.
\smallskip

\NI{\rm (ii)}
More generally, if  every point 
in the 
superlevel set $\{K(x) > K(F)\}$
has at most twofold isotropy, then  $K(F) = m(F) = 0$
and $F$ is semifree.
\end{thm}

Even if the conditions above are not satisfied,
the existence of a semifree fixed point component still gives some information, though this is hard to interpret unless one has 
some global information about the action.  The reason 
is that, although semifree fixed components $F$ with nonzero 
$K(F)$ or $m(F)$ always  make nontrivial contributions to the
Seidel element, these could in general be cancelled by something else.
One special case is when the circle contracts in a 
compact Lie subgroup $G$ of the diffeomorphism group.
In this case, we can apply  Theorem~1.3 in McDuff--Tolman~\cite{McTLie} which
states that any semifree fixed component $F$ has a {\bf reversor}, that is there is
$g \in G$  which fixes $F$ but reverses
$\Lambda$ in the sense that 
$g^{-1} \La g = \La^{-1}$.  In the Hamiltonian case, a 
symplectomorphism $g$ is a reversor if and only if $K\circ g = -K$.  
The following result is an immediate consequence.

\begin{prop}\label{prop:Lie}
Consider a Hamiltonian circle action $\Lambda_K$  on
a compact connected
symplectic manifold $(M,\omega)$  with
moment map $K$. 
Let $F$ be a semifree component of the fixed point set.
Let $G \subset \Ham(M,\omega)$ be a compact Lie group
which contains  $\Lambda_K$.
If $\Lambda_K$ is inessential in  $G$, then there is an element $g\in G$ such that $g(F) = F$ and $K\circ g = -K$.  In particular,
$K(F) = m(F) = 0$.
\end{prop}

Another special case is when
the underlying manifold is toric and the circle $\La_K$ 
is a subgroup of the $n$-torus $T$ that acts on $M$ (where $\dim M = 2n$).
We denote by  $\Phi: M \to \ft^*$ the mean normalized moment map;
the moment image is the polytope  $\Delta = \Phi(M) \subset \ft^*$.
Given a face $f \subset \Delta$ of dimension $k$, the preimage $\Phi\inv(f)$ 
represents a homology class on $M$ of degree $2k$.  If $\La$ is any 
circle in $T$, the moment map $K$ for $\La$ is the composite of $\Phi$
with the associated  projection.

\begin{prop}\labell{prop:torsim}  Fix a  symplectic  toric
manifold $(M, \om,\Phi)$  with moment image $\Delta$. 
Consider a circle  subgroup $\La_K \subset T$ which is inessential in 
$\Ham(M,\omega)$.  Let $F$ be a semifree fixed point component 
for this circle action.
Then  $K(F) = m(F) = 0$.
Moreover, let $f^+$ and $f^-$ be the largest faces of $\Delta$ whose minimum
and maximum, respectively, are $\Phi(F)$.
Then,  $[\Phi\inv(f^+)] = [\Phi\inv(f^-)] \in H_*(M)$.
\end{prop}

This result does go further than Proposition~\ref{prop:Lie} 
since there are toric manifolds $M$ that admit circles 
$\La_K$ which are inessential in $\Ham(M,\om)$ but are 
essential in the natural compact Lie groups that act on 
$M$: see McDuff--Tolman~\cite{McTtor}.

Another tractable situation is when the action $\La_K$ 
has at most  twofold isotropy.  
In this case, even if there is no semifree fixed point component,  
the homology of the manifold $M$ has the symmetry properties 
that one would expect if the action had a reversor $\phi$.  
(Again, this echoes a result in  McDuff--Tolman~\cite{McTLie} 
which states that if a Hamiltonian circle action is
inessential in a compact Lie subgroup $G$ of $\Ham(M,\om)$,
then it can be reversed.)
Here is a precise statement.

\begin{prop}\labell{prop:2iso0} 
Consider a Hamiltonian circle action $\Lambda_K$ on a compact
symplectic manifold $(M,\omega)$ with moment map $K$.
Assume that $\Lambda_K$ is inessential in $\Ham(M,\omega)$
and that $\Lambda_K$ has at most twofold isotropy.    
\smallskip

\NI{\rm (i)}
For all $\mu\in \R$, a homology class 
of $M$
can be represented in the sublevel set 
$\{K(x)<\mu\}$
if and only if it can be represented in the superlevel set 
$\{K(x)>-\mu\}$.

\NI{\rm (ii)} For any connected component $N \subset M^{\Z/(2)}$,  any
integers $j$ and $n$, and any $\mu \in \R$,
$$
 \bigoplus_{\stackrel{\scriptstyle{F \subset N}}{K(F) = \mu,\ m(F) = n}}
H_{j - \alpha_F}(F) 
\;\; \cong \;\;
\bigoplus_{\stackrel{\scriptstyle{F' \subset N,}}{K(F') = -\mu,\ 
 m(F') = -n}}
 H_{j - \beta_{F'}}(F'), $$ 
where the sums are over fixed components,
$\alpha_F$ is the Morse index of $F$ with respect to $K$,
and $\beta_{F'}$ is the Morse index of $F'$ with respect to $-K$.
\smallskip

\end{prop}

Let $F_{\max}$ and $F_{\min}$ denote the maximal and minimal 
components of $K$, and define $K_{\max} = K(F_{\max})$, $K_{\min} =
K(F_{\min})$, $m_{\max} = m(F_{\max})$, and $m_{\min} = m(F_{\min})$.
The proposition above  immediately implies that  $M$ 
is symmetrical with respect
to $K$, in the sense that  $K_{\max }= - K_{\min}$,
$m_{\max} = - m_{\min}$, and $H_i(F_{\max}) = H_i(F_{\min})$ for
all $i$.
Similarly, every component $N$ of $M^{\Z/(2)}$ is symmetrical
with respect to $K$ in the same sense.

More generally one can show that the moment image is not too skew
whenever 
the isotropy weights are not too large.  
Given two fixed point components $A$ and $B$, let $q = q(A,B)$ denote the
largest integer so that $A$ and $B$ lie in the same component of 
$M^{\Z/(q)}$.

\begin{prop}\labell{prop:Kdist}  
Consider a Hamiltonian circle action $\Lambda_K$ on a compact
symplectic manifold $(M,\omega)$ with moment map $K$.
Assume that $\La_K$ is inessential in $\Ham(M,\om)$.
Then there exists some sequence
of fixed components   $F_{\max} = F_0, F_1, \ldots, F_j = F_{\min}$
with $K(F_i) \neq K(F_{i-1})$ for all $i$
so that
$$ K_{\max}  \geq   \sum_{i= 1}^{j}  \frac
{|K(F_{i-1}) - K(F_i)|}
{q(F_{i-1},F_i)}.$$
Moreover,  we can choose the sequence so  either 
the inequality above is strict, or
$$ m_{\max}    = \sum_{i = 1}^{j}  \frac
{(m(F_{i-1}) - m(F_i))}
{q(F_{i-1},F_i)} \cdot
\frac {K(F_{i-1}) - K(F_i)} {|K(F_{i-1}) - K(F_i)|}
.$$
\end{prop}

Note that if $M$ has at most $k$-fold isotropy, then
for any sequence of fixed point components 
  $F_{\max} = F_0, F_1, \ldots, F_j = F_ {\min}$ with
$K(F_{i-1}) \neq K(F_i)$ for all $i$,
$  \sum  \frac
{|K(F_{i-1}) - K(F_i)|}
{q(F_{i-1},F_i)} \geq \frac{K_{\max} - K_{\min}}{k}$.  
Moreover, the inequality is strict unless 
 $F_{\max}$ and $F_{\min}$ (and in fact all the $F_i$) lie  
in the same component of $M^{\Z/(k)}$,
and $K(F_{i-1}) >  K(F_{i})$ for all $i$.
In this case,
$ \sum \frac
{(m(F_{i-1}) - m(F_i))}
{q(F_{i-1},F_i)} \cdot
\frac {K(F_{i-1}) - K(F_i)} {|K(F_{i-1}) - K(F_i)|}
= \frac{m_{\max} - m_{\min}}{k}$.
Therefore, the proposition has the following corollary.

\begin{corollary}
Suppose  in the situation of Proposition~\ref{prop:Kdist} that
 $M$ has at most $k$-fold isotropy.
Then, after possibly reversing the circle action,
$$
K_{\max} \leq |K_{\min}| \leq (k-1) K_{\max},
$$ 
and the 
second
inequality is strict unless
$m_{\min} = (k-1) |m_{\max}|,$
and
$F_{\max}$ and $F_{\min}$ lie in the same
component of $M^{\Z/(k)}$. 
\end{corollary}

Observe finally that the inequality in Proposition~\ref{prop:Kdist} is 
sharp:  just consider the action on $\CP^{k}$ given by
$[z_0:z_1:z_2: \cdots :z_k]\mapsto 
[\lambda^{k} z_0:\lambda z_1:z_2:\cdots: z_k]$.

Similar but more complicated statements can be 
made about the other fixed point components: 
see Prop~\ref{prop:xdist}.  

\subsection{Results on the Seidel representation}\labell{ss:Srep}
  
We now state our main results on the Seidel
representation and use them to deduce 
Theorems~\ref{thm:simplemax} and \ref{thm:asimplemax}, and
Propositions~\ref{prop:torsim} and
\ref{prop:2iso0}.

The Seidel representation is a quantum 
version of a  classical homomorphism defined 
by Weinstein~\cite{Wei} using the action functional.
Let $H_*^S(M): = H_*^S(M;\Z)$ denote the spherical
homology of $M$.
Let   $I_\om, I_c:H_2^S(M)\to \R$ denote the homomorphisms induced by 
evaluating the classes $[\om]$ and  $c_1 = c_1(TM)\in H^2(M,\Z)$.
Weinstein's  homomorphism 
$
\Aa_\om: \pi_1(\Ham(M)) \to \R/(\im\, I_\om)
$
takes the circle $\La_K$ to the value $K(x)$ of the generating 
moment map at any 
critical point.\footnote{
This homomorphism itself contains quite a bit of information: see
for example McDuff--Tolman~\cite{McTtor}.}   
As we show in~\S\ref{ss:se} this extends
to take the weights into account. 

\begin{lemma}\labell{le:actind} Let $(M,\omega)$ be a compact
symplectic manifold. There is a homomorphism 
$$
\Aa_{\om,c}:   \pi_1(\Ham(M,\om)) \to \R\oplus \Z/\im (I_\om 
\oplus I_c)
$$
whose value at a Hamiltonian circle action $\La_K$ is $[K(x), -m(x)]$, 
where $x$ is any critical point of $K$.
\end{lemma}

The Seidel representation 
$$
\Ss:\pi_1(\Ham(M, \om)) \to \QH_\ev(M;\La)^\times
$$
is a lift of $\Aa_{\om,c}$
to the group of even units $\QH_\ev(M;\La)^\times$ 
of the quantum homology ring $\QH_*(M) = \QH_*(M;\La): = H_*(M)\otimes \La$ of
$M$: see~\cite{Sei,LMP,Mcq}.  Here, following~\cite{MS2},
we use coefficients $\La: = 
\La^{\rm univ}[q,q^{-1}]$ where $q$ is a variable of degree $2$ and
$\La^{\rm univ}$ is a generalized Laurent series ring in a 
variable $t$ of degree $0$:
$$
\La^{\rm univ}: = \Bigl\{\sum_{\ka\in \R} r_\ka t^\ka\; \Big| \; r_{\ka}\in \Q,\;\;
\#\{\ka > c\,|\,r_\ka \ne 0\}< \infty,\forall \; c\in \R\Bigr\}.
$$
We shall order the elements  $\sum_{d,\ka} 
a_{d,\ka} \otimes q^d t^{\ka}$ in $\QH_*(M;\La)$ by the valuation\footnote
{
One must treat this ordering with some care.  Although
$v(a*b) \le v(a) + v(b)$
for all $a,b\in \QH_*(M)$ with equality only
if the usual intersection product of the highest order terms
is nonzero, in the case when this intersection product is zero
the term of highest order in $a*b$ may not
be equal to the product of the highest order terms in $a$ and $b$.
}
$$
v\Bigl(\sum_{d,\ka} a_{d,\ka} \otimes q^d t^{\ka}\Bigr) = 
\max \{\ka\;|\;\exists \, d :  a_{d,\ka}\ne 0\}.
$$
For more details see \S\ref{ss:se}.

The image $\Ss(\La)\in \QH_\ev(M;\La)^\times$ of the Hamiltonian  loop $\La$  
is called the {\bf Seidel element} of $\La$.  
It has degree $\dim M$ and gives rise to
a degree preserving automorphism of $\QH_*(M)$ by quantum multiplication:
$$
\Ss(\La)(a): = \Ss(\La) * a.
$$
Thus $\Ss(\La) = \Ss(\La)(\1)$ where 
$\1$ denotes the unit $[M]$ in
$\QH_*(M)$.
Our results are based on a partial calculation of 
$\Ss(\La_K)$.

If a circle action $\La$ is inessential, then $\Ss(\La) = \1$, so
$\Ss(\La)(a) = a$ for all homology classes $a$.
Therefore, if we write
$$ \Ss(\La)(a) = \sum_{d,\kappa} c_{d,\kappa}  \otimes q^d t^\kappa,$$
there are potentially two techniques to show
that $\La$ is essential.
First,  we can show that $c_{d,\kappa} \neq 0$ for some pair $(d,\kappa)
\neq (0,0)$;  this is the approach we take for homologically visible
fixed points.
Second, we can find a homology class
$b$ so that  $a \cdot b \neq 0$ but
that $c_{0,0} \cdot b = 0$; this is the approach that 
we take in all other cases.

The following theorem is proved in Section~ \ref{ss:chain}.
It is a considerable generalization of
Proposition 7.11 in Seidel~\cite{Sei}, which applies only in the case of monotone $(M,\om)$ and does not discuss the structure of the higher order terms. 
By analogy with the 
algebraic case, we say that the almost complex manifold $(M,J)$  is 
Fano (resp. NEF) if there are no $J$-holomorphic spheres in classes 
$B$ with $c_1(B) \le 0$ (resp. $c_1(B) < 0$).

\begin{thm}\labell{thm:max}  
Consider a Hamiltonian circle action $\La_K$  on a
compact symplectic manifold $(M,\omega)$ with normalized moment map $K$.
Assume that the maximal fixed point component $\Fmax$ is semifree.
Then:
\smallskip

\NI{\rm (i)}
$$
    \Ss(\La_K) = [\Fmax]\otimes q^{-m_{\max}}\, t^{K_{\max}} +  
     \sum_{B \in H_2^S(M): \omega(B) > 0} a_B \otimes q^{-m_{\max}-c_1(B)}\, 
     t^{K_{\max}-\om(B)},
$$ 
where $a_B$  is the contribution from the section class $\si_{\max} + B$.
Moreover, if $a_B \neq 0$ then $\deg(a_B)  =  \dim F_{\max} + 2c_1(B)$.
\smallskip

\NI{\rm (ii)}
If $(M,J)$ is Fano (resp. NEF) for  some $S^1$-invariant 
$\om$-compatible 
almost complex structure
$J$  then   $a_B = 0$ unless  $c_1(B) > 0$ (resp. 
$c_1(B) \geq 0$).
\smallskip

\NI{\rm (iii)}
Assume that $(M,J)$ is NEF  for some $S^1$-invariant $\om$-compatible 
almost complex structure $J$.   
If $2c_1(B') \geq \codim \Fmax$
for every $J$-holomorphic sphere $B'$ which intersects 
$\Fmax$ then all the lower order terms vanish.  If 
the latter
hypothesis holds
except for spheres which lie in $\Fmax$ itself,  
then  $a_B = 0$ unless $2c_1(B) < \codim \Fmax$ and $B$  lies in the image of
the spherical homology $H^S_2(\Fmax)$.\end{thm}

The last sentence in part (i)  of Theorem \ref{thm:max} expresses the fact
that $\Ss(\La_K)$ preserves degree. 
If  $a_B \neq 0 $ then 
$$
\deg\bigl(a_B \otimes q^{-m_{\max}-c_1(B)}\,
t^{K_{\max}-\om(B)}\bigr)  = \deg a_B  - 2m(F_{\max})  - 2c_1(B) = \dim M.
$$
Moreover,  $\dim F_{\max} = \dim M + 2m(F_{\max}) $ because $F_{\max}$ is semifree.
Therefore 
$$
\deg a_B  = \dim F_{\max} + 2c_1(B).
$$

This theorem gives the most information when 
$\codim \Fmax = 2$, for example in the case of a circle 
action on a toric variety that
fixes one facet.  
If, in addition $(M,\omega,J)$ is Fano for some 
$S^1$-invariant $\om$-compatible $J$,  then, by part (ii)
all the lower order terms vanish.
That is, $\Ss(\Lambda_K)= [\Fmax] \otimes q^{-m_{\max}} t^{K_{\max}}$.
In \S\ref{ss:chain} we shall give a more precise description of 
the lower order terms in $\Ss(\La)$.  
These remarks have consequences for the structure of the quantum cohomology 
of toric manifolds that are explained in \S\ref{ss:toric}.
\MS

\begin{example}\label{ex:s2}\rm
Consider the rotation of $S^2$ generated by the
height function $K$ and let $A= [S^2]$.
Then $\Ss(\La_K) = [pt]\otimes q t^{\om(A)/2}$.   
\end{example}

\NI
{\bf Proof of Theorem \ref{thm:simplemax}.}
If $\Lambda_K$ is inessential in $\Ham(M,\omega)$,  
then $\Ss(\Lambda_K) = \1$.
Hence Theorem \ref{thm:simplemax} follows immediately from  
the first claim of Theorem \ref{thm:max}.\QED \MS   

\NI
{\bf Proof of Proposition \ref{prop:torsim}.}
Pick any $x \in F$.
There is a neighborhood of $x$ and an isomorphism
of $T$ with $(S^1)^n$ so that
the action of $T$ is equivariantly symplectomorphic
to the standard action of $(S^1)^n$ on $\C^n$.
In these coordinates, the action of $\Lambda_K$ on $\C^n$ is given by
$ \lambda z = (\lambda^{m_1} z_1,\ldots,\lambda^{m_n} z_n)$, where 
$m_1,\dots,m_n$ are the weights at $x$.

For $ 1 \leq i \leq n$, let $D_i$ be the  facet of $\Delta$ which corresponds 
to $z_i = 0$.  Let $\eta_i \in \ell$ denote the outward primitive
normal vector to $D_i$, where $\ell \subset \ft$ is the integral lattice. 
Note that $K_i: = \langle\eta_i,\Phi(\cdot)\rangle$ 
is the moment map for a circle action $\Lambda_i$, and that $\Phi\inv(D_i)$
is a semifree maximum 
for this action.  By Theorem~\ref{thm:max},
$\Ss(\Lambda_i) = y_i\otimes q\,t^{\eta_i(D_i)}$,
where $y_i = [\Phi\inv(D_i)] + $ lower order terms.
Since $\Lambda_K$ is inessential, $\Ss(\Lambda_K) = \1$.
On the other hand, by looking at the action near 
the fixed point $x$ one sees that
$\Lambda_K = \prod \Lambda_i^{-m_i}$.  Therefore
$$
\prod\; [y_i]^{-m_i}\otimes q^{m_i}\, t^{-m_i \eta_i(D_i)}  = \1,
$$
where the product is taken in $\QH_*(M;\La)$.

Let $m_i = 1$ for $1 \leq i \leq r$, $m_i = -1$
for $r < i \leq r + s$, and $m_i=0, i>r+s$.
Then
$$
y_{1} * \cdots * y_{r}\otimes q^{-r}\;t^{\eta_1(D_1) +  \cdots +  \eta_r(D_r)} 
= y_{r+1} * \cdots * y_{r+s} \otimes q^{-s}\; t^{\eta_{r+1}(D_{r+1}) +  
\cdots +  \eta_{r+s}(D_{r+s})} 
$$
In particular, the highest order terms must agree. 
Since $D_1 \cap \cdots \cap D_r = f^+$ and
$D_{1+r} \cap \cdots \cap D_{r+s} = f^-$,
we have $[\Phi\inv(D_1)] \cap \cdots \cap [\Phi\inv(D_r)] = [\Phi\inv(f^+)]$
and
$[\Phi\inv(D_{1+r}) ] \cap \cdots \cap [\Phi\inv(D_{r+ s})] = [\Phi\inv(f^-)]$.
Since these intersections are nontrivial,
the highest order terms in the quantum product of the corresponding 
$y_i$ are given by these intersections. 
The result follows.
\QED

\begin{remark}\labell{rmk:torsim}\rm
More generally, let $F$ be a fixed component of any inessential Hamiltonian
loop $\Lambda_K \subset T$.  Let $(m_1,\ldots,m_n)$ be the weights at 
$x \in F$ with corresponding facets $D_i$.
Define homology classes in $M$ by $ X^+ = \cap_{m_i > 0} 
 [\Phi\inv(D_i)] ^{m_i}$
and $X^- = \cap_{m_i < 0}  [\Phi\inv(D_i)]^{-m_i}.$
(Here, we are taking the ordinary cap product in homology.)
If both $X^+$ and $ X^-$ are nonzero, 
then, by an argument similar to the one above,
$K(F) = m(F) = 0$ and $X^+ = X^-$.
\end{remark}

\subsubsection{Semifree actions and canonical bases for 
homology}\label{sss:canon}

For a general action our methods do not give 
any information about $\Ss(\La_K)(a) : = \Ss(\La_K)*a$ for $a\ne \1$.  
However, when the action is semifree, it is possible to  
describe the top order term in $\Ss(\La_K)(a)$ for any $a\in H_*(M)$.   
This formula is best written in terms of some canonical bases $\{c_i^-\}$ 
and $\{c_i^+\}$  for $H_*(M)$.

Before explaining this, we introduce more notation.
Given $\mu\in \R$, define 
$$
\begin{array}{ll}
M_\mu: = K\inv\bigl([\mu,\infty)\bigr),&
M_{> \mu}: = K\inv\bigl(\mu,\infty)\bigr),\\
M^\mu := K\inv\bigl((-\infty, \mu]\bigr),&
M^{< \mu} := K\inv\bigl((-\infty, \mu)\bigr).
\end{array}
$$
The inclusions $M_\mu \to M$ and $M^\mu \to M$
induce maps $H_*(M_\mu) \to H_*(M)$ and $H_*(M^\mu) \to H_*(M)$ in
rational homology.  
We call the images of these maps $F_\mu H_*(M)$ and $F^\mu H_*(M)$, 
respectively. 

We now give a brief review of equivariant cohomology:
Let $S^1$ act on a space $N$.  The equivariant cohomology $H^*_{S^1}(N)$ of $N$ is defined to be the 
cohomology of the total space 
$$
N_{S^1} : = S^{\infty}\times _{S^1} N
$$  
of the universal $N$-bundle over the 
classifying space $BS^1 = \CP^{\infty}$.  Thus $H^*_{S^1}(N)$  is a 
module over $H^*_{S^1}(pt)\cong H^*(\CP^{\infty})$, 
which  is 
a polynomial ring with one 
generator $u$ of degree $2$.  Moreover, there is a natural map
from $H_{S^1}^*(N)$ to $H^*(N)$, given by restricting to any fiber.

If $S^1$ acts trivially on $F$  there is a natural identification 
$H_{S^1}^*(F) = H_{S^1}^*(pt) \otimes H^*(F)$.
Given $\TY \in H_{S^1}^*(F)$, we say that the 
{\bf degree} of $\TY$ in $H_{S^1}^*$ is 
$j$ if $j$ is the smallest integer such that
$$
\TY \in \bigoplus_{i=0}^jH^i_{S^1}(pt) \otimes H^*(F).
$$
We now explain a procedure for producing a natural
set of generators for the equivariant cohomology of $M$
given a set of generators for the cohomology of each fixed component.

\begin{lemma}\label{le:unique0}
Let $S^1$ act on a compact symplectic manifold $(M,\om)$ with moment map $K$.
Let $F \subset M$ be any fixed component of index $\alpha$;  
let $e^-_F\in H_{S^1}^\alpha(F)$ be the equivariant 
Euler class of the negative normal bundle to $F$.
Given any cohomology class $Y \in H^i(F)$, there exists a unique 
cohomology class $\TY^+ \in H_{S^1}^{i + \alpha}(M)$ so that
\smallskip

\NI{\rm (a)}
 The restriction of $\TY^+$ to $M^{< K(F)}$ vanishes, 
\smallskip

\NI{\rm (b)}
 $\TY^+|_F = Y \cup e^-_F$,  and 
\smallskip

\NI{\rm (c)}
 the degree of $\TY^+|_{F'}$ in $H_{S^1}^*(pt)$ is less than the 
index $\alpha_{F'}$  of $F'$ for all fixed components $F' \neq F$.
\smallskip

\NI
Moreover, these classes generate $H^*_{S^1}(M)$ as a $H^*_{S^1}(pt)$ module.
\end{lemma}

We can use this lemma, which we prove in section \ref{ss:canclasses}, 
to create a set of generators for the homology of $M$.
Let $F$ be a fixed component, and let
$\alpha_F$ and $\beta_F$ denote the index of $F$ with respect to
$K$ and $-K$, respectively.
Given a homology class  $c \in H_i(F)$, we define
the {\bf upwards extension}  $c^+ \in H_{i + \beta_F}(M)$ as follows:
\begin{itemize}
\item Let $Y \in H^{\dim F - i}(F)$ be the Poincar\'e dual to $c$. 
\item Let  $\TY^+ \in H_{S^1}^{\dim F + \alpha_F - i}(M)$ be
the unique equivariant cohomology class which satisfies the
conditions of Lemma~\ref{le:unique0}.
\item Let ${Y}^+ \in H^{\dim F + \alpha_F - i}(M)$ be the restriction
of $\TY^+$ to ordinary cohomology.  
\item Let $c^+ \in H_{i + \dim M - \dim F - \alpha_F}(M) =
H_{i + \beta_F}(M)$ be the Poincar\' e dual to ${Y}^+$.
\end{itemize}
Note that, by construction, $c^+$ lies in $F_{K(F)} H_*(M)$.
The {\bf downwards extension} $c^- \in  H_{i + \alpha_F}(M)$,
which lies in $F^{K(F)} H_*(M)$ is defined analogously; 
simply replace $K$ by $-K$.

Since the classes $\TY^+$ generate $H^*_{S^1}(M)$ as 
a $H^*_{S^1}(pt)$ module and the restriction $H^*_{S^1}(M) \to H^*(M)$
is surjective, the classes ${Y}^+$ generate
$H^*(M)$ as  a 
(rational) 
vector space. Hence, the classes $c^+$
(or, alternatively, the classes $c^-$) generate $H_*(M)$ as a vector space.

When the action is semifree, the classes $c^+$ and $c^-$ have a nice
geometric description.
Assume that $c$ can be represented by an $i$-dimensional 
submanifold $C \subset F$.
By Lemma~\ref{sfpseudo},
if $g_J$ is the metric associated to a generic
$S^1$-invariant
$\om$-compatible almost complex structure $J$
and we choose $C$ generically,
the stable manifold $W^s(C)$ is an $(i + \beta_F)$-dimensional 
pseudocycle.
(See section \ref{pseudosec}.)
Hence, it represents a homology class $[W^s(C)] \in  H_{i + \beta_F}(M)$.
By Proposition~\ref{prop:equivar}, $[W^s(C)] = c^+$.
Similarly, $c^-$ is represented by the unstable manifold $[W^u(C)]$.

\begin{rmk}\rm 
We may define an automorphism $\Dd_K:H_*(M)\to H_*(M)$ by  
$$
\Dd_K (c^-) = c^+, \quad c\in H_*(M^{S^1}).
$$
For example,  if $c=[F_{\max}]\in H_*(F_{\max})$ 
is the maximal fixed point set of $K$, 
then $c^- = \1$ while $c^+ = [\Fmax]$.  Therefore
$$
\Dd_K(\1) = [F_{\max}].
$$
If $K$ is 
Morse then $\Dd_K$ can be interpreted as a form of  
duality.  If $\{c_i\}$ is given by  
the set of critical points of $K$, then the bases
$\{c_i^-\}$ and $\{c_i^+\} = \{\Dd_K(c_i^-)\}$ are 
dual with respect to the intersection pairing.
Although it is tempting to think that $\Dd_K$ is an involution, in 
fact the correct relation is $\Dd_{-K}\circ\Dd_K = \1$.\QED 
\end{rmk}

The following theorem is proved in section \ref{ss:semifree}.

\begin{thm}\labell{thm:semifree}  Let $S^1$ act semifreely on 
a compact symplectic manifold $(M, \om)$.  Let $F$ be a component of the fixed
point set, and choose a homology class $c \in H_*(F)$.  
Then
\begin{eqnarray}\labell{eq:simp}
\Ss(\La_K)(c^-) &=& c^+ \otimes q^{-m(F)}\,\,t^{K(F)} + \\\notag
&&\qquad \sum_{B\in 
H_2^S(M):
\om(B)> 0} a_B\otimes q^{-m(F)-c_1(B)}\,\,t^{K(F)-\om(B)}.
\end{eqnarray}
Moreover if $a_B \neq 0$ then $\deg a_B = \deg c^+ + 2c_1(B)$.
\end{thm}

Since every element $a \in H_*(M)$ can
be written  as a linear combination of such $c^-$,
this theorem gives the leading order term of $\Ss(\La_K)(a)$
for every $a \in  H_*(M)$.

The last claim of the theorem follows 
from the fact that $\Ss$ preserves degree.
This implies that if $a_B \neq 0$, then 
$$
\deg(a_B) - 2m(F) - 2 c_1(B) = \deg (c^-).
$$ 
Since the action is semifree, $m(F)$ is the number of positive weights
minus the number of negative weights.  But
the degree of $c^+$ is the degree of $c$ plus twice the number of positive
weights, and the degree of $c^-$ is the degree of $c$ plus twice the number
of negative weights.
Hence $\deg(a_B)  = \deg (c^+) + 2c_1(B)$.

\begin{example}\rm Think of $\CP^2$ as
the manifold obtained from the closed unit ball in $\C^2$ by
identifying its boundary to a complex
line via the Hopf map, and consider the action
$$
(z_1,z_2)\mapsto (\lambda^{-1} z_1, \lambda^{-1} z_2).
$$
Then $K(z_1,z_2) = \pi (c-|z_1|^2 - |z_2|^2)$ where $c = 2/3$,
$F_{\max} = \{pt\}$ and all critical points are semifree.  
Since $c_1(L) = 3$ where $L = [\CP^1]$,
there can be no lower order terms in the formula for
$\Ss(\La_K)(a)$ since the dimensional condition can never be 
satisfied.  Hence,
since $\om(L): = \pi$,  we find that $\Ss(\La_K)$ acts by:
$$
\1 \mapsto [pt]\otimes q^2\,t^{2\pi/3},\;\;
L \mapsto \1\otimes q^{-1}\,t^{-\pi/3},\;\;
[pt] \mapsto L\otimes q^{-1}\,t^{-\pi/3},
$$
which is consistent with the formula $\Ss(\La_K)(a) = \Ss(\La_K)*a$.
The  above results also agree with the formulas found in~\cite{Mch}\,\S4
for rotations of the one point blow up of $\CP^2$: see
Example~\ref{ex:CPB}.  In this example we shall also see that
although Theorem~\ref{thm:max} implies that there are no 
lower order terms in the Seidel element $\Ss(\La)$ itself if $(M, \om)$ is 
Fano and $F_{\max}$ has codimension $2$, there may be lower order terms
in $\Ss(\La)(a)$ for such actions.
\end{example}

\subsubsection{Actions with at most twofold isotropy}\label{sss:twofold}

We can also obtain some information about $\Ss(\Lambda_K) (a)$,
though considerably less than before,
when  $\Lambda_K$ acts
with at most twofold isotropy.  
Throughout the following discussion we denote by $\cdot_Y$ the 
intersection pairing $H_k(Y)\times H_{m-k}(Y)\to \Q$ on the homology 
of an oriented $m$-dimensional manifold $Y$.  For convenience we set
$a\cdot_Y b= 0$ whenever the dimensional condition $\deg(a) + 
\deg(b) = m$ is not satisfied.

The following theorem is proved in section \ref{ss:twofold}.

\begin{thm}\labell{thm:2iso}
Consider a Hamiltonian circle action $\La_K$ on a compact
symplectic manifold $(M,\omega)$ with at most twofold 
isotropy. 
Let $F$ and $F'$ be  components of the fixed point set. 
Choose  homology classes $c \in H_*(F)$ and $c' \in H_*(F')$,
and  write
$$
\Ss(\La_K)(c^-) = \sum_{d,\ka} c_{d,\ka}\otimes q^d\,\,t^\ka. 
$$
If $K(F') \leq -K(F)$, 
then $c_{0,0} \cdot_M (c')^- = 0$
unless $K(F) = -K(F')$, $m(F) = -m(F')$, and
$F$ and $F'$ lie in the same component of $M^{\Z/(2)}.$
\end{thm}

\NI
{\bf Proof of Proposition~\ref{prop:2iso0}.}
Let $F$ and $F'$ be  components of the fixed point set 
so that $K(F') \leq -K(F)$. 
Choose homology classes $c \in H_*(F)$ and $c' \in H_*(F')$.
Since $\Lambda_K$ is inessential,  
$\Ss(\Lambda_K)(c^-) = c^- \otimes \1$.
Hence, by Theorem~\ref{thm:2iso}
$c^- \cdot_M (c')^- = 0$ unless $K(F') = -K(F)$,
$m(F') = -m(F)$, and $F$ and $F'$ lie in the same component of $M^{\Z/(2)}$

For any $\mu \in \R$,
$F^{\mu}H_*(M)$ is generated by elements $c^-$, 
where $c \in H_*(F)$ and $F$ is a fixed component with $K(F) \leq \mu$.
The paragraph above implies 
that every such $c^-$ lies in $F_{-\mu} H_*(M)$.
Hence
$F^{\mu}H_*(M)\subseteq F_{-\mu}H_*(M)$.
Similarly, applying the theorem to the moment map $-K$,
$F^{\mu}H_*(M) \supseteq F_{-\mu}H_*(M)$.
Hence $F^\mu H_*(M) = F_{-\mu} H_*(M)$.
This proves the first claim.

Since  both  $K$ and $-K$ are perfect Morse functions,
both $H_*(M_\mu) \to H_*(M)$ and $H_*(M^\mu) \to H_*(M)$ are injections.
Hence,  the  arguments above imply that
$$
H_j(M^\mu, M^{<\mu})  = H_j(M_{-\mu}, M_{>-\mu}).
$$
By the Thom isomorphism theorem, this is equivalent to 
$$
\bigoplus_{K(F) =\mu} H_{j -  \alpha_F}(F) =
\bigoplus_{K(F')=-\mu} H_{j -  \beta_{F'}}(F'),
$$
where the sums are over fixed components,
$\alpha_F$ is the index of $F$ with respect to $K$,
and $\beta_{F'}$ is the index of $F'$ with respect to $-K$.

Now suppose that  $F$ and $F'$ 
lie in different components of the isotropy submanifold $M^{\Z/(2)}$
and satisfy $K(F) = -K(F') = \mu$.
Consider $c \in H_i(F)$ and $c' \in H_*(F')$.
We saw above that $c^- \cdot_M (c')^- = 0$.
Hence,
the $F'$ component of the image of $c$ under the isomorphism above must
be zero.
Hence, the isomorphism above  is still an isomorphism when
restricted to any component of $M^{\Z/(2)}$.
A similar argument applies if $m(F) \neq -m(F')$.
\QED

\MS

Finally, let us consider the contribution of a homologically visible 
component.  The following theorem is also proved in~\S\ref{ss:twofold}.

\begin{thm}\labell{thm:pfasimplemax}  Consider  a Hamiltonian
circle action $\Lambda_K$ on a compact symplectic manifold.
Let $F$ be a fixed component. 
Assume that all the positive weights at $F$ are $+1$.
Suppose further that  the superlevel set  $\{K(x) > K(F)\}$ 
has at most twofold isotropy.
Choose a homology class $c \in H_*(F)$, and write
$$
\Ss(\Lambda_K)(c^-) = \sum_{d,\ka}\; c_{d,\ka}\otimes 
q^{-m(F)+d}\,\, \,t^{K(F)+\ka}.
$$
Then $c_{0,0} \in F_{K(F)} H_*(M)$.
Moreover,  for any $c' \in H_*(F)$,
$$
c_{0,0} \cdot_M (c')^-  = (e(\Ee) \cap_F  c) \cdot_F c'
$$
where $e(\Ee)$ denotes the Poincar\'e dual of the
Euler class of the obstruction bundle $\Ee\to F$.
(See equation (\ref{eq:ee}).)
\end{thm}

\MS

\NI
{\bf Proof of Theorem~\ref{thm:asimplemax}.}
Let $F$ be a homologically visible fixed component
and assume that every point in the superlevel set
$\{K(x) > K(F)\}$ has at most twofold isotropy.
Apply Theorem~\ref{thm:pfasimplemax} with $c = \1_F \in H_*(F)$ and 
$c'\in H_*(F)$ chosen so that $e(\Ee)\cdot_F c' = k\ne 0$.  Then 
$c_{0,0} \cdot_M (c')^- = e(\Ee)\cdot_F c' = k \neq 0$. 
Therefore the coefficient 
$c_{0,0}$ of $q^{-m(F)}\,t^{K(F)}$ in
$\Ss(\Lambda_K)\bigl((\1_F)^-\bigr)$ is nonzero. 
Since $\Lambda_K$ is inessential,
$\Ss(\La_K)\bigl((\1_F)^-\bigr) = (\1_F)^-$.
Therefore  
$K(F) = m(F) = 0$ 
and  
$c_{0,0} = (\1_F)^-$, so  $e(\Ee) = (\1_F)$.
This proves (ii).
Item (i) is simply a special case.
\QED 
\MS

\NI
{\bf Acknowledgements}  The first author thanks
Paul Seidel and Yong-Geun Oh for useful conversations and
the Ellentuck Foundation and
the Institute for Advanced Study for their generous support during 
Spring 2002.

\tableofcontents

\section{Quantum homology and the Seidel representation}\labell{sec:qh}

This section reviews the necessary background material.
The main geometric idea  behind
our results, symplectic bundles over the two sphere,
is explained in \S\ref{ss:sb}.
 We  review (small) quantum homology  in \S\ref{ss:qh}
to fix notational
conventions, and then describe the Seidel representation in \S\ref{ss:se}. 

\subsection{Symplectic bundles over the two sphere}\label{ss:sb}

Throughout we shall use the following notational/sign conventions.
If $ H_{t}, 0\le t \le 1,$ is a (time dependent) Hamiltonian
then we define the corresponding vector field $X_H$  by
the identity
\begin{equation}\labell{eq:XH}
\om(X_H,\cdot) = -dH_t.
\end{equation}
Thus $X_H = J (\grad H_t)$, where $J$ is an $\om$-compatible almost
complex structure and the gradient is taken with respect to the metric
$g_J$ given by $g_J(x,y) = \om(x,Jy)$.  As an example,
consider the unit sphere $S^2$ in $\R^3$, oriented via stereographic
projection from the north pole.\footnote
{
This means that the vertical
projection from the tangent space at the south pole to the $(x_1,x_2)$-plane
preserves orientation.  Hence this orientation is the opposite of
its orientation as the boundary of the unit ball.}
Then its area form is $dx_3\wedge d\theta$ and the vector field 
$X_K$ generated by
the normalized height function $K = 2\pi x_3$  is
$X_K = 2\pi \p_\theta$.  Thus the corresponding flow
is the anticlockwise rotation of $S^2$
about the axis from the south to the north pole.
 Note that this flow is positive (i.e.
anticlockwise) at the
south pole $s$ (the minimum of $K$) and negative at the north pole $n$ (the
maximum of $K$), which agrees with the usual conventions for
defining the {\bf moment map}.

Consider the locally trivial
bundle $P_\La\to S^2$ constructed by using $\La = \{\phi_t\}\in
\pi_1(\Ham(M))$ as a clutching function:
\begin{equation}\label{eq:Pla}
P_\La = (D_0\times M)\; \cup \;(D_\infty\times M)/\sim,\quad
\mbox{where }\; (e^{2\pi it},\phi_t(x))_0 \sim
(e^{2\pi it},x)_\infty.
\end{equation}
Here we are thinking of $D_0$ as the closed unit disc centered at $0$
in the Riemann sphere
$S^2 = \C \cup\{\infty\}$ and of $D_\infty$
 as another copy of this disc, embedded in $S^2 = \C\cup \{\infty\}$
via the orientation reversing map $r\, e^{i\theta}\mapsto
r^{-1}\, e^{i\theta}$.
Correspondingly, we denote the fibers over $0,
\infty$ by $M_0, M_\infty$. Note that
our definition of $P_\La$ agrees with that in~\cite{Mch} but
differs in orientation from
the convention used in~\cite{LMP,Mcq}.

The fact that $\La$ is Hamiltonian implies that there is
a closed $2$-form $\Om$ on $P_\La$ extending the fiberwise
symplectic forms: see~\cite{Sei} or~\cite{MS1} Chapter~6 for example.
Conversely, every pair consisting  of a smooth bundle $\pi: P\to S^2$
with fiber $M$ together with a closed $2$-form $\Om$ on $P$
that is nondegenerate on  each fiber arises in this way
from a loop in $\Ham(M, \om)$.
By adding to $\Om$ the
pullback of a suitable area form on the base, we may assume that $\Om$ is
nondegenerate.  Any such symplectic extension of the fiberwise forms will
be called $\om$-{\bf compatible.}  The set of these forms is contractible.
Note that each such $\Om$ gives rise to a connection on $P$ with
Hamiltonian holonomy, whose horizontal distribution consists of the
$\Om$-orthogonals to the fibers.  

Each such triple $(P,\pi,\Om)$ admits a contractible family 
$\Jj(P,\pi,\Om)$ of 
$\Om$-compatible almost complex structures $\TJ$  such that 
$\pi:(P,\TJ)\to (S^2,j_0)$ is holomorphic.  Each  
$\TJ\in\Jj(P,\pi,\Om)$ preserves 
the tangent bundle to 
the fibers and hence also the horizontal distribution.

Now observe that the
bundle $(P_\La,\Om) \to S^2$ supports two canonical
cohomology classes.
The first is the {\bf first Chern class of the vertical tangent bundle}
$$
c_{\ver} = c_1(TP_\La^{\ver})  \in H^2(P_{\La},\Z).
$$
The second is the {\bf  coupling class}, which is the unique class 
$u_\La \in H^2(P_\La,\R)$ such that
$$
i^*(u_\La) = [\om],\qquad u_\La^{n+1} = 0,
$$
where $i: M\to P_\La$ is the inclusion of a fiber.

Another important geometric fact about Hamiltonian bundles over $S^2$ 
is that they always have sections.  A direct geometric argument
shows that this is equivalent to saying that the map 
$$
\pi_1(\Ham(M,\om))\to \pi_1(M,x):\quad \{\phi_t\}\mapsto \{\phi_t(x)\}
$$
given by evaluation at the base point $x$ is trivial.  The latter 
statement follows  from the proof of the Arnol'd conjecture or by the 
very existence of the Seidel representation: see~\cite{LMP}.
Therefore, in particular, there always exists
a {\bf section class}, that is,  a class $\sigma \in  H^S_2(P,\R)$ 
that projects onto the positive generator of $H_{2}(S^{2}, \Z)$.
We shall denote by $H_2^\sect(P;\Z)$ the affine subspace of $H_2(P;\Z)$ 
consisting of such section classes.

\subsubsection*{Circle actions}

We now assume that $\Lambda_K$ is a circle action with moment map $K$,
and show how the ideas above simplify in this case. 
Let $S^1$ act on $S^3\times M$ by the diagonal action
$$
(z_1,z_2; x)\mapsto (e^{2\pi it}z_1, e^{2\pi it}z_2; \phi_t x).
$$
We claim that $P_\La$ can  be identified with the quotient $S^3\times_{S^1} M.$
To see this,  write $[z_1,z_2; x]$ for the equivalence class containing 
the point $(z_1/r,z_2/r; x) \in S^3 \times M$, where $r^2 = |z_1|^2 + |z_2|^2$.
In these coordinates,
$$
D_0\times M = \{[z,1;x]: |z|\le 1,\; x\in M\}
\quad \mbox{and} \quad   D_\infty\times M 
= \{[1,z;x] : |z|\le 1, \,x\in M\} .$$
In particular, $M_0$ is the fiber at $[0:1]\in \CP^1 = S^3/S^1$, and
$M_\infty$ is the fiber at  $[1:0]$.  Both have natural 
identifications with $M$.
Since the 
orientation on $D_\infty$ was reversed, the gluing map  is given by
$$
[1,e^{-2\pi it}; x] \sim [e^{2\pi it},1; \phi_tx],
$$
as required.

Let $\al \in \Omega^1(S^3)$ be the usual contact form on the unit sphere,
normalized so that  $d\al= \chi^*(\tau)$ where 
$\chi:S^3\to S^2$ is the Hopf map and $\tau$ is 
the standard area form on $S^2$ with total area $1$.
Given any $c \in \R$, the form 
$\omega + c\,d\alpha - d(K \alpha) \in \Omega^2(S^3 \times M)$ 
is closed and
 basic, and  hence descends under the projection
$\pr: S^3\times M\to S^3\times_{S^1} M$ to a closed two form
$\Om_c$ on $P_\Lambda$ which extends the fiberwise symplectic 
form.  Thus,
\begin{equation}\label{eq:omc}
    \Om_c = \pr_*(\om + d((c-K)\al)) = \pr_*(\om + c\,d\al - d(K\al)).
\end{equation}
If $c > \max K$, then $\Om_c$ is symplectic.
The coupling class is simply $[\Om_0]$:
 \begin{equation}\label{eq:om0}
u_\La =  [\Om_0] = [\pr_*(\om - d(K\al))]. 
 \end{equation}

\begin{remark}\labell{rmk:monotone}\rm
Note that, if $c_1$ and $[\omega]$ are linearly dependent,
then, since the $S^1$-orbit of an arc going
from the minimum to the 
maximum of 
$K$ is a sphere on which both $\om$ and $c_1$ are positive, $(M,\omega)$
must be monotone, that is,  $I_c = \mu I_\om$ for some $\mu> 0$. 
\end{remark}

Each  fixed point $x$ of the
$S^1$-action gives rise to a  section of $P$
$$
\si_x: = S^3 \times_{S^1} \{x\} = 
D_0\times \{x\} \cup D_\infty\times \{x\}.
$$
We will sometimes write $\si_F$ or $\si_{\max}$ instead of 
$\si_x$, when $x\in F$ or $x\in F_{\max}$, respectively.
Here are some useful facts about these sections.

\begin{lemma}\label{le:Ax} 
If $x$ is a fixed point of a Hamiltonian circle action $\Lambda_K$ on a
symplectic manifold $(M,\omega)$, then
$$
c_\ver(\sigma_x) = m(x) \quad \mbox{and} \quad u_\Lambda(\sigma_x) = -K(x) .
$$
Moreover, if $B$ is the class of the sphere formed by the $\La$-orbit 
of an arc from $x$ to another fixed point $y$, then
$B = \si_x - \si_y$.
\end{lemma}

\begin{proof}\, 
The normal bundle of $\si_x$ can be identified with a sum of 
holomorphic line bundles $L_i\to \CP^1$, one for each weight $m_i$ at $x$.  
Moreover, $c_1(L_i) = m_i$.  Thus $c_{\ver}(\si_x) = m(x)$.
Further, by Equation~(\ref{eq:om0}) 
$u_\La(\si_x) = [\om - d(K\al)](\si_x) = -K(x)$.  
This proves the first claim. 

Using the sign conventions explained at the beginning of \S\ref{ss:sb},
one finds by an easy calculation 
that $\omega(B) = K(y) - K(x) = u_\Lambda(\sigma_x - \sigma_y) =
\omega(\sigma_x - \sigma_y)$.
This identity
holds for all $\La$-invariant symplectic forms $\om'$ on $M$. 
But, after averaging,  any closed $2$-form sufficiently close to $\om$ 
is a $\La$-invariant symplectic form.
Hence the classes $[\om']$ fill out an open neighborhood of $[\om]$ in $H^2(M)$,
and so $B = \sigma_x - \sigma_y$.
 \end{proof}

Let $J$ be any $S^1$-invariant almost complex structure on $M$.
The standard complex structure $J_0$ on $\C^2$ is also $S^1$-invariant (under the
diagonal action),
and its restriction to $S^3$ preserves the contact planes $\ker\al$.
Moreover, each vector 
$\xi \in T_pP_{\Lambda}$ can be considered as an
equivalence class of vectors on $T ( S^3 \times M)$; each such equivalence
class has a unique representative in 
$\ker \al \oplus TM$
at each point in the $S^1$-orbit $\pr^{-1}(p)$.
Therefore, the product complex structure
$J_0 \times J$ on $\ker\al \oplus TM$
descends to an almost complex structure $\TJ$  on $P_\Lambda$.  
By construction, if $J$ is compatible with $\omega$,
then $\TJ$ is compatible with $\Omega_c$ for all $c > \max K$.
Moreover, $\TJ$ preserves the tangent spaces
to the fibers, and the section $\sigma_x$ is holomorphic for all fixed $x$.

\begin{defn}\label{def:JS}  
We define $\Jj_S(M)$ to be the set of all 
$S^1$-invariant
$\om$-compatible almost complex structures on $M$, and
denote by $\Jj_S(P)$ the space of almost
complex structures on $P$ constructed as above from the elements 
$J\in \Jj_S(M)$.  Note that $\Jj_S(P)\subset \Jj(P,\pi,\Om_c)$ for 
all $c>\max K$.
\end{defn}

\subsection{Small quantum homology}\labell{ss:qh}

We shall work with quantum homology with coefficients in the ring
$\La: = 
\La^{\rm univ}[q,q^{-1}]$ where $q$ is a variable of degree $2$ and
$\La^{\rm univ}$ is a generalized Laurent series ring in a 
variable $t$ of degree $0$:
$$
\La^{\rm univ}: = \Bigl\{\sum_{\ka\in \R}
 r_\ka t^\ka\; \Big| \; r_{\ka}\in \Q,\;\;
\#\{\ka > c\,|\,r_\ka \ne 0\}< \infty,\forall c\in \R\Bigr\}.
$$
Correspondingly, quantum cohomology  has coefficients in the dual ring
$$
\check\La: = 
{\check\La}^{\rm univ}[q,q^{-1}]
$$
 where $q$ is as before and
$$
{\check\La}^{\rm univ}: = \Bigl\{\sum_{\ka\in \R}
 r_\ka t^\ka\;\big|  \; r_{\ka}\in \Q,\;\;
\#\{\ka < c\,|\,r_\ka \ne 0\}< \infty,\forall c\in \R\Bigr\}.
$$
Thus we define 
$$
\QH_*(M;\La)= H_*(M, \Q)\otimes_{\Q}\La,\qquad
\QH^*(M;\check\La)= H^*(M, \Q)\otimes_{\Q}{\check\La}.
$$
These rings are $\Z$-graded  in the obvious way:
\begin{equation}\label{eq:grad}
\deg(a\otimes q^d\,t^{\ka}) = \deg(a) + 2d,
\end{equation}
where $a\in H_*(M)$ or $H^*(M)$.
They also have $\Z/2\Z$-gradings in which the even part is 
strictly commutative; for example,
$$
\QH_{\ev}: =
H_{\ev}(M)\otimes\La, \quad \QH_{\odd}: =
H_{\odd}(M)\otimes\La.
$$

Recall that the {\bf quantum intersection product}
$$
a*b\in \QH_{i+j - \dim M}(M;\La), \qquad \mbox {for }\; a\in H_i(M), b\in
H_j(M)
$$
 is defined as follows:
$$
a* b = \sum_{B \in H_2^S(M;\Z)} (a*b)_B\otimes q^{-c_1(B)}\,t^{-\om(B)},
$$
where  $(a*b)_B\in H_{i+j- \dim M+2c_1(B)}(M)$ is defined by the
requirement that
$$
(a*b)_B\,\cdot_M\, c = \GW_{B,3}^M(a,b,c) \quad\mbox{ for all }\;c\in
H_*(M).
$$
Here  $\GW_{B,3}^M(a,b,c)\in \Q$ denotes the Gromov--Witten invariant that counts
the number of spheres in $M$ in the class $B$
that meet cycles representing the
classes $a,b,c\in H_{*}(M)$. 
The product  $*$ is
extended to $\QH_*(M)$ by linearity over $\La$, and is associative.
Moreover, it  respects the $\Z$-grading.

This product $*$ gives $\QH_{*}(M;\La)$
the structure of a
graded commutative ring with unit $\1 = [M]$. Further, the invertible
elements in $\QH_{\ev}(M;\La)$ form a commutative group
$\QH_\ev(M;\La)^\times$ that acts on  $\QH_*(M;\La)$ by  quantum multiplication.

We shall work mostly with quantum homology since this is more geometric.
However, some examples mention quantum cohomology.  The
multiplication  (quantum cup product) is defined via Poincar\'e
duality: given $\al, \be \in H^*(M)$ with Poincar\'e duals $a =
\PD(\al),b = \PD(\be)$
$$
\al *\be = \PD(a*b) = \sum_{B\in H_2^S(M;\Z)} \PD((a*b)_B) 
        \otimes q^{c_1(B)}\,t^{\om(B)}.
$$
Note that the coefficient is $q^{c_1(B)}\,t^{\om(B)}$ rather than
$q^{-c_1(B)}\,t^{-\om(B)}$: in general 
the Poincar\'e duality map $\PD: \QH^*(M) \to \QH_*(M)$ is given by
$\PD(\al\otimes q^d\,t^\ka) = \PD(\al)\otimes q^{-d}\,t^{-\ka}$.   Thus in 
cohomology we must use the dual $\check{v}$ of the valuation $v$, namely
\begin{equation}\label{eq:vcheck}
\check v \Bigl(\sum_{d,\ka} a_{d,\ka} \otimes q^d t^{\ka}\Bigr) = 
\min \{\ka\;|\;\exists \, d :  a_{d,\ka}\ne 0\}.
\end{equation}

\subsection{The Seidel representation}\label{ss:se}

In this paper, we will study Hamiltonian loops $\Lambda$
by examining the geometry of the symplectic bundle $P_\Lambda$.  
On the classical level, this can be done by examining
the (generalized) Weinstein homomorphism, 
which we mentioned in Lemma~\ref{le:actind}.
\MS

\NI
{\bf Proof of Lemma~\ref{le:actind}.}\,\,  
Define the map
$$
\Aa_{\om,c}:   \pi_1(\Ham(M,\om)) \to \R\oplus \Z/\im (I_\om 
\oplus I_c)
$$
by
$$
\Aa_{\om,c}(\La) = -[u_\La(\si), c_{\ver}(\si)],
$$
where $\si: S^2 \to P_\La$ is any  section.
Note that this map is well defined.
From the construction of
the classes $u_\La$ and $c_{\ver}$,  and from the fact that $P_{\La_1+\La_2}$ is 
the fiber sum $P_{\La_1}\sharp_M P_{\La_2}$, it is easy
to see that $\Aa_{\om,c}$ is a homomorphism: see~\cite[Lemma 3.E]{LMP}.
Finally, if $\La$ is a circle with moment map $K$,
then  Lemma \ref{le:Ax} implies that  
$u_\La(\si_x) = -K(x)$ and $c_{\ver}(\si_x) = m(x)$
for each fixed point $x$.
\QED

\begin{definition}
We 
define the {\bf Seidel element} $\Ss(\Lambda) \in \QH_{\dim M}(M;\La)$  by
\begin{equation}\label{eq:psi}
\Ss(\La) = \sum_{\si\in H_2^\sect(P)} a_{\si}\otimes q^{-c_{\ver}(\si)} 
\;t^{-u_{\La}(\si)} 
\end{equation}
where 
$
a_{\si}\cdot_{M} c = \GW^{P_{\La}}_{\si,1}(c)
$
for all $c\in H_{*}(M)$.\footnote
{
We use a $1$-point Gromov--Witten invariant here.
Similarly, we define $\Ss(\La)(a)$ using a $2$-point invariant. 
Because $[M]\cdot[\si] = 1$ for any section class, the divisor axiom 
for GW invariants implies that the $1$-point 
invariant $GW^P_{\si,1}(c)$ equals the more usual $3$-point invariant
$GW^P_{\si,3}([M],[M],c)$.  However, it is sometimes more convenient 
to use the $1$-point 
invariant because the moduli space $\Mm_{0,1}$ can be compact
while $\Mm_{0,2}$ never is because the two marked points must always 
be distinct.
}
Here $H_2^{\sect}(P)$ denotes the affine subspace
of $H_2(P;\Z)$ that is represented by sections.
\end{definition}

Intuitively, $a_{\si}$ is represented by the class
$$
\ev_*\bigl(\Mm_{0,1}(P_\La,\TJ;\si)\bigr)\cap [M]
$$
where $\Mm_{0,1}(P_\La,\TJ;\si)$ is the  moduli space of
all $\TJ$-holomorphic sections in class
$\si$ with one marked point, $\ev$ is the obvious evaluation map
to $P_\La$ and $[M]$ denotes the homology class represented by a fiber; 
see~\cite{Mcq}.
This moduli space  has formal dimension $\dim M + 2c_\ver(\si) + 2.$
We find that $a_\si = 0$ unless 
$$
\deg(a_\si \otimes q^{-c_{\ver}(\si)}) =  \deg(a_\si) - 2c_{\ver}(\si) = \dim M.
$$
Because all dimensions are even,
$\Ss(\La)$ belongs to the strictly commutative part
$\QH_{\ev}$ of $\QH_{*}(M)$.
Moreover, $\Ss(\Lambda)$ is independent of 
the choice of symplectic extension form
$\Om$ since all of these are deformation equivalent.
It is shown in~\cite{Mcq} (using ideas from~\cite{Sei,LMP}) that
 $\Ss(\Lambda)$ lies in $\QH_\ev(M;\La)^\times,$ 
the group of multiplicative units in the ring $\QH_{\ev}(M)$,
and that the correspondence $\Ss$ induces a group homomorphism
$$
\Ss: \pi_{1}(\Ham(M,\om))\;\;\to\;\; \QH_\ev(M;\La)^\times.
$$
It is immediate from the definition that it lifts the Weinstein 
homomorphism.

It is often useful to identify $\QH_\ev(M;\La)^\times$ with 
$\Aut(\QH_*(M;\La))$,
the group of automorphisms of  $\QH_*(M;\La)$ as a right $\QH_*(M;\La)$-module, since 
every such automorphism is determined by its value at $\1$.
Correspondingly we define
$$
\Ss(\La)(a): = \Ss(\La) * a \ \quad \forall \ a \in \QH_*(M;\La).
$$ 
Since the Seidel element has degree $\dim M$, this endomorphism 
preserves degree.

\begin{defn} The {\bf Seidel representation} is the group homomorphism
$$
\Ss:\pi_1(\Ham)\to \QH_\ev(M;\La)^\times = \Aut(\QH_*(M;\La)).
$$
\end{defn}

It is shown in~\cite{Mcq} (see also~\cite{Sei,LMP,MS2}) that
\begin{equation}\label{eq:psila}
\Ss(\La)(a) = \sum_{\si\in H_2^\sect(P)} b_{\si}\otimes q^{-c_{\ver}(\si)} 
\;t^{-u_{\La}(\si)} 
\end{equation}
where 
$
b_{\si}\cdot_{M} c = \GW^{P_{\La}}_{\si,2}(a,c)
$
for all $c \in \QH_*(M)$.
Here one should think of $a$ as
represented by a cycle in the fiber $M_0$ over
the center of the disc $D_0$, 
and $b_\si$ and $c$ as represented by cycles
in the fiber $M_\infty$ over the center of
$D_\infty$. Then the element $\Ss(\La)$ induces a ring isomorphism from
$\QH_*(M_0)$ to $\QH_*(M_\infty)$.
Intuitively, the class $\Ss(\La)(a)$ is represented by
the intersection of $M_\infty$
and the space of all $J$-holomorphic sections of $P_\La$ that
meet the cycle in $M_0$ which represents $a$.
Since the connection in  the bundle $(P,\Om)\to S^2$ provides an
identification of $M_0$ with $M_\infty$ that is well defined up to
symplectic isotopy, 
$\Ss(\La)$  gives rise to a
well defined  element of $\Aut(\QH_*(M;\La))$ as claimed.

\begin{example}\label{ex:seid2}\rm
Consider the rotation of the unit sphere
$S^2$ with $K = 2\pi x_3$.
Then the fibration $P_{\La}$
can be identified with the nontrivial fibration from the one point
blow up $M_*$ of $\CP^2$ to $S^2$. By Lemma~\ref{le:Ax}  the section $\si_{\max}$
corresponding to the maximum (the north pole)
has normal bundle of Chern number $m(n) = -1$, and so 
is the exceptional divisor, while the section
$\si_{\min}$ corresponding to the minimum (the south pole) has Chern
number $1$, and so lies in the class of a line. 
Since the the Seidel element $\Ss(\Lambda_K)$
has degree $\dim M = 4$, a section $\sigma$ can only contribute 
to it if $0 \geq 2 c_\ver(\sigma) = 2c_1(X) \geq -4$.
Therefore, $\sigma_{\max}$ is the only holomorphic section of 
$P_{\La}$ that can contribute to the Seidel element.
It follows easily that $\Ss(\La_K) = [pt]\otimes q\,t^{\om(A)/2}$, 
as claimed in Example~\ref{ex:s2}.
\QED\end{example}

\section{Computing the Seidel element}

This section contains the main proofs.
We begin by calculating the contribution to the Seidel element 
$\Ss(\La)$ of the sections $\si_{\max}$ through points on
the maximal fixed set $\Fmax$.  In Proposition~\ref{prop:amax}
we show that this is nonzero precisely when
$\Fmax$ is homologically visible.
These arguments use easy results on the
behavior of $J$-holomorphic spheres. To go further, we need 
a version of the localization theorem: 
there is a $T^2$-action on the moduli spaces of stable maps and only 
the invariant elements contribute to $\Ss(\La)$.
This theorem is stated in \S\ref{ss:chain}.  We  defer the proof to 
\S\ref{ss:local}, devoting the rest of this section to 
an investigation of the invariant elements.  We first prove 
Theorem \ref{thm:max}.   Then, in \S\ref{ss:semifree}, we consider
the semifree case, and prove Theorem \ref{thm:semifree}.
Finally, in \S\ref{ss:twofold}, we consider the  
case where the isotropy is at most twofold, and prove
Theorems \ref{thm:2iso} and \ref{thm:pfasimplemax}.

\subsection{The contribution of the maximal fixed 
set}\label{ss:maxim}

We begin with some preliminary remarks about $\TJ$-holomorphic 
sections of $P: = P_\La$. Throughout we assume $\TJ\in 
\Jj_S(P_\La)$, the space  of
almost complex structures on $P_\Lambda$ that are 
constructed from $S^1$-invariant almost complex structures on $M$ using the 
identification of $P_\Lambda$ with a quotient of $S^3\times M$.
See Definition~\ref{def:JS}.

Let $\Mm_{0,k}(P,\TJ;\si)$ denote the space of equivalence classes
$[u,\bz]$ of $\TJ$-holomorphic  maps $u:S^2\to P$ in class $A$ 
with $k$ pairwise distinct marked 
points $\bz:= \{z_1,\dots,z_k\}$. 
Here,
two such pairs $(u,\bz)$ and $(u',\bz')$
 are equivalent if there is $\psi\in \PSL(2,\C)$ such that
$$
u' = u\circ \psi,\qquad \psi(z_i') = z_i, \quad i = 1,\dots,k.
$$
The  compactification $\oMm_{0,k}(P,\TJ;\si)$ consists of equivalence classes
$\tau= [\Si(\bu),\bu,\bz]$ of $\TJ$-holomorphic stable maps $\bu:\Si(u)\to P$
with $k$ marked points.  Here $\Si(\bu)$ is a  union of copies of 
$S^2$ attached via a tree graph, and the equivalence relation 
is given by all reparametrizations that respect the special points, 
i.e. the attaching (or nodal) points and the marked points.  
If $\si$ is a section class, each element $\tau$ in 
$\oMm_{0,k}(P,\TJ;\si)$ projects via $\pi:P\to S^2$
to an equivalence class of  holomorphic maps $\pi\circ \bu: \Si(u)\to S^2$
of total degree $1$.
Such a map has just one  component of degree $1$; on all 
the other components $\pi\circ \bu$ is constant.
Thus $\tau$ has a distinguished component that is a section, 
called the {\bf root}. 
The other components are mapped into fibers.

We shall be specially interested in the case when $k=2$ and 
the first marked point is mapped to $M_0$, the other to $M_{\infty}$.  In this 
 case, there is a unique chain of spheres joining the component that 
 contains the first marked point $z_0$ to the component
that contains the second marked point $z_{\infty}$;  
we call the components of this chain the {\bf principle components}. 
 The other 
 spheres are called {\bf bubbles}.  The root is always a principle 
 component. 
 
For most of the results in this paper, we will need to look at invariant
chains, as described in the next subsection.  
However, the  following observation, due to Seidel\footnote
{Private communication.}, allows us to give a simpler argument
when we are studying curves in a class $\si_{\max} + B$
with $\omega(B) \leq 0$.  

\begin{lemma}\labell{le:Fmax} Consider a Hamiltonian circle action
$\Lambda_K$ on a compact symplectic manifold $(M,\omega)$,  and
let $\TJ\in \Jj_S(P_\La)$ be constructed from $J\in \Jj_S(M)$. 
Fix $B \in H_2^S(M;\Z)$, and consider the moduli space
$$
\oMm_{0,0}(P_\Lambda,\TJ;\si_{\max} + B).
$$ 
\begin{itemize}
\item If $B \neq 0$ and $\om(B)\le 0$, the moduli space is empty.
\item If $B = 0$, the moduli space is compact and
can be identified with $F_{\max}$ itself.
\end{itemize}
\end{lemma}
\begin{proof}\,  
The symplectic form $\Omega_c$ defined in~(\ref{eq:omc})
is compatible with $\TJ$
for any $c > \max K$.
Fix $[z,w]  \in P_\Lambda = S^3 \times_{S^1} M$.
Recall that any non-zero tangent vector  $\xi \in T_{[z,x]} P_\Lambda$
can be  uniquely represented by  a vector $h + v \in T_{(z,x)}(S^3 \times M)$,
where $h \in \ker \alpha \subset T_z S^3$ and $v \in T_x M$.
Now
\begin{eqnarray*}
    \Omega_c(\xi,\TJ\xi) & = & 
\bigl(\omega - dK \wedge\alpha + (c - K) d\alpha\bigr)(h + v, 
J_0h + Jv) \\
& = & \omega(v,Jv)  + (c - K) \,d\alpha(h, J_0h) \\
&\geq &  (c - K_{\max}) \,d\alpha(h,J_0h)\\
& = & (c - K_{\max})\chi^*\tau(\xi,\TJ\xi)
\end{eqnarray*}
with  equality impossible unless $v = 0$ and $K(x) = \max K$.
Since  $\chi^*\tau$ is the pullback by the Hopf map of the area form 
on $S^2$ with area $1$, it follows that for any 
 $\TJ$-holomorphic section $\si$
$$
 \Omega_c(\sigma) \geq  
c - K_{\max} = \Omega_c(\sigma_{\max}),
$$
with equality occurring exactly if $\sigma$ is a constant section
$\sigma_x$ for some $x \in F_{\max}$. 

 Since every stable map in a section class $\si$ either consists
 of a section, or is the union of a section with other spheres $A_i$
 which lie in the fibers and satisfy $\omega(A_i)  > 0$, 
the only stable maps that represent a section class $\si$ with 
$\Om_c(\si)\le c-\max K$ are the constant sections $\si_x, x\in \Fmax$.
The result follows.
\end{proof}

\begin{lemma}\labell{le:reg}  Let $x$ be any fixed point of the $S^1$-action.
For each $\TJ\in \Jj_S(P_\La)$, 
the $\TJ$-holomorphic curve $\si_x$ is regular 
precisely when the negative weights at $x$ are all equal to $-1$.
\end{lemma}

\begin{proof}\, 
Recall that $\si_x$ is regular if and only if the 
linearization $D_u$ of the corresponding 
Cauchy--Riemann operator is surjective. 
When $J\in \Jj_S(P_\La)$ the normal bundle of $\si_x$ is holomorphic
and splits into a sum of line bundles $\oplus_i L_i$ that are 
preserved by $D_u$.
Moreover, $D_u$ restricts on each $L_i$ to the 
usual Dolbeault delbar
operator. Thus $D_u$ is surjective precisely when $c_1(L_i) 
\ge -1$ for all $i$.
 \end{proof}

The next proposition generalizes part (i) of Theorem~\ref{thm:max}.

\begin{prop}\label{prop:amax}
Consider a Hamiltonian circle action $\La_K$  on a
compact symplectic manifold $(M,\omega)$ with normalized moment map $K$.
Let 
$e(\Ee_{\max})\in H_*(\Fmax)$ denote the
Poincar\'e dual of the Euler class of the obstruction 
bundle at $\Fmax$  (see equation~(\ref{eq:ee})), 
denote the inclusion $H_*(\Fmax)\to H_*(M)$ by $\io$,
and set  $K_{\max}: = K(\Fmax)$ and $m_{\max}: = m(\Fmax)$. 
Then:
$$
\Ss(\La_K) = \io\bigl(e(\Ee_{\max})\bigr)
  \otimes q^{-m_{\max}}\, t^{K_{\max}} +  
\sum_{B \in H_2^S(M): \omega(B) > 0} a_B \otimes 
q^{-m_{\max}-c_1(B)}\, 
t^{K_{\max}-\om(B)}.
$$
\end{prop}

\begin{proof}\,
By Lemma~\ref{le:Ax}, 
$u_\Lambda(\sigma_{\max}) = -K_{\max}$ and 
$c_{\ver}(\sigma_{\max}) = m_{\max}.$
Therefore we may write
$$
\Ss(\La_K) = \sum_{B \in H_2^S(M;\Z)} a_B\otimes q^{-m_{\max} - 
c_1(B)}\,t^{K_{\max}-\om(B)},
$$
where $a_B$  is the contribution from the section class $\si_{\max} + B$.
Fix  any $c \in H_*(M)$.
It is enough to show that 
$a_0 \cdot_M c = \io(e(\Ee_{\max})) \cdot_M c$,
and that $a_B \cdot_M c = 0$ for every 
nonzero $B \in H_2^S(M;\Z)$ such that $\omega(B) \leq 0$.
By definition, $a_B \cdot_M c = \GW^{P_\Lambda}_{\sigma_{\max} + B,1}(c)$.
Choose an almost complex structure $\TJ \in \Jj_S(P)$.
By Lemma~\ref{le:Fmax},
if $\omega(B) \leq 0$ and $B \neq 0$
then the moduli space $\Mm_{0,0}(P_\Lambda,\TJ,\si_{\max} + B)$ 
is empty, and so $a_B \cdot_M c = 0$, as required.
On the other hand, $\Mm_{0,1}(P_\Lambda,\TJ,\si_{\max})$ 
can be identified with the compact manifold $S^2\times F_{\max}$.  
(The $S^2$-factor
is the locus of the single marked point.)

If $\Fmax$ is semifree, then  Lemma \ref{le:reg} implies
that $\sigma_x$ is regular for every $x \in F_{\max}$.
Hence the intersection of the evaluation pseudocycle
$
\ev:\Mm_{0,1}(P,\TJ,\si_{\max})\to M_{\infty}$
with any class $c$ in the fiber 
$M_\infty$ is precisely $ [F_{\max}] \cdot c$.  Thus $a_0 = [\Fmax]$ 
in this case.

If any of the negative
weights $-k_i$ at $\Fmax$ is less than $-1$, the
  elements of the compact manifold $\Mm: = \Mm_{0,1}(P,\TJ,\si_{\max})$
  are not regular.  Rather, for each $i$, the cokernel of the restriction
  $D_{u_x}: C^{\infty}(S^2, E_i)\to \Om^{0,1}(S^2,E_i)$ is a vector 
  space of dimension  $\dim E_i\otimes \C^{k_i-1}$, and as $x$ varies 
  in $\Fmax$ these cokernels fit together to form the bundle $E_i\otimes 
  \C^{k_i-1}$ over $\Mm$.  Thus the total obstruction bundle is the 
  bundle
  $\Ee\to\Mm$ of equation~(\ref{eq:ee}).  It follows from the standard theory (see for 
  example~\cite[\S5.3]{Mcv} or~\cite[Chapter~7.2]{MS2}) that the 
  regularized moduli space corresponds to the zero set of a generic 
  section of $\Ee=:\Ee_{\max}$.   Therefore 
  $\GW^{P_\Lambda}_{\sigma_{\max},1}(c)=  
 \io(e(\Ee))\cdot_M c$ 
  for each $c\in 
  H_*(M)$,
  and the result follows.
  \end{proof}

\subsection{Invariant beads and chains}\label{ss:chain}

In order to understand the moduli spaces of sections in 
an arbitrary class $\si$ we exploit the fact that $T^2$ acts 
on $(P_\La, \Om, \TJ)$
when $\TJ\in \Jj_S(P_\La)$. 
Here the first factor $S^1\times \{1\}$ acts 
on $P_\La$
by rotating the fibers via  $\phi_t$ while the second factor 
$\{1\}\times S^1$ acts by rotating the base as 
follows:
$$
\theta\cdot[z,1;x] = [e^{2\pi i\theta}z,1; \phi_tx],\quad 
\theta\cdot[1,z;x] = [1,e^{-2\pi i\theta}z; x].
$$
Note that
the only points of $P_\La$ 
fixed by the whole group are the points in $M_0$ and 
$M_\infty$ that are fixed by the original $S^1$-action $\phi_t$.  
Because the elements of $\Jj_S(P_\La)$ are constructed from $S^1$-invariant
almost complex structures on $M$ (see Definition~\ref{def:JS}),
this action preserves $\TJ$.  Hence $T^2$
acts on the moduli 
spaces of stable maps  via postcomposition.

The next result is a version of the localization principle for 
$T^2$-actions; it is well known in the algebraic case 
and is proved in the symplectic situation 
in \S\ref{ss:local}. 
Given two (weighted)
pseudocycles $f:Z\to P_\Lambda$ and $f':Z'\to P_\Lambda$
(see Definition~\ref{def:psc})
and a section class $\si$, we 
define
$$
\oMm_{0,2}(P_\La, \TJ, \si; Z,Z'): = \ev^{-1}(\ov{f(Z)}\times \ov{f'(Z')})
$$
where $\ev:\oMm_{0,2}(P_\La,\TJ,\si)\to P_\La\times P_\La$ is the 
evaluation map. 
The pseudocycles are said to be  $\bf S^1${\bf-invariant}
if the  images $f(Z)$ and $f'(Z')$ 
are closed under the action of $S^1$.   In this case, 
if $\TJ \in \Jj_S(P_\La)$,
then clearly there is an induced action of 
$T^2$ on this cutdown moduli space.

\begin{prop}\labell{prop:invar}
Suppose that $f:Z\to M_0$ 
and $f':Z'\to  M_{\infty}$ 
are $S^1$-invariant 
weighted
pseudocycles
which represent the classes $a$ and $a'$ in $H_*(M)$, respectively.
Given  $\TJ\in \Jj_S(P_\La),$  write
$$
\Ss(\La)(a) = \sum_{\si\in H_2^{\sect}(P)} a_\si\otimes 
q^{-c_{\ver}(\si)}t^{-u_{\La}(\si)}.
$$
Then $a_\si\cdot_M a' =  0$ unless 
the moduli space ${\oMm}\,\!^{cut}: = 
\oMm_{0,2}(P_\La, \TJ,\si; Z,Z')$ contains a 
$T^2$-invariant element.
Moreover, $a_\si\cdot_M a'$ is a sum of contributions, one from each connected 
component of the space $({\oMm}\,\!^{cut})^{T^2}$ of invariant elements.
\end{prop}

Note that most $T^2$-invariant elements in  
$\oMm_{0,2}(P_\La, \TJ,\si; Z,Z')$ are not regular.
Therefore it would be a nontrivial task to calculate their actual
contributions to the invariant.  In this paper we do not attempt such 
calculations.

The next task is to figure out the structure of the $T^2$-invariant 
elements in $\oMm_{0,2}(P_\La, \TJ,\si)$.  Note that each principal 
component has $2$ special points joining it to the other principal
components. 
We will place these at $0$ and $\infty$ and then identify 
$S^2\oursetminus \{0,\infty\}$ with the cylinder $(s,t)\in \R\times S^1$ with
complex structure $j_0$ defined by $j_0(\p_s) = \p_t$.  

\begin{lemma}\labell{le:invar}
    Let $\TJ\in \Jj_S(P_\La)$ be constructed from $J\in \Jj_S(M)$ and 
    denote  by $g_J$ the metric on $M$ defined by $J$ and $\om$.
\smallskip

\NI
{\rm (i)}  If $A$ is a section class the only elements in 
$\Mm_{0,2}(P_\La, \TJ, A)$ that are fixed by the $T^2$-action
have the form $[u;0, \infty]$ where $u:S^2\to P_{\La}$ is parametrized 
as a section and has as image some constant sphere $\si_x$ where
$x\in M^{S^1}$.  
\smallskip

\NI
{\rm (ii)} If $A\in H_2(M)$ then the only elements in 
$\Mm_{0,k}(P_\La, \TJ, A)$ that are fixed by the $T^2$-action lie 
in either  $M_0$ or $M_{\infty}$.   If 
such an element does not lie entirely in $M^{S^1}$,
then $k\le 2$ and there exists a parametrization 
$u: \R \times S^1 \to M$ and a path $\gamma:\R \to M$ which joins
two fixed points $x$ and $y$ in $M$ so that the marked points lie in 
$u^{-1}(\{x,y\}),$ 
and
\begin{equation}\label{eq:pqbead}
u(s,t) = \phi_{pt/q}\ga(s),\quad \mbox{and} \quad \ga'(s) = \frac{p}{q} 
\grad_{g_J} K,
\end{equation}
where $p \neq 0$ and where  $q > 0$ is the 
order of the isotropy group of the points 
in the image of $\ga$.  
There is a unique choice of parametrization such that 
\begin{equation}\label{eq:pqlimit}
\lim_{s\rightarrow -\infty} u(s,t) = x \quad \mbox{and} \quad
\lim_{s\rightarrow \infty} u(s,t) = y.
\end{equation}
\end{lemma}
\begin{proof}\,  Statement (i) is clear, as is the first claim in (ii). 
Thus, identifying $M_0$ and $M_\infty$ with $M$, we just need to 
understand the spheres  $u : S^2 \to M$ that are fixed by $\La$ 
and do not lie 
entirely in $M^{S^1}$. 
Suppose first that $u$ is simple, i.e. not multiply covered.  
For each $\theta\in S^1 = \R/\Z$, the composite
$\phi_\theta\circ u$ must be a reparametrization of $u$, that is,
there is a unique $\psi_{\theta}\in  \PSL(2,\C)$ such that
$\phi_\theta\circ u = u\circ \psi_\theta$.
It follows easily the assignment $\theta
\mapsto \psi_\theta$ defines 
 a homomorphism $S^1\to \PSL(2,\C)$.  Since the only circle subgroups of 
$\PSL(2,\C)$ consist of  rotations about a fixed axis,
there are two 
points in $S^2$, say $0,\infty$, that are mapped by $u$ into $M^{S^1}$.
If $\im\,u\cap M^{S^1} = \{x,y\}$ we may suppose $u(0)= x$ 
and
$u(\infty) = y$. If there are marked points on $u$ they must form a 
subset of $\{0,\infty\}$.
Using coordinates $(s,t)$ on $S^2 \oursetminus\{0,\infty\}=\R\times S^1$ as above, 
we find that for some  $q\ne  0$ 
$$
\psi_\theta(s,t) = (s,t + q\theta),\qquad 
\phi_\theta\circ u(s,t) = u(s, t + q\theta).
$$
Thus, $\im\,u$ lies in the set of points with isotropy group 
$\Z/(|q|)$.   
Denoting $\ga(s): = u(s,0)$, we have $u(s,\theta) = \phi_{\theta/q}\ga(s)$.
Moreover
$$
0 = \p_s u + J\p_t u= (\phi_{t/q})_*\bigl(\ga'(s) + 
\frac{1}{q}JX_{K}(\ga(s))\bigr),
$$
where $X_K$ is the  Hamiltonian flow induced by $K$.
Thus $\ga' = \frac{1}{q} \grad K$ because $-JX_K =\grad K$.  (Here we take the 
gradient with respect to the metric $g_J$.)
Since every sphere is the $|p|$-fold cover of a simple sphere, this 
proves (ii).  To get the stated result,
we absorb any negative sign into $p$ rather 
than $q$.
\end{proof}

\begin{defn}
Let $x$ and $y$ be two   
fixed points in $M$.

For  $q > 0$ and $p \neq 0$,  a {\bf bead from 
$\bf x$ to $\bf y$ of type $\bf (p,q)$} is a
map $u: \R \times S^1 \to M$ 
which satisfies equations~(\ref{eq:pqbead})  and ~(\ref{eq:pqlimit}).

For $q = 0$ and $p > 0$, 
a {\bf bead from $\bf x$ to $\bf y$ of type $\bf (p,q)$} is a
is a $p$-fold cover of a simple $J$-holomorphic sphere
$u: \R \times S^1 \to M$ that lies entirely in one component
of the fixed point set $M^{S^1}$ and 
which satisfies equation~(\ref{eq:pqlimit}).
\end{defn}

\begin{defn}
Given $x,y,z\in M^{S^1}$ an {\bf invariant principal chain from $\bf x$ to $\bf y$ 
in class $\bf \si_z+A$ and with root $\bf z$}  is a sequence of
critical points $x = x_1,x_2,\dots,x_k = y$ of $K$ joined by invariant 
$\TJ$-holomorphic spheres with the following properties:
\smallskip

\NI{\rm (a)} there is  $1\le i_0 \le k$ such that  $x_{i_0} = x_{i_0+1} = z$
and these points are joined by the section $\si_z$;
\smallskip

\NI{\rm (b)}
for each $1\le i < k$ where $i\ne i_0,$ the points
$x_{i}, x_{i+1}$ are joined by a $(p_i,q_i)$-bead in class $A_i$;
\smallskip

\NI{\rm (c)}  $\sum_{i\ne i_0} A_i = A$.
\smallskip

\NI
Further an {\bf invariant chain from $\bf x$ to $\bf y$ in class 
$\bf \si_z+A$ and with root $\bf z$} 
is a chain as above with additional ghost components 
at each of which a $T^2$-invariant tree of $(p,q)$ beads
is attached.  In this case,  $A$ is the sum of the classes represented by
the principal spheres and the bubbles.
\end{defn}

The next lemma is an immediate consequence of Lemma~\ref{le:invar} and the 
above definitions.

\begin{lemma}\labell{le:invariants}  Let 
$f:Z\to P_\La$ and $f':Z'\to P_\La$ be $T^2$-invariant 
pseudocycles, $\si$ be  a section class  and 
choose $\tilde{J} \in \Jj_S(P_\La)$.
Then every $T^2$-invariant element of the cut down
 moduli space $
\oMm_{0,2}(P_\La, \TJ, \si; Z,Z')$ is an invariant chain from a 
point $x\in \ov{f(Z)}$ to a point $y\in \ov{f'(Z')}$.
\end{lemma}

We will need the following useful facts about beads.
 
\begin{lemma}\labell{le:bead} 
Choose $J\in \Jj_S(M)$ and
 consider a 
 ($J$-holomorphic)
 $(p,q)$-bead from $x$ to $y$ in class $A$.  If $q \neq 0$, then  $A = 
 p(\si_x - \si_y)/q$.  Further:
\smallskip

\NI
{\rm (i)} If $K(y) > K(x)$ then $p> 0$,  $\om(A) 
 = p|K(x) - K(y)|/q$ 
 and $c_1(A) =  p(m(x) - m(y))/q$.
\smallskip

\NI
{\rm (ii)} If $K(y) < K(x)$ then $p< 0$, $\om(A) 
 = |p| |K(x) - K(y)|/q$ and $c_1(A) =  p(m(x) - m(y))/q$.
\smallskip

\NI
{\rm (iii)} If  $K(y) = K(x)$  then $\om(A) > |K(y) - K(x)|$.
 \end{lemma}

\begin{proof}\, 
    We saw in Lemma~\ref{le:Ax} that
the homology class of the sphere formed
    by the $\La$-orbit of an arc going from $x$ to $y$ is 
    $\si_x-\si_y$.  Hence 
    each $(p,q)$ bead from $x$ to $y$ lies in the class $A=p(\si_x-\si_y)/q$
    where $\om(A) = p(K(y)-K(x))/q$.  Statements (i), (ii) and (iii) 
    now follow from the fact that $\om(A) > 0$
    and $c_1(A) = p(m(x) - m(y))/q$.
\end{proof}

The proofs of our other results are based on a more careful 
study of the structure of the  
$T^2$-invariant elements in $\oMm_{0,2}(P_\La, \TJ, \si_{\La} + B).$
We begin by slightly strengthening the conclusion of
Proposition~\ref{prop:invar}.

\begin{lemma}\label{le:repr} 
Consider a Hamiltonian circle action $\La_K$  on a
compact symplectic manifold $(M,\omega)$ with normalized moment map $K$.
Let $\Fmax$ be the  maximal fixed component 
and choose $J \in \Jj_S(M)$.
Given $B \in H^S_2(M)$, let $a_B$ denote the contribution  of
 $\sigma_{\max} + B$ to
the Seidel element
 $\Ss(\Lambda_K)$. 
Then $a_B = 0$ unless $B$ can be represented by an invariant $J$-holomorphic
stable map that intersects $\Fmax$.
More generally, 
if $f':Z' \to M$ is an invariant pseudocycle representing the class 
$a'$, 
then $a_B\cdot_M a' = 0$  unless
$B$ can be represented
by an invariant $J$-holomorphic stable map that intersects both $F_{\max}$ and 
$\ov{f'(Z')}$. 
\end{lemma}

\begin{proof}{}
Assume $a_B \cdot_M  a'\neq 0$.
By Proposition~\ref{prop:invar}, there must be a 
$T^2$-invariant element in 
$$
\oMm_{0,2}(P_\La, \TJ, \si_{\max} + B; M_0, \ov{Z'}),
$$
where $\TJ\in \Jj_S(P_\La)$ is constructed from $J$ in the usual way.
Hence, 
by Lemma~\ref{le:invariants},
there is an invariant chain from $x \in M_0$ to $y \in 
\ov{f'(Z')}$
in the class $\sigma_{\max} + B$.
Let $z$ denote its root.
Let $A'$ be the sum of the homology classes represented by the
subchain of spheres in $M_0$ from $x$ to $z$, and let $A''$ be the sum of the
homology classes represented by the subchain of spheres in $M_{\infty}$
from $z$ to $y$.
Then $A' + A'' + \sigma_z = \sigma_{\max} + B$.
Since the  orbit of an upward gradient flow line from $z$ to $F_{\max}$ is 
$J$-holomorphic, the class  $\sigma_z - \sigma_{\max}$ is also represented 
by a $J$-holomorphic sphere.  
\end{proof}

\NI
{\bf Proof of Theorem~\ref{thm:max}.}
Part (i) is included in Proposition~\ref{prop:amax}.  
Part (ii) follows immediately from Lemma~\ref{le:repr}.

Now assume that $(M,J)$ is NEF
and that $2c_1(B') \geq  \codim \Fmax$ 
for all $J$-holomorphic spheres $B'$
that do not lie entirely in 
$\Fmax$.
Assume also that $a_B \neq 0$.
By Lemma~\ref{le:repr}, $B$  can be represented by
a $J$-holomorphic stable map which intersects $\Fmax$.
We must show that all components of this stable map lie in $\Fmax$.  
Suppose the contrary.
Then the assumptions imply
that $2 c_1(B)  \geq \codim \Fmax$.
On the other hand, $0 \leq \deg(a_B) =  \dim F_{\max} + 2c_1(B)  \leq \dim M$.
Therefore, $2c_1(B) = \codim F_{\max}$.
Therefore, $\deg(a_B) = \dim M$. 
Since $a_B$ is  not zero, this implies that it is a 
multiple of the generator of $H_{\dim M}(M)$; 
hence $a_B \cdot [pt] \neq 0$.
Choose 
$y \in F_{\min}$. Then since $a_B \cap [y] \neq 0$,
Lemma~\ref{le:repr} implies
$B$ can be represented by a $J$-holomorphic stable
map which intersects $\Fmax$ and $y$.
Let $B_1$ be a sphere in the corresponding
stable map which intersects $\Fmax$ at $x_1$
but does not lie entirely in $\Fmax.$
Let $x_2$ denote the second marked point in $B_1$. 
Let $B_2,\ldots B_k$ be the remaining $J$-holomorphic spheres in  $B$.
Then 
$$
\codim(\Fmax) = 2c_1(B) = 2\sum_{i=1}^k c_i(B_i).
$$
Since the assumptions imply that $c_1(B_i) \geq 0$ for all $i$
and  $2c_1(B_1) \geq \codim (\Fmax)$,  we conclude that
$2c_1(B_1) = \codim (\Fmax)$ and
$c_1(B_i) = 0$ for all $i \neq 1$.
Since $F_{\max}$ is semifree, $B_1$ is a bead of type  
$(p,q)$ with $q = 1$.  
By  Lemma~\ref{le:bead}(ii), $p< 0$ and
$$
2c_1(B_1) = 2p\, \bigl(m(F_{\max}) -m(x_2)\bigr) 
= -2p\, m(x_2) -p \,\codim(F_{\max}).
$$
Since $2c_1(B_1) = \codim (\Fmax)$,  $m(x_2) \leq 0$. 
Since $c_1(B_i) = 0$ for all $i \neq 1$, the next bead on the 
principal chain must
connect $x_2$ to another point, $x_3$, which also satisfies 
$m(x_3)  \leq 0 $.
Proceeding inductively, we see that 
$m(y) \leq  0$. 
But this is impossible, because
$m(y)> 0$ for all $y\in F_{\min}$.
\QED

\NI{\bf Proof of Proposition~\ref{prop:Kdist}.}
Since $\Ss(\La_K) = \1$, there is a class $B$ with 
$$
\om(B) =
- u_\La(\si_{\max}) = K_{\max} \quad \mbox{and} \quad c_1(B) = -c^\ver(\si_{\max}) = -m_{\max},
$$ 
such that $a_B\cdot_M [pt] \neq 0$.  By Lemma~\ref{le:repr} there is
an invariant $J$-holomorphic stable map  in class $B$ that intersects
$F_{\max}$ and $F_{\min}$.  Let $B_1,B_2,\ldots,B_j$ be the beads
in the principal chain 
which do not lie in a single fixed component.
Note that 
$
\omega(B) \geq \sum_{i=1}^j \om(B_i),
$
with equality impossible unless $B = \sum B_i$.

Let the second marked point of the bead $B_i$ lie in the
fixed component $F_i$.  
Since $B_i$ does not lie in
a single fixed component,  $K(F_i) \neq K(F_{i-1})$.
Obviously,
$F_{i-1}$ and $F_i$ cannot be joined by a bead 
of type $(p,q)$, where $q > 0$,
unless they lie in the same component of $M^{\Z/(q)}.$
Hence, by Lemma~\ref{le:bead}, 
$\omega(B_i) \geq \frac{|K(F_{i-1}) - K(F_{i})|}{q(F_{i-1},F_i)}.$
If $K(F_{i-1}) - K(F_i) >0 $, then equality is impossible unless $p = -1$,
in which case 
$-c_1(B_i) = \frac{m(F_{i-1}) - m(F_{i})}{q(F_{i-1},F_i)}$.
If $K(F_{i-1}) - K(F_i) <0 $, then equality is impossible unless $p = 1$,
in which case 
$-c_1(B_i) = -\frac{m(F_{i-1}) - m(F_{i})}{q(F_{i-1},F_i)}$.
Thus, 
$$
\sum_{i=1}^j \om(B_i) \geq
\sum_{i=1}^j \frac{|K(F_{i-1} - K(F_i)|}{q(F_{i-1},F_i)},
$$
with equality impossible 
unless 
$$
m_{\max} = -\sum_{i = 1}^{j}  c_1(B_i) =  \sum_{i = 1}^{j}  \frac
{m(F_{i-1}) - m(F_i)}
{q(F_{i-1},F_i)} \cdot
\frac {K(F_{i-1}) - K(F_i)} {|K(F_{i-1}) - K(F_i)|}.
$$
\QED\MS

One can formulate analogous results for the intermediate fixed components $F$.  
Consider the function $R:\Cc\to \Cc$, where $\Cc\subset \R$ is the set of 
critical values of $K$.  For $\mu\in \Cc$ we define $R(\mu)$ to be the
infinimum of the set of $\mu'$ such that
$F^\mu H_*(M) \subset  F_{\mu'}H_*(M)$. (For notation, see the discussion
after Remark~\ref{rmk:torsim}.)  In other words, every class 
$c^-$, 
where $c\in H_*(F)$ for some $F$ with $K(F)\le \mu,$
 is a linear combination of classes
$(c')^+$ where $c'\in H_*(F')$ for some $F'$ with $ K(F')\ge R(\mu)$.

\begin{prop}\labell{prop:xdist} 
 Suppose that $\La_K$ is inessential. 
Then for every critical value $\mu\in \Cc$, there is an invariant chain in some class $\si_z+B$  with $u_\La(\si_z+B) = c^{\ver}(\si_z+B)=0$  from a critical point $x$ with $K(x) \le \mu$ to 
another critical point $y$ with $K(y)\le R(\mu)$.
\end{prop}
\begin{proof}{}  If $c^-\in F^\mu H_*(M)$ then $\Ss(\La_K)(c^-) = c^-$ is represented in $F_{R(\mu)}H_*(M)$.  Hence $\Ss(\La_K)(c^-)\cdot_M c'\ne 0$ for some $c'\in F^{R(\mu)} H_*(M)$.  By Proposition~\ref{prop:invar} this is possible only if 
there is a chain in class $\si_z+B$ with the given properties.\end{proof}

\subsection{The semifree case}\labell{ss:semifree}

Suppose that the moment map $K$ generates a semifree $S^1$-action.
By Lemma~\ref{le:KJreg},
for a  generic almost complex structure 
$J \in \Jj_S(M)$,
the pair $(K,g_J)$ is Morse regular 
where $g_J$ denotes the metric associated to $J$
(see Definition \ref{Morseregular}). 
Let $F$ and $F'$ be  fixed components.
By Lemma~\ref{sfpseudo}, for any generic submanifolds $C$ of  $F$  
and $C'$ of $F'$, the unstable manifolds $W^u(C)$ and $W^u(C')$ 
are pseudocycles.  
We shall denote them by $W_J^u(C)$ and $W_J^u(C')$
to emphasize that they depend on the choice of $J$.
By construction, these
unstable manifolds are $S^1$-invariant. 
Hence, to prove Theorem~\ref{thm:semifree}, we only need
analyze the invariant chains in the moduli space
$\oMm_{0,2}(P_\La,\TJ,\si_F + B; W_J^u(C), W_J^u(C'))$, 
where $\omega(B) \leq 0$.

\begin{lemma}\labell{le:semifree} 
Consider a semifree Hamiltonian circle action $\Lambda_K$ on a compact symplectic
manifold $(M,\omega)$.
Let $\TJ$ be a generic almost complex
structure  in  $\Jj_S(P).$
Let $F$ and $F'$ be connected components of the 
fixed point set and let
$C \subset F$ and $C' \subset F'$ be generic submanifolds.
Fix $B \in H_2^S(M)$ such that $\omega(B) \leq 0$, and
consider the moduli space
$$
\oMm_{0,2}(P_\La,\TJ,\si_F + B; W_J^u(C), W_J^u(C')).
$$
\NI
{\rm (i)} If $B \neq 0$, the moduli space contains no  invariant chains.
\smallskip

\NI
{\rm (ii)} If $F \neq F'$, 
there are no  invariant chains unless 
$\dim(W_J^s(C)) + \dim(W_J^u(C')) > \dim M$.
\smallskip

\NI
{\rm (iii)} If $B = 0$ and  $F = F'$, the only  invariant chains are the constant
sections $\sigma_x$ for $x \in C \cap C'$. 
\end{lemma}

\begin{proof}\, 
Assume that there is an invariant chain from $x \in \ov{W_J^u(C)}$ to  
$y \in \ov{W_J^u(C')}$ with root $z$ in the class $\sigma_F + B$.
Note immediately that
$$ 
K(x)\; \leq\; K(F),
$$
with equality if and only if $x \in C$.
Let $A'$ and $A''$ denote the classes represented by 
the invariant subchains from $x$ to $z$, and from $z$ to $x$, respectively.
Then $A' + A'' + \sigma_z =  \sigma_F + B$, and so
by Lemma \ref{le:Ax}
$$ 
\omega(A') + \omega(A'') - K(z) + K(F) \;=\; \omega(B)\; \leq\; 0.
$$
Because the action is semifree, every bead of type $(p,q)$
in the invariant chain from $x$ to $z$ has $q = 1$.
Hence, by Lemma \ref{le:bead} 
$$
K(z) - K(x)\;\le\; \om(A')
$$
with equality impossible unless 
$K(x)  \le K(z)$,  $A'$ is the class of a chain of $(1,1)$ beads from $x$ to $z$,
and $A' = \sigma_x - \sigma_z$.
This implies both that there is a broken $K$-trajectory from 
$z$ to $x$,
and that
$$
0\;\le\; \omega(A''),
$$ 
with equality if and only if $A'' = 0$.
In this case, $z = y \in \ov{W_J^u(C')}$.
Therefore, $K(z) \leq K(F')$.

Considering all four displayed inequalities together, it is clear that
in fact they must all be equalities.  
This implies
that $A' = \sigma_x - \sigma_z$ and  $A'' = 0$,
and also that $x \in C \subset F$, so
$\sigma_x = \sigma_F$.  
Therefore $B = A' +A''+ \si_z - \si_F = 0$.
This proves (i).

Next, since it implies both that  $x \in C$ and that there is is a broken
$K$-trajectory from $z$ to $x$, there is
a broken $K$-trajectory from $z$ to $C$.  
Additionally,  since $z \in \ov{W_J^u(C')}$, by Lemma~\ref{closurebroken},
there  is a  broken $K$-trajectory from $C'$ to $z$.
Therefore, there is a broken $K$-trajectory from $C'$ to $C$.
If $F \neq F'$, then by Lemma~\ref{nobroken},  
this implies that $\dim W_J^s(C) + \dim W_J^u(C') > \dim M $.  This proves (ii).  

Finally, assume that $F = F'$. 
Then  since $K(x) = K(F)$,
$K(x)\le K(z)$,
and $K(z) \leq K(F')$,
it follows that $K(x) = K(z)$. Thus
$z\in F$ and
$A' = \sigma_F - \sigma_z = 0$. 
Since also $A'' = 0$, the last claim follows. \end{proof}

Before proving the rest of the theorems from the first section,
we need to consider the contributions of fixed point sets
other than $\Fmax$.
To simplify the proof of Lemma~\ref{le:transneg} below,
it is convenient to work with
almost complex structures on $M$ that are well behaved near the fixed components.
Each fixed component $F$ has a neighborhood
$\Nn_F$ that can be identified with a
neighborhood of the zero section in a
sum of Hermitian vector bundles
$\pi_F:E_1\oplus\dots\oplus E_k\to F$ in such a way
that the moment map $K$ is given by
\begin{equation}\label{eq:K}
K(v_1,\dots,v_k) = \sum_j\; \pi m_j\|v_j\|^2,\quad m_j\in 
\Z \oursetminus\{0\}
\end{equation}
and $S^1$ acts in  $E_j$ by rotation by $e^{2\pi i m_j}$.
The symplectic connection with horizontal spaces $\Hor_x$ equal to the
$\om$-orthogonals to the fibers is also $S^1$-invariant.
Therefore, starting from any $\om$-compatible $J_F$ on the components $F$, we 
may extend $J_F$ to an $S^1$-invariant
$\om$-compatible almost complex structure $J_M$ on $M$ whose restriction to
each set 
$\Nn_F$  agrees with
the complex structure on the fibers of $\pi_F$, leaves the horizontal
distribution invariant and is such that $\pi_F$ is holomorphic.
\begin{defn}\label{def:JSn} 
Fix once and for all such an almost complex structure  $J_M$ on $M$. Define
$\Jj_S^n(M)$ to be the set of all 
$S^1$-invariant
$\om$-compatible almost complex structures on $M$ that equal $J_M$ 
near the fixed point components $F$.   
Let $\Jj_S^n(P_{\La})$ denote the subspace of 
$\Jj_S(P_{\La})$ constructed from $J\in \Jj_S^n(M)$ as in 
Definition~\ref{def:JSn}.
\end{defn}

Thus when $J\in \Jj_S^n(M)$ 
each fixed point component $F$ has a neighborhood $\Nn_F$
 that can be identified with a neighborhood of the
 zero section in the complex vector bundle $\pi_F: E^+\oplus E^-\to F$, 
where $E^+$ (resp. $E^-$) is the subbundle of the normal bundle
with positive (resp. negative) weights.  Moreover, $\pi_F$ is $J$-holomorphic.
Hence  a neighborhood of the
submanifold $S^2\times F$ in $P_\La$ can be identified with a neighborhood of
the zero section in 
$$
\Tilde{\pi}_F: \Tilde E^+ \oplus \Tilde E^-\to S^2\times F,
$$
where the bundle $\Tilde E^\pm\to S^2\times F$ is induced in the obvious
way from the $S^1$-action.
Moreover, $\Tilde\pi_F$ is 
$\TJ$-holomorphic.
Denote by $\Mm^{\sect}_{0,2}(P_\La,\TJ,\si_F)$ the moduli space of 
$\TJ$-holomorphic maps $u:S^2\to P$ in class $\si_F$ and 
parametrized as sections
of $P_\La\to S^2$, and by $\Ff\subset \Mm^{\sect}_{0,2}$ 
the subspace of constant sections.

\begin{lemma}\labell{le:transneg} 
Fix $\TJ\in \Jj_S^n(P)$ and a
fixed point component $F$.  Then
the evaluation map
$$
\ev: \Mm^{\sect}_{0,2}(P_\La,\TJ,\si_F)\to M_0\times M_{\infty}
$$
is transverse to 
$(\Nn_F\cap E^-)\times (\Nn_F\cap E^-)\subset M_0\times M_{\infty}$
at all constant maps $u_x\in \Ff$.  Moreover, if all the positive 
weights at $F$ are $+1$, $\ev$ is a local diffeomorphism 
onto $(\Nn_F\cap E^+)\times (\Nn_F\cap E^+).$
\end{lemma}

\begin{proof}\,  
 The definitions imply that   $\Ff$ has a neighborhood $\Nn(\Ff)$ 
 consisting of all holomorphic maps $\tu: S^2\to \Tilde E^+\oplus \Tilde E^-$
whose composite with the projection
$\Tilde E^+\oplus \Tilde E^-\to S^2\times F$
is a holomorphic section of $\pi: S^2\times F\to S^2$ in the class $[S^2\times pt]$.
Since $S^2\times F\subset P_\La$ has the product complex structure, $\tu$ must 
project to some sphere $S^2\times \{x\}$ for $x\in F$.  
Therefore $\tu$ is a holomorphic 
section of the bundle $\Tilde E^+\oplus \Tilde E^-|_{S^2\times \{x\}}$.
But this  bundle is 
a sum of line bundles whose Chern classes 
are the nonzero weights of the $S^1$ action at $F$.  The line bundles with 
negative Chern classes have no sections, but those with positive 
Chern classes have plenty.   
Moreover if  all the positive weights  are $+1$ then there is precisely 
one section of $\Tilde E^+|_{S^2\times \{x\}}$
through any pair of points lying in distinct fibers.
The result follows.
 \end{proof}

\NI
{\bf Proof of Theorem~\ref{thm:semifree}.}
By Lemma~\ref{le:Ax}, 
$u_\Lambda(\sigma_F) = -K(F)$ and $c_{\ver}(\sigma_F) = m_F$.
Therefore we may write
$$
\Ss(\Lambda_K)(c^-) =  \sum_{B\in 
H_2^S(M)} a_B\otimes q^{-m(F)-c_1(B)}\;t^{K(F)-\om(B)},
$$
where $a_B$ is the contribution from the section $\sigma_F + B$.
Let $F'$ be a fixed component, and consider $c' \in H_*(F')$.
Fix $B \in H_2^S(M)$ so that $\omega(B)  \leq 0$.
We want to show that $a_B = 0$ unless
$B = 0$ and $F' = F$, 
in which case $a_B = c^+$.  Since the classes $(c')^-$ form a basis for
$H_*(M)$ it
is enough to show that $a_B \cdot_M  (c')^- = 0$ unless $B = 0$ and $F' = F$, 
in which case
$a_B \cdot_M (c')^- = c \cdot_F c'$.

Choose a generic almost complex structure $J \in \Jj_S^n(M)$.
By Lemma \ref{le:KJreg}, the pair $(K,g_J)$ is Morse regular,
where $g_J$ is the metric associated to $J$.
We may assume without loss of generality that $c$ and $c'$ can
can be represented by  generic submanifolds $C \subset F$ and $C' \subset F'$.
By Lemma~\ref{sfpseudo}, 
the unstable manifolds  $W_J^u(C)$ and $W_J^u(C')$ are pseudocycles.
Moreover, by Proposition~\ref{prop:equivar}, 
$[W_J^u(C)] = c^-$, and $[W_J^u(C')] = (c')^-$.
 Let $\TJ \in \Jj_S^n(P)$ be the associated almost complex
structure on $P_\Lambda$.

Assume first that $B \neq 0$.
By Lemma \ref{le:semifree}, 
the moduli space  $\oMm_{0,2}(P_\La,\TJ,\si_F + B; W_J^u(C), W_J^u(C'))$
contains no invariant chains. 
By Proposition~\ref{prop:invar}  
and Lemma~\ref{le:invariants}, 
this implies that $a_B \cdot_M (c')^- = 0$.

Now suppose that $B = 0$ but $F \neq F'$.
If $a_0 \ne 0$,
then  $\deg(a_0) = \deg(c^+) = \dim W_J^s(C)$.
Since $\deg((c')^-) = \dim W_J^u(C')$, we see immediately
that $a_0 \cdot_M (c')^- = 0$ for dimensional reasons unless
$\dim W_J^u(C) + \dim W_J^s(C') = \dim M$.  
However, if $\dim W_J^u(C) + \dim W_J^s(C') = \dim M$, then
by Lemma \ref{le:semifree}, 
the moduli space  $\oMm_{0,2}(P_\La,\TJ,\si_F; W_J^u(C), W_J^u(C'))$
contains no invariant chains. 
By Proposition~\ref{prop:invar}
and Lemma~\ref{le:invariants}, 
this implies that $a_B \cdot_M (c')^- = 0$.
 
Finally,
assume that $B = 0$ and $F = F'$.
By Lemma \ref{le:semifree},  
the space  
$$
\oMm_{0,2}(P_\La,\TJ,\si_F ; W_J^u(C), W_J^u(C'))
$$
contains no invariant chains 
except the constant sections
$\sigma_x$ for $x \in C \cap C'$.
By Proposition~\ref{prop:invar}
and Lemma~\ref{le:invariants},
this implies that only these
elements contribute to $a_0 \cdot_M (c')^-$. 
Now consider the full moduli space $\oMm_{0,2}(P_\La,\TJ,\si_F)$.  
It follows from Lemma~\ref{le:transneg} that the 
evaluation map
$$
\ev: \oMm_{0,2}(P_\La,\TJ,\si_F) \to M_0\times M_\infty
$$
intersects $W_J^u(C)\times W_J^u(C')$ transversally in $c \cdot_F c'$  
points. Hence $\GW^{P}_{\si_F,2}(c,c') = c \cdot_F c'$, as required.
 
\subsection{The case of at most twofold isotropy}\labell{ss:twofold}

When the action is not semifree, 
the unstable manifolds given by an $S^1$-invariant metric
may not be pseudocycles.  Therefore, we shall need
to consider more general objects.

\begin{defn}\label{def:Zc}
Let $C$ be a submanifold of a fixed component
$F$ with index $\alpha_F$.
A {\bf downwards pseudocycle from $\bf C$}
is an $S^1$-invariant 
weighted 
pseudocycle $f:Z_C^-\to M$, or $Z_C^-$ for short, 
of dimension
$\dim C + \al_F$ such that 
$f(Z_C^-)$ lies in $M^{K(F)}$, 
$f(Z_C^-) \cap K\inv(K(F)) = C$, and $[Z_C^-] = [C]^-$. Here, 
$[C]^-\in H_*(M)$ is the downwards extension of $[C]\in H_*(F)$ 
constructed in Section~\ref{sss:canon}.
\end{defn}

We show in Lemma~\ref{genpseudo} 
and Proposition \ref{prop:equivar}  
that given any generic submanifold
$C \subset F$, there exists a downwards pseudocycle $Z^-_C$.

As in the semifree case, 
we investigate $T^2$-invariant elements of  the cut moduli space.

\begin{lemma}\label{le:2isoa} 
Consider a Hamiltonian circle action $\Lambda_K$ 
on a compact symplectic
manifold $(M,\omega)$ with at most twofold isotropy. 
Let $\TJ$ be a generic almost complex structure in 
$\Jj_S(P_\La)$.
Let $F$ and $F'$ be (not necessarily distinct) fixed point components, 
and let $C \subset F$ and $C' \subset F'$ be generic submanifolds.
Fix a section class $\sigma$ 
and consider the moduli space
$$
\oMm_{0,2}(P_\La,\TJ,\sigma; Z_C^-, Z_{C'}^-).
$$
If $u_\Lambda(\sigma) \leq -\frac{1}{2}(K(F) + K(F'))$,
then the moduli space contains no invariant chains unless 
$\sigma = \frac{1}{2}(\sigma_F + \sigma_{F'})$ and
$F$ and $F'$ lie in the
same component of $M^{\Z/(2)}$.
\end{lemma} 

\begin{proof}\,
Assume that there is  an invariant chain from $x \in \ov{Z}_C$ to 
$y \in \ov{Z}_{C'}$ with root $z$  in the class 
$\sigma$.
We see immediately that
$$
K(x) \;\leq \;K(F) \quad \mbox{and}  \quad K(y)\; \leq\; K(F'),
$$
with equality if and only if $x \in C \subset F $ and $y \in C' \subset F'$.

Let $A'$ and $A''$ denote the  classes represented by 
the invariant subchains from $x$ to $z$, and from $z$ to $y$, respectively.
Since 
$A' + A'' + \sigma_z  \;=\; \sigma$,
by Lemma \ref{le:Ax} 
$$
\omega(A') + \omega(A'') - K(z) \; \leq  \;  -\frac{1}{2} (K(F) +  K(F'))
$$

Since the action has at most twofold isotropy, 
every bead of type $(p,q)$ in the invariant chain from $x$ to $z$ 
has $q \leq 2$.
Hence, by Lemma~\ref{le:bead} 
$$
\frac 12 \bigl(K(z) - K(x)\bigr)\;\le\; \om(A'),
$$
with equality impossible unless   $A'$ is
the class of a chain of $(1,2)$ beads from $x$ to $z$
and hence $A' = \frac 12 (\sigma_z - \sigma_x)$.
In particular, 
in this case
$x$ and $z$ lie in the same component of $M^{\Z/(2)}$.
By similar reasoning,
$$
\frac 12 \bigl(K(z) -  K(y)\bigr)\;\le \;\om(A''),
$$ 
with equality impossible unless  
$A'' = \frac 12 (\sigma_y -\sigma_z)$ and $y$ and $z$ 
lie in the same component of
$M^{\Z/(2)}$.

Considering all five displayed inequalities together, 
it is clear that they must
all be equalities.  
First, this means that $A' = \frac{1}{2}( \sigma_x - \sigma_z)$,  
$A'' = \frac{1}{2} (\sigma_y - \sigma_z)$,
$\sigma_x = \sigma_F$, and  $\sigma_y = \sigma_{F'}$; hence, 
$\sigma = \frac{1}{2}(\sigma_F + \sigma_{F'}).$
Second, it
implies that $F$ and $F'$ lie in the same component of $M^{\Z/(2)}$. 
\end{proof}

\NI
{\bf Proof of Theorem \ref{thm:2iso}.}
Choose a generic almost complex structure $J \in \Jj_S(M)$.
Let $\TJ \in \Jj_S(P_\La)$ be the associated almost complex
structure on $P_\Lambda$.
We may assume without loss of generality that $c$ and $c'$ can
be represented by  generic submanifolds $C \subset F$ and 
$C' \subset F'$, respectively.
By Lemma~\ref{genpseudo}
and Proposition \ref{prop:equivar}  we can find 
downwards  pseudocycles $Z_C^-$ and $Z_{C'}^-$ from $C$ and $C'$
as described above.
Fix a section class $\sigma$ so that
$c_{\ver}(\si) = u_\La(\si)=0$.  
If $K(F') \leq -K(F)$, then
$0 = u_\Lambda(\sigma) \leq - \frac 12 (K(F) + K(F'))$.
so by  Lemma~\ref{le:2isoa}, 
the moduli space 
$\oMm_{0,2}(P_\La,\TJ,\si; Z_C^-, Z_{C'}^-)$
contains no invariant chains
unless $K(F) = -K(F')$, $m(F) = -m(F')$ and $F$ and $F'$ lie in the same
component of $M^{\Z/(2)}$.
The result now follows
from  Proposition~\ref{prop:invar} and Lemma~\ref{le:invariants}.
\QED\MS

The previous result concerns the term $a_{0,0}\otimes 1$ in 
$\Ss(\La)(a)$ and hence gives information only in cases when we know that 
$a_{0,0}\ne 0$, for example if $\La$ is inessential.  
We next investigate the contribution from homologically visible components $F$ for 
general $\La$.  Again our arguments work only if the isotropy 
at levels above $F$ is at 
most twofold.

\begin{lemma}\label{le:2iso}  
Consider a  circle action $\Lambda_K$ on a compact symplectic
manifold $(M,\omega)$ with moment map $K : M \to \R$.
Let $\TJ$ be a generic almost complex
structure  in $\Jj_S(P).$
Let $F$ and $F'$ be connected components of the fixed point
set
and let $C \subset F$ and $C' \subset F'$ be generic submanifolds.
Assume $K(F') \leq K(F)$, that every positive weight at $F$
is $+1$,  and also that the isotropy for points
$w$ with $K(w)> K(F)$  is at most twofold.
Consider the moduli space
$$
\oMm_{0,2}(P_\La,\TJ,\si_F + B; Z_C^-, Z_{C'}^-)
$$
where $\om(B)  = 0$.
\begin{itemize}
\item 
If $B\ne 0$ or 
$F \neq F'$, there are no invariant chains in the moduli space.
\item 
If $B=0$ and 
$F = F'$, the only invariant chains
are the constant sections $u_{x}$ for $x\in C \cap C'$.
\end{itemize}
\end{lemma}

\begin{proof}\,
Assume that there is an invariant chain 
in class $\si_F + B$
from $x \in \ov{Z^-_C}$ 
to $y \in \ov{Z_{C'}^-}$
with root $z$.  
We see immediately that
$$ 
K(x) \;\leq\; K(F) \quad \mbox{and}  \quad K(y) \;\leq\; K(F'),
$$
with equality if and only if $x \in C \subset F$ and $y \in C \subset F'$.

Let $A'$ and $A''$ denote the  classes represented by 
the invariant subchains from $x$ to $z$, and from $z$ to $y$, respectively.
Since $A' + A''+ \sigma_z =   
\sigma_F + B$,
by Lemma \ref{le:Ax} 
$$
\omega(A') + \omega(A'') - K(z) +  K(F)\;  =\; \om(B)\; =\; 0.
$$
Since the isotropy for points $w$ with $K(w) > K(F)$
is at most twofold, every $(p,q)$ bead in the part of the invariant
chain from $x$ to $z$ which lies at least partially above $K(F)$ has
$q \leq 2$.
Hence, by Lemma \ref{le:bead},
$$ 
\frac 12 \Bigl(K(z) - \max\{K(F),K(x)\}\Bigr)\;\le\; \om(A'), 
$$
with equality impossible unless 
 $K(z) \geq K(F)$, 
$A' = \frac{1}{2}( \sigma_x - \sigma_z)$, 
and  $x$ and $z$ lie in the same component of $M^{\Z/(2)}$.
A similar reasoning applies to $A''$. 
Hence,
$$
\frac 12 \Bigl(K(z) -  \max \{K(F),K(y)\}\Bigr)\;\le\; \om(A''),
$$ 
with equality impossible unless 
$K(z) \geq K(F)$,
$A'' = \frac{1}{2}( \sigma_y - \sigma_z)$, 
and  $y$ and $z$ lie in the same component of $M^{\Z/(2)}$.

Considering all five displayed equations together  
with the hypothesis $K(F') \leq K(F)$, 
it is clear that they must
all be equalities.  

This implies that $x \in F$ and $y \in F',$ 
that $K(z) \geq K(F)$, that $K(y) = K(F)$,
and that $x, y$ and $z$
lie in the same connected component of $M^{\Z/(2)}$.  
Since all the positive weights at $F$ are $+1$, this implies that in fact 
$z \in F$ and $y \in F$, so that $\sigma_x = \sigma_y = \sigma_z = \sigma_F$.
Since $A' = \frac{1}{2}(\sigma_x - \sigma_z)$ and $A'' = \frac{1}{2}(\sigma_y - \sigma_z)$,
this implies that $B = 0$, and the result follows.
\end{proof}

\NI
{\bf Proof of Theorem \ref{thm:pfasimplemax}.}
Let $F'$ be a fixed component with $K(F') \leq K(F)$ and 
consider $c' \in H_*(F')$.
Choose a generic almost complex structure $J \in \Jj_S^n(M)$.
Let $\TJ \in \Jj_S^n(P)$ be the associated almost complex
structure on $P_\Lambda$.
We may assume without loss of generality that $c$ and $c'$ can
be represented by  
generic transversally intersecting 
submanifolds $C \subset F$ and 
$C' \subset F'$, respectively.
By Lemma~\ref{genpseudo} and Proposition~\ref{prop:equivar}. 
we can find downward pseudocycles 
$Z_C^-$ and $Z_{C'}^-$ as in Definition~\ref{def:Zc}.
These are constructed to coincide with the unstable manifolds 
$W^u(C)$ and $W^u(C')$ near $F$.
Since $J \in \Jj_S^n(M)$ is normalized near $F$, these unstable 
manifolds agree with neighborhoods of the zero section in the 
restrictions of $E^-\to F$  to $C$ and $C'$ respectively.

To show that 
$c_{0,0} \in F_{K(F)} H_*(M)$,
it is enough to show that if $K(F') < K(F)$, then $c_{0,0} \cdot_M (c')^- = 0$ 
for all $c'$. 
By Proposition~\ref{prop:invar}
and Lemma~\ref{le:invariants},
it suffices
to show that 
$\om(B) = c_1(B) = 0$ then the moduli space
$\oMm_{0,2}(P_\La,\TJ,\si_F+ B; Z_C^-, Z_{C'}^-)$
contains no  invariant chains.
Since $F' \neq F$, this is immediate from 
Lemma~\ref{le:2iso}.  

Now suppose $F = F'$.
 Consider the moduli space $\oMm\,\!^{cut}= 
 \oMm_{0,2}(P_\La,\TJ,\si_F+ B; Z_C^-, Z_{C'}^-)$.  By 
  Lemma~\ref{le:2iso},  
  this is nonempty only if $B=0$.  Further, in this case
  the only invariant chains 
 in $\oMm\,\!^{cut}$ are the constant sections $u_x$ for $x \in C \cap C'$.
   These 
  sections form one of the connected components of  $\oMm\,\!^{cut}$, and by 
  Proposition~\ref{prop:invar} we may ignore any other components.  Thus we 
  may suppose that $\oMm\,\!^{cut}$ reduces to  
 the compact manifold $C \cap C'$.  
If any of the negative weights along $F$ are less than $-1$
the
  elements of $\oMm\,\!^{cut}$
  are not regular.  As in the proof of
  Proposition~\ref{prop:amax}, their cokernels fit together to form the
 obstruction bundle
  $\Ee\to\oMm\,\!^{cut}$ of Equation~(\ref{eq:ee}). 
  Because all positive weights along $F$ are $+1$ and because the 
  sets $\ov{Z_C^-}, \ov{Z_{C'}^-}$ 
  coincide  near $F$ with the 
  bundles $E^-|_C, E^-|_{C'}$, it follows from
  Lemma~\ref{le:transneg} that the full moduli space  intersects 
  $\ov{Z_C^-}\times \ov{Z_{C'}^-}$ transversally in 
  $\oMm\,\!^{cut}$.
  Hence standard theory implies  that
  the regularized cut down moduli 
  space represents the class  
$e(\Ee)\cdot_F [C\cap C'] = (e(\Ee)\cap_F c)\cdot_F c'$ and that
 $$
\GW^{P_\Lambda}_{\sigma_{F},2}(c^-,(c')^-)=  
(e(\Ee)\cap_F c)\cdot_F c'.
 $$
The result follows. 
 \QED

\section{Proofs of main technical lemmas}

We now establish
the main technical results used in the paper.

\subsection{Invariant cycles in $M$.}\labell{ss:Morse}

This section establishes the properties of the canonical extension
classes $c^{\pm}$ used in 
Theorems~\ref{thm:semifree}, \ref{thm:2iso}, and \ref{thm:pfasimplemax}. 
We first prove Lemma~\ref{le:unique0} 
which is used to construct
canonical downwards 
and upwards extensions of the homology  classes
of the fixed point set, and then show how to define
 representing cycles for these classes that have the properties 
 claimed in Definition~\ref{def:Zc}.

\subsubsection{Canonical classes}\labell{ss:canclasses}

Let $S^1$ act on a symplectic manifold $(M,\omega)$,
with a moment map $K$ which is proper and bounded below.
Then $K$ is  Morse--Bott function with extraordinary properties.
(For background information see~\cite{TW}.)
First,  $K$ is equivariantly perfect, that is, the restriction
map $H_{S^1}^*(M) \to H_{S^1}(M^{< \mu})$ is surjective for
all $\mu \in \R$, where $M^{<\mu} := K\inv(-\infty,\mu)$.
The same proof shows that the restriction to the fixed
point set is injective.  More specifically,
given any $\TY \in H_{S^1}(M)$, then
$\TY|_{M^{<\mu}} = 0$ if and only if 
$\TY|_{F'} = 0$ for all fixed components $F'$ with $K(F') < \mu$.
The same argument also
shows that $H^*_{S^1}(M)$ is equivariantly formal, that is,
the restriction $H^*_{S^1}(M) \to H^*(M)$ is surjective.
\MS

\NI
{\bf Proof of Lemma~\ref{le:unique0}.}
Let $S^1$ act on a compact symplectic manifold $(M,\om)$ with moment map $K$.
Let $F \subset M$ be any fixed component of index $\alpha$;  
and let $e^-_F\in H_{S^1}^\alpha(F)$ be the equivariant 
Euler class of the negative normal bundle to $F$.
Given any cohomology class $Y \in H^i(F)$, 
we must show that
there exists a unique 
cohomology class $\TY^+ \in H_{S^1}^{i + \alpha}(M)$ so that
\begin{description}
\item[(a)] the restriction of $\TY^+$ to $M^{<K(F)}$ vanishes, 
\item[(b)] $\TY^+|_F = Y \cup e^-_F$,  and 
\item[(c)] the degree of $\TY^+|_{F'}$ in $H_{S^1}^*(pt)$ is less than the 
index $\alpha_{F'}$  of $F'$ for all fixed components $F' \neq F$.
\end{description}
Moreover, we claim that
 these classes generate $H^*_{S^1}(M)$ as a $H^*_{S^1}(pt)$ module.

Since $K$ is equivariantly perfect,
we can find $\TY^+$ satisfying (a) and (b).
In fact, in general there will be many such $\TY^+$.

Enumerate the fixed sets other than $F$ by $F_1,\ldots,F_k$ so that 
$K(F_j) \leq K(F_{j+1})$ for all $j$.
Assume that $\TY^+$ satisfies (c) for all $F_j$ such that $j < i$.
Let $\alpha_i$ denote the index of $F_i$,  
and let $m_-(F_i)$ denote the product of the negative weights at $F_i$.
Then $e^-_{F_i}$,
the equivariant Euler class of
the negative  normal bundle to $F_i$,
is equal to $e^-_{F_i} = m_-(F_i) \,u^{\al_i/2} +$ 
terms of degree $< \al_i$ in 
$H_{S^1}^*(pt)$, where  $m_-(F_i) \neq 0$.
Therefore, $\TY^+|_{F_i} \in  H_{S^1}^*(F_i)$ can be written uniquely as a
sum $\TX + \TX'$, where $\TX$ is a multiple of $e^-_{F_i}$ and the degree 
of $\TX'$ in 
$H_{S^1}^*(pt)$ is less than $\alpha_i$.
Since $K$ is equivariantly perfect, there exists $\TY' \in H^*_{S^1}(M)$ so
that $\TY'|_{F_i} = \TX$ and $\TY'_{F_j} = 0$ for all $j < i$. 
 After subtracting $\TY'$,  we find a new $\TY^+$
that  satisfies (c) for all $F_j$ such that $j \leq i$.

To see that $\TY^+$ is unique, let $\widehat{Y}$ be the difference of
two classes that satisfy (a), (b), and (c).
Then  the degree of $\widehat{Y}|_{F'}$ in $H_{S^1}^*(pt)$ is less than the 
index of $F'$ for
{\em every} fixed component $F'$
(and $\widehat{Y}|_F = 0$).
Let $F_j$ be the smallest $j$ such that
$\widehat{Y}|_{F_j} \neq 0$.
Then, since the restriction to the fixed point set is injective, 
$\widehat{Y}$ vanishes when restricted to $H_{S^1}^*(M^{< K(F_j)}).$ 
Hence, $\widehat{Y}|_{F_j}$ is a multiple of $e_-(F^j)$.  But this is impossible,
so $\widehat{Y}|_{F_i} = 0$ for all $i$. Hence, $\widehat{Y} = 0$.

Finally, for any $Y \in H^*(F)$, since $\TY^+|_{M^{< K(F)}} = 0$
the restriction of $\TY^+$ to $M_{K(F)}$ 
is an element of 
$H^*_{S^1}(M^{K(F)},M^{< K(F)})$.
By injectivity, as $Y$ ranges over $H^*(F)$,
these classes generate
$H^*_{S^1}(M^{K(F)},M^{< K(F)})$ as a
a $H^*_{S^1}(pt)$  module.  Hence,
if we also let  $F$ vary over all fixed components, 
then they generate $H^*_{S^1}(M)$. 
\QED 

\subsubsection{Morse cycles and equivariant cohomology}
\labell{pseudosec}

We shall work throughout with pseudocycles, and
begin by recalling their definition from~\cite{MS2}.
A {\bf pseudocycle}
of dimension $d$ in a manifold $M$ 
is a smooth map $f:V \to M$ from an oriented smooth $d$-dimensional 
manifold $V$ to $M$ whose 
$\Om$-limit set 
$$
V^{\infty} : = \{x\in M\;:\; x=\lim_{i \rightarrow \infty}  f(y_i),
\mbox{ where $\{y_i\}_{i=1}^\infty$ has no 
limit point in } V\}
$$
has codimension at least $2$, i.e. it is in the 
image of a smooth map $g:W^{d-2}\to M$.  Two pseudocycles 
$f_i: V_i\to M$ are bordant
if they  can be extended over a manifold $W$  with boundary $V_1\cup 
-V_2$  by a map whose $\Om$-limit set has dimension at most 
$d-1$.   Any $f:V^d\to M$ is bordant to a map that intersects a given 
codimension $d$ submanifold $X$ of $M$ transversally, i.e. $X\cap 
V^{\infty} = \emptyset$ and $f:V\to M$ meets $X$ transversally. 
Moreover, because the boundary has codimension at least $2$, 
each bordism class $[f,V]$ of pseudocycles 
defines  a unique 
rational homology class $c(f,V)$. (In fact, it defines a unique integral 
class: see Schwarz~\cite{Sch}.)
We say that two such cycles of 
complementary dimension  meet 
transversally  if the closures of their images  
$f(V)$ and $f'(V')$
intersect only along their top strata 
$f(V)$ and $f'(V')$
and if these 
intersections are transverse in the usual way.  It is shown in~\cite{MS2} 
that the intersection number $c(f,V)\cdot c(f',V')$ can be calculated 
by perturbing $(f,V)$ to be transverse to $(f',V')$ and then counting 
the points of intersection of $f$ with $f'$ in the usual way.

In this paper it is convenient to work with rational combinations of 
such pseudocycles.  Therefore we make the following definition.

\begin{defn}\label{def:psc} A {\bf weighted pseudocycle} is a finite sum 
    $\sum_iq_i f_i$, where  $q_i\in \Q$ and $f_i: Z_i\to M$ 
is a pseudocycle as above.  For short we sometimes forget the 
weights $q_i$ and denote this pseudocycle by $f:Z\to M$, where 
$Z: = \cup_i Z_i$ and $f|_{Z_i}: = f_i$.  We say that $(f,Z)$ 
is $S^1$-invariant iff 
the closure  $\ov{f(Z)}$ of the image  is
invariant under the $S^1$-action.  We also sometimes omit the map $f$ 
from the notation, denoting the cycle by $Z$ and the closure of its 
image by  $\ov{Z}$.
\end{defn}

All the cycles considered in this 
paper (except for the virtual moduli cycle)
are weighted pseudocycles.

Let $K$ be a Morse-Bott function, and let $g$ be any metric.
Consider the {\bf negative gradient flow}
$$
\psi : \R \times M \to M$$
such that 
$$ 
\frac{\partial}{\partial t} \psi(t,x) = -\grad f(\psi(t,x)) \quad \mbox{and} \quad
\psi(0,x) = x \quad \mbox{for all} \ x \in M, \; t \in \R.
$$
A {\bf gradient trajectory} is a map $\gamma: \R \to M$ such
that 
$$
\frac{d}{dt} \gamma(t) = - \grad f(\gamma(t)) \quad \mbox{for all } t \in \R.
$$
More generally, a {\bf broken gradient trajectory} is a set of gradient
trajectories $\gamma_1,\ldots,\gamma_n$ such that
$\lim_{t \rightarrow \infty} \gamma_i = \lim_{t \rightarrow -\infty} \gamma_{i+1}$
for all $i$. (By convention, we allow the case $n = 1$.)

Given any critical component  $F$, 
we define the {\bf stable manifold} and the {\bf unstable manifold}, respectively,
by 
$$ W^s(F) = \{ x \in M \mid \lim_{t \rightarrow \infty} \psi(t,x) \in 
F\},\quad 
W^u(F) = \{ x \in M \mid \lim_{t \rightarrow -\infty} \psi(t,x) \in F\}.
$$
Define  maps
$$
\pi_+: W^s(F) \to F \quad \mbox{and} \quad \pi_-:W^u(F) \to F
$$
 by
$$\pi_+(x) =  \lim_{t \rightarrow \infty} \psi(t,x) \quad \mbox{and} \quad 
\pi_-(x) =  \lim_{t \rightarrow -\infty} \psi(t,x).$$
Both $\pi_+$ and $\pi_-$ are submersions; see \cite{AB}.

Given a collection of distinct
critical components $F_1,\ldots,F_k$,
we define a natural map 
$$
f: W^u(F_1) \times \cdots \times W^u(F_{k-1}) \times W^s(F_2) \times \cdots 
W^s(F_k) \to (M^{k-1} \times F_2 \times \cdots \times F_{k-1})^2
$$
by 
\begin{eqnarray}\labell{eq:trans}
&&f(a_1,\ldots,a_{k-1},b_2,\ldots,b_k) = \\\notag 
&&\qquad (a_1,\ldots,a_{k-1},\pi_-(a_2),\ldots,\pi_-(a_{k-1}), 
b_2,\ldots,b_k,\pi_+(b_2),\ldots,\pi_+(b_{k-1})).
\end{eqnarray}
Define $\calM (F_1,\ldots,F_k) = f\inv(\Delta)$.
We can (and will) identify $\calM (F_1,\ldots, F_k)$  with tuples
$$
(x_1,\ldots,x_{k-1}) \subset M^{k-1}
$$
such that 
$x_i \in W^s(F_i) \cap  W^u(F_{i+1})$ and
$\pi^+(x_i) = \pi^-(x_{i+1})$ for all $i$.
More geometrically,
$\calM (F_1,\ldots,F_k)$  consists of all
tuples $(x_1,\ldots,x_{k-1})$ for which there is a
broken gradient trajectory $\ga_1,\ldots, \gamma_{k-1}$  from $F_1$ to $F_k$ through
$F_2, F_3, \ldots, F_{k-1}$  so that $\ga_i$ contains the point $x_i$ 
for all $i$.
Define maps
$$\pi_- : \calM (F_1,\ldots,F_k)  \to F_1 
\quad \mbox{and} \quad 
\pi_+ : \calM (F_1,\ldots,F_k)  \to F_k$$
by 
$$\pi_-(x_1,\ldots,x_{k-1}) = \pi_-(x_1) 
\quad \mbox{and} \quad 
\pi_+(x_1,\ldots,x_{k-1}) = \pi_+(x_{k-1}).$$

\begin{definition}\labell{Morseregular}
We say that the  pair $(K,g)$ is {\bf Morse regular} if $f$ is transversal
to the diagonal 
$$\Delta \subset (M^{k-1} \times F_2 \times \cdots \times F_{k-1})^2$$
for every collection of critical components $F_1,\ldots,F_k.$
\end{definition}

In general, this is stronger than assuming that $W^s(F)$ and 
$W^u(F')$ intersect transversally for all critical sets $F$ and $F'$, but it is 
equivalent if $K$ is a Morse function.

If $(K,g)$ is Morse regular, then
by transversality,  $\calM (F_1,\ldots,F_k)$  is a manifold
of dimension $f_1 + \alpha_1  - \alpha_k$,
where $f_i$ is the dimension of $F_i$ and $\alpha_i$ is the index of $F_i$.
Note that the reparametrization group $\R^{k-1}$ acts on
the elements $(\ga_1,\cdots,\ga_{k-1})$ in
$\calM (F_1,\ldots,F_k)$ so that the 
set of points in $M$ that lie on a broken
trajectory in $\calM (F_1,\ldots,F_k)$ has dimension $\le 
f_1 + \alpha_1  - \alpha_k - (k-1)$.

\begin{lemma}\labell{nobroken}
Let $M$ be a  compact manifold.
Let $K: M \to \R$ be a Morse-Bott function and let $g$ be a metric 
so that the pair $(K,g)$ is Morse regular.
Let $F$ and $F' $ be distinct critical components.
If $C \subset F$ and $C' \subset F'$ are generic submanifolds,
there is no  broken gradient trajectory 
from $C'$ to $C$ unless
$$
\dim  W^s(C) + \dim W^u(C') > \dim M.
$$
\end{lemma}

\begin{proof}{}
Assume that there is a broken trajectory from $F' = F_1$ to $F = F_k$
through critical components $F_2,\ldots,F_{k-1}$.
By genericity, we may assume that the
maps $\pi_-: \calM (F_1,\ldots, F_k,) \to F_1$ and
$\pi_+: \calM (F_1,\ldots, F_k) \to F_k$ are transverse to
$C'$ and $C$, respectively. 
Therefore, the set $X = C' \times_{\pi_-} \calM (F_1,\ldots,F_k) \times_{\pi_+}
C$ is a manifold of dimension $c' + \alpha' + c - \alpha-f$,
where $c',c,$ and $f$ denote the dimensions of $C',C,$ and $F$, and
$\alpha',\alpha$ denote the 
index of $F',F$ respectively.
There is a proper effective action of $\R$ on $\calM (F_1,\ldots,F_k)$
which moves $x_1$ along the gradient trajectory on which it lies;
This induces an action on $X$.
Hence, $X$ is empty unless 
$c' + \alpha' + c -\alpha -f> 0$.
Since $\dim W^u(C') = c' + \alpha'$, and $\dim W^s(C) = c -f - \alpha$,
$X$ is empty unless 
$\dim  W^s(C) + \dim W^u(C') > \dim M$ as claimed.
\end{proof}

We will also need the following lemma, which can be easily proved by
a slight variation of the proof for the analogous fact in
the Morse case.
 
\begin{lemma}\labell{closurebroken}
Let $M$ be a  compact manifold.
Let $K: M \to \R$ be a Morse-Bott function and let $g$ be a metric 
so that the pair $(K,g)$ is Morse regular.
Let $C$ be a submanifold of a critical component $F$.
Every point in $\overline{W^u(C)}$, the closure of the 
unstable manifold of $C$, lies on a 
broken trajectory beginning in $C$.
\end{lemma}

Since we want the unstable manifold $W^u(C)$ to be $S^1$-invariant,
we next investigate gradient flows with respect to invariant metrics.
In general, due to the presence of isotropy spheres,
there may be no $S^1$-invariant metric $g$ so that the pair $(K,g)$ 
is Morse regular, even if the moment map $K$ is Morse.
For example, consider 
the action $[z_0 \colon z_1 \colon z_2]\mapsto 
[e^{2\pi it}z_0 \colon z_1 \colon e^{-2\pi 
it}z_2]$ on $\CP^2$ and blow up the point $[0 \colon 1 \colon 0]$.
The exceptional divisor 
$\Si$ has isotropy group $\Z/(2)$ and contains two critical points, 
both of index $2$. Any $S^1$-invariant vector field must be tangent to
$\Si$ since if $\phi$ denotes the generator of the isotropy 
subgroup $d\phi$ acts as $-1$ in the directions normal to 
$\Si$.   In particular, the gradient $\grad K$ of $K$ with respect to an 
invariant metric  must be tangent to $\Si$ and hence 
have trajectories joining two critical 
points of equal index. The following lemma shows that this
is the only obstruction to finding a Morse regular pair 
$(K,g_J)$.  
Recall that  $\Jj_S(M)$ is the space of smooth invariant $\om$-compatible 
almost complex structures on $M$ that are normalized near the fixed 
point components $F$ as described in Definition~\ref{def:JS}.

\begin{lemma}\labell{le:KJreg} 
Let $S^1$ act semifreely  on a compact symplectic manifold $(M,\om)$
with moment map $K$.
For a  generic almost complex structure  
$J\in \Jj_S^n(M)$,
the pair  $(K,g_J)$ is Morse regular,
where $g_J$ is the metric associated to $J$.
 \end{lemma}

\begin{proof}\,  Salamon and Zehnder show in ~\cite[Theorem~8.1]{SZ}
that the gradient flow of any Morse function $H$  on $(M, \om)$
is Morse--Smale with respect to a generic metric of the form $g_J$, where   
$J$ ranges over the set of all $\om$-compatible almost complex 
structures.  
We simply need to check that their argument continues to 
hold for Morse-Bott functions in the presence of a semifree $S^1$-action.

Inspection of the proof of \cite[Theorem~8.1]{SZ} shows that
the map $f$ in equation~(\ref{eq:trans}) satisfies the
required transversality condition 
provided that the 
tangent space $T_J(\Jj_S)$ of the space $\Jj_S: = \Jj_S(M)$ 
of allowable $J$ is large enough. (See also
Austin--Braam~\cite[Proposition~B.2]{AB}.)
This tangent space $T_J(\Jj_S)$ is contained in
the space of
$S^1$-invariant
sections of the
bundle $\End$ of anti-$J$-holomorphic endomorphisms of $TM$ over $M$,
and we need each gradient flow line $\ga$  to go through
a point $x\in M$ such that
there are elements $Y\in T_J(\Jj_S)$ whose value $Y(x)$ is an
arbitrary element in $\End_x$
and whose support intersects $\ga$ in an arbitrarily small set. 
Since
the isotropy group of $x$ is trivial for all points on $\ga$ this
is clearly the case; $\ga$ is transverse to the level sets of
the moment map $K$  and there are elements in $T_J(\Jj_S)$ with support
in $K^{-1}(a,a+\eps)$ for arbitrarily small $\eps$ and arbitrary 
value at $x$.  If there were isotropy at $x$ then this argument would 
fail because $Y(x)$ would have to be fixed by
$d\phi$ for all $\phi$ in the isotropy group at $x$.
\end{proof}

The following lemma is  adapted from  Schwarz~\cite{Sch}.  
\begin{lemma}\labell{sfpseudo}
Let $S^1$ act semifreely on a compact manifold $M$.
Let $K: M \to \R$ be an $S^1$-invariant  Morse-Bott function and let $g$ be an
$S^1$-invariant  
metric so that the pair $(K,g)$ is Morse regular.
Given a generic submanifold $C$ of a fixed component $F$,
the unstable manifold $W^u(C)$ is a pseudocycle.
\end{lemma}

\begin{proof}{}
The unstable manifold  $W^u(C)$ is a submanifold of dimension $c + \alpha$,
where $c$ is the dimension of $C$ and $\alpha$ is the index of $F$.
Hence, we must show that $\ol{W^u(C)} \ssminus W^u(C)$ has dimension
at most $c + \alpha - 2$. 

Because $C$ is generic, we may assume that
$C$ is transverse to the map $\pi^- : \calM(F_1,\ldots,F_n) \to F_1$
for every sequence of critical points $F = F_1,F_2,\ldots,F_n$.
Therefore, 
$X = C \times_{\pi^-} \calM(F_1,\ldots,F_n) \subset \calM(F_1,\ldots,F_n)$ is
a manifolds of dimension $c + \alpha   - \alpha_n$,
where  $\alpha_n$ is the  index of $F_n$. 

There is a smooth proper action of 
$\R$ on $X$,
which moves the first coordinate $x_1$ along the gradient
trajectory on which it lies.  There is another smooth proper action of $S^1$ on
$X$, which is given by the circle action on $x_1$.
If $n > 1$,   
the evaluation map $\ev: X \to M$ 
defined  by $\ev(x_1,\ldots,x_n) = x_n$ is
constant along the orbits of
these actions.
Hence, the image of the evaluation map has dimension at most
$c + \alpha - 2$.

By Lemma \ref{closurebroken}, every point in the closure $\ov{W^u(C)}$ 
lies on a broken trajectory  that begins on  $C$,
that is, it lies in the image of the evaluation map for 
$X = C \times_{\pi^-} \calM(F_1,\ldots,F_n) \subset \calM(F_1,\ldots,F_n)$ 
for some sequence of fixed points $F = F_1,\ldots,F_n$.
Moreover, if the point does not lie $W^u(C)$ itself, then 
$n$ must be greater than one.
\end{proof}

\begin{lemma}\labell{genpseudo}
Let $S^1$ act on a compact manifold $M$.
Let $K: M \to \R$ be an $S^1$-invariant  Morse-Bott function.
Given a generic submanifold $C$  of a fixed component $F$
of index $\alpha_F$, there
exists an $S^1$-invariant 
weighted
pseudocycle $Z_C^-$ in $M^{K(F)}$ 
of dimension $\dim C + \alpha_F$
such that $Z_C^- \cap K\inv(K(F)) = C$.
\end{lemma}
\begin{proof}{}
In this case, as illustrated by the example
after Lemma~\ref{closurebroken}
there may be no $S^1$-invariant metric $g$ so that
the pair $(K,g)$ is Morse regular.  Instead, we begin with  any
$S^1$-invariant metric $g$, and then consider an $S^1$-invariant
multivalued perturbation.

Briefly, the idea is this.  
Consider the space of all $S^1$-invariant smooth
multivalued vector fields $Y$ on $M$.  We will suppose for simplicity 
that $Y$ is single valued everywhere except on a finite number of 
disjoint slices 
of the form $M^{\mu}\oursetminus M^{\mu-\eps}$ that contain no critical points of 
$K$, and that at each point $x$ in such a slice $Y(x)$ is a finite set 
that is invariant under the action of the isotropy group at $x$.  The 
smoothness condition means that the graph $\{(x,v); v\in Y(x)\}$ of 
$Y$ is a union of smoothly embedded open subsets of Euclidean space. 
For example, in the case of the 
blow up of $\CP^2$ discussed at the beginning of this section, we allow
$Y$ to take two values $\pm v(x)$ when $x\in M^{\mu}\oursetminus M^{\mu-\eps}$ 
is near the isotropy submanifold. It is easy to check that 
there are enough perturbations of this kind so that for generic 
small $Y$ each solution
$\ga:\R\to M$ of the corresponding perturbed gradient flow relation
\begin{equation}\labell{eq:gradreln}
\frac d{ds}\ga(s) \in\{ -(\grad_{g}K + Y)(\ga(s))\}
\end{equation}
is  transverse to the level sets $K=const$
and  regular in the sense of Salamon--Zehnder~\cite{SZ}. 
To keep the structure of the 
solution set as simple as possible we may assume that the number of 
elements in each set $Y(x)$ is constant and equal to $N$  for all $x$ lying 
in the interior of a slice, where $N$ is the l.c.m. of the 
orders of the stabilizer subgroups of the $S^1$-action.  
Then, a trajectory $\ga$ that goes from 
a point $\ga(-\infty)\in F$ to $\Fmin$ passes through some number  $k$ 
of slices and hence satisfies  one out of a set of $N^k$ 
possible equations. Moreover, 
because  the set of trajectories that do not reach $\Fmin$ lie in a 
closed subset of codimension at least $2$, there is a neighborhood $U$
of $\ga(-\infty)$ in $F$  such that
the set of trajectories that start in $U$ form a disjoint union of 
$N^k$ submanifolds.  Thus for each generic submanifold $C$ in $F$
the set $W^{u,Y}_C$ of solutions to~(\ref{eq:gradreln}) 
that start at $C$ and end 
in $\Fmin$ is a manifold.  As before, the transversality condition 
means that $W^{u,Y}_C$ intersects the corresponding stable manifolds 
transversally.  (These are solutions to the relation
$
\frac d{ds}\ga(s) \in\{(\grad_{g}K + Y)(\ga(s))\}$.)
Hence the previous arguments apply to show that $W^{u,Y}_C$  is a 
pseudocycle.  It is $S^1$-invariant by construction. Therefore we  
define
$
Z^-_C $ to be the weighted pseudocycle
$$
Z^-_C: = \frac 1{N^k} W^{u,Y}_C.
$$
This completes the construction.\end{proof}

 Repeating the above construction for $-K$ we obtain upwards
pseudocycles $Z_C^+$. 
It remains to prove that these extensions represent the canonical 
extensions $[C]^{\pm}.$

\begin{prop}\labell{prop:equivar}  
Let $S^1$ act on a compact symplectic manifold 
$(M,\om)$ with moment map $K: M \to \R$. 
Let $C$ be a  generic  submanifold  of a fixed component $F$.
If the action is semifree, let $g$ be an $S^1$-invariant metric so that
the pair $(K,g)$ is Morse regular. Then $[W^s(C)] = [C]^+.$
More generally, 
construct the 
weighted
pseudocycle  $Z_C^+$ as in Lemma~\ref{genpseudo}. 
 Then $[Z_C^+] = [C]^+$. 
 \end{prop}

 \begin{proof}{}  
Let $Z^+$ denote either the pseudocycle $W^s(C)$ or
the weighted pseudocycle $Z_C^+$, as appropriate.
Let $f$ be the dimension of $F$, $i$ be the dimension of $C$, and let $\alpha$
be the index of $F$.
Let $Y \in H^{f - i}(F)$ be the Poincar\'e dual to $C$.
Let ${Y^+} \in H^{f - i + \alpha}(M)$ denote the restriction to ordinary
cohomology of the unique equivariant cohomology class  
$\TY^+ \in H_{S^1}^{f - i + \alpha}(M)$ described in 
Lemma~\ref{le:unique0}.
Recall from Section~\ref{sss:canon} that the upwards extension 
$[C]^+$ is defined to
be the Poincar\'e dual of the restriction of $\TY^+$ to $H^{f-i+\alpha}(M)$.

Fix $N > d: = f - i + \alpha$, and note that 
$\dim Z^+ = \dim M-d$.
Since $Z^+$ is $S^1$-invariant, it
can be extended to a cycle $(Z^+)^N: = S^{2N+1}\times_{S^1} Z^+ $ 
in the finite dimensional 
approximation $
M_{S^1}^N: = S^{2N+1}\times_{S^1} M
$
to $M_{S^1}$.
Denote by 
$$
{\TX}\,\!^N \in H^{d}(M_{S^1}^N)
$$
the Poincar\'e dual of $(Z^+)^N$ in $M_{S^1}^N$.  
 Clearly, the restriction of ${\TX}\,\!^N$ to $M$ 
 is Poincar\'e dual to $[Z^+]$.
Therefore, it is enough to show that ${\TX}\,\!^N$ is  the restriction 
of $\TY^+$ to $M_{S^1}^N$.
By  the injectivity of the restriction maps
$$
H^*_{S^1}(M)\to H^*_{S^1}(M^{S^1}),\qquad  H^d_{S^1}(M)\to 
H^d_{S^1}(M_{S^1}^N),
$$
it is enough to show that the restriction  ${\TX}\,\!^N|_F$
of ${\TX}\,\!^N$  to $S^{2N+1} \times_{S^1}  F$ 
satisfies 
the conditions (a), (b), and (c) of Lemma~\ref{le:unique0}.

Because $Z^+$ 
has standard form near $C\subset F$, it 
is represented by the restriction over $C$ of
the positive normal bundle of $F$.  Therefore 
${\TX}\,\!^N|_F = Y \cup e^-_F$ as 
required by property (b) of Lemma~\ref{le:unique0}.
Clearly, ${\TX}\,\!^N|_{F'} = 0$ for all fixed components 
$F'$ such that $K(F') < K(F)$.
Therefore it suffices  to check that ${\TX}\,\!^N$ has property (c).
Let $F'$ be any fixed component other than $F$, and
let $\alpha'$ be the index of $F'$.
We wish to show that the degree of ${\TX}\,\!^N|_{F'}$ in
$H^*((BS^1)^N): = H^*(\CP^N)$ is less than $\alpha'$, or equivalently 
that the degree  of ${\TX}\,\!^N|_{F'}$ in $H^*(F')$ is
greater than $d -\alpha'$.
To prove this, it is enough to show that if 
$X \subset F'$
is a generic submanifold of dimension 
$d - \alpha'$,
then $(S^{2N+1}/S^1) \times X \subset  S^{2N+1} \times_{S^1} F$
does not meet $S^{2N+1} \times_{S^1}\ov{Z^+}$. Hence it suffices to 
check that $X$ does not meet $\ov{Z^+}$.

In the semifree case, $Z^+$ is the stable manifold 
$W^s(C)$ with respect to a generic metric $g_J$.  By 
Lemma~\ref{closurebroken} every element 
in the closure $\ov{W^s(C)}$ lies on a broken geodesic ending at $C$.
Therefore, by Lemma~\ref{nobroken},  $X\cap \ov{W^s(C)}\ne \emptyset$ only if 
$\dim X + \al'+ \dim W^s(C) > \dim M$.
Since by construction $\dim X + \al' + \dim W^s(C) 
= \dim M$, the
intersection is empty.

In the general case, $Z^+$ is the sum of pseudocycles
that are arbitrarily $C^0$-close to $\ov{W^s(C)}$.  Therefore, 
the argument above shows that
it can be constructed so that its closure  
$\ov{Z^+}$ is disjoint 
from any finite set of manifolds  
$X_i\subset F'$ 
that span the homology group $H_{d-\al'}(F')$.
The result follows.
\end{proof}
 
\subsection{Localization}\label{ss:local}

In this section we show that when calculating Gromov--Witten invariants on a 
manifold with $S^1$-action  only the $S^1$-invariant stable maps 
contribute.  Here is a formal statement of our results.

Let $(P,\om)$ be a closed symplectic manifold,  and, given classes
$a_1,\dots,a_k\in H_*(P)$, let $\al: Z\to P^k$ be a 
(possibly weighted)
pseudocycle that 
represents their exterior product  $a_1\times\dots\times a_k\in 
H_*(P^k)$.  Define  
$$
\oMm_{0,k}(P,J,A;Z): = \ev^{-1}(\ov{\al(Z)}),
$$
where $\ev:\oMm(P,J,A) \to P^k$ is 
the evaluation map.   First, we show that 
the calculation of the corresponding
Gromov--Witten invariant can be localized in $P$ in the following sense.

\begin{lemma}\label{le:cutdown}  The  invariant
$GW_P(a_1,\dots,a_k; A)$ is a sum of contributions, one from each    
connected component of 
the cutdown moduli space $\oMm_{0,k}(P,J,A;Z)$.
\end{lemma}

Now consider the situation when $(P, \om)$ carries an
$S^1$-action.  We assume that $J$ and the 
cycle $\al$ are 
$S^1$-invariant (what this means for cycles is explained in 
Definition~\ref{def:psc}) so that
the cutdown space $\oMm_{0,k}(P,J,A; Z)$ 
 also has an $S^1$-action.

\begin{prop}\label{prop:invarv}  Let  $(P, \om)$ be a
closed symplectic manifold with
 an $S^1$-action $\{\phi_t\}_{t\in \R/\Z}$,
 and  suppose that $J$ and $\al:Z\to P^k$
 are $S^1$-invariant,
where $\al$ represents $a_1\times\dots\times a_k$ as above.
Then a connected component $\Cc$ 
of  $\oMm_{0,k}(P,J,A; Z)$ makes no contribution to
$GW_P(a_1,\dots,a_k; A)$ unless it contains an $S^1$-invariant element.
\end{prop}

The next argument shows that this proposition is precisely what we need.
\MS

\NI
{\bf Proof of Proposition~\ref{prop:invar}.}
 Since the cycles $Z,Z'$ are $S^1$-invariant, the 
torus $T^2$  acts on the cutdown moduli space.  Choose $N$ greater 
than the order of any of the isotropy subgroups of the 
$S^1$-action 
on $M$.  Then the only sections of the 
bundle $P\to S^2$ that are invariant under the action of the subgroup 
$\{(Nt,t): t\in S^1\}$ of $T^2$ are the constant sections $\si_x$
at the fixed points $x\in M^{S^1}$.
Therefore, it
 follows from Lemma~\ref{le:invar} that the fixed points of this 
 circle subgroup
are the same as those for the action of the full torus.
Hence by Lemma~\ref{le:cutdown}
and Proposition~\ref{prop:invarv} the only components of the cutdown 
moduli space ${\oMm}\,\!^{cut}: = \oMm_{0,k}(P,J,A;Z)$ that contribute to the GW 
invariant are those containing $T^2$-invariant elements.

The second statement in Proposition~\ref{prop:invar} goes one step 
further, and claims that the GW invariant is a sum of contributions 
one from each component of the space of invariant elements 
$({\oMm}\,\!^{cut})^{T^2}$ in the cutdown moduli space.  This is 
proved by applying the proof of Lemma~\ref{le:cutdown} to the 
components of $({\oMm}\,\!^{cut})^{T^2}$.  The details are 
straightforward, and are left to the reader.
\QED\MS

We now explain the idea of the proof of Proposition~\ref{prop:invar}
assuming for 
simplicity that $\oMm_{0,k}(P,J,A; Z)=:\Cc$ is connected.  
If $\Cc$ contains no 
$S^1$-invariant elements,  $S^1$ acts with finite 
stabilizers on $\Cc$ 
and hence also on some neighborhood  $\Nn(\Cc)$ of $\Cc$ in $\oMm: = 
\oMm_{0,k}(P,J,A)$.  
We will see that we may 
give the quotient $\Nn(\Cc)/S^1$ an orbifold structure and hence
 construct the regularized moduli cycle
$\ev^{\nu}: \oMm\,\!^\nu\to P^k$ so that its subset
$$
\Cc^\nu: = (\ev^{\nu})^{-1}\bigl(\ov{\al(Z)}\bigr)
$$
has a neighborhood  $\Nn(\Cc^{\nu})$ that supports a {\it free} $S^1$-action.
Moreover, $\ev^{\nu}:\Nn(\Cc^{\nu})\to P^k$ is $S^1$-equivariant.
We will explain below the precise nature of the regularization 
$\oMm\,\!^\nu$, but suppose for now that it is a closed manifold.
It then suffices to apply  the following fact.
Suppose that a closed oriented manifold $X$ supports a free 
$S^1$-action and that $f:X\to P$ is equivariant.  Then $f:X\to P$
may be perturbed to an equivariant map whose image is disjoint from 
the closure of the image of any 
invariant pseudocycle $\al:Z\to P$ of complementary dimension.
This holds because locally $X$ is the product of a transverse slice $Y$ 
with $S^1$, and it suffices to perturb the restriction $f|_Y$ so that 
it is disjoint from $\ov{\al(Z)}$ and then extend by equivariance. 
It is essential here that the action on $X$ is free; otherwise one 
could not extend an arbitrary perturbation of $f|_Y$ to $X$.

 Similar arguments have 
been used by many authors, for example in connection with the 
calculation of the Floer homology of a time independent small function: 
cf. Fukaya--Ono~\cite{FO} and 
Liu--Tian~\cite{LIUT}.  The only difference is we are here dealing 
with an external $S^1$-action (i.e. one on the range of the stable maps) 
rather than a reparametrization  action which lives on the domain.  

To carry out the details of the proof  we will first describe how to 
construct the regularized (or virtual) moduli cycle $\ev^{\nu}: 
\oMm\,\!^\nu\to P^k$.  We shall then prove Lemma~\ref{le:cutdown}, and 
finally the proposition. 
As in McDuff~\cite{Mcv,Mcq}, 
we will use the regularization
process of Liu--Tian~\cite{LIUT}; readers can substitute their
preferred constructions.

\subsubsection{Branched pseudocycles}

To start, we describe what kind of object 
the virtual moduli cycle is.  For more details see~\cite{LIUT,Mcv}.

First, it is a $d$-dimensional 
partially smooth space $\io_X:X_{sm} \to X$.  Here $X$ 
is a compact Hausdorff space (this is called the first topology), 
$\io$ is a bijective continuous map,
and $X_{sm}$ is a union of a finite number of disjoint  smooth 
manifolds $X^i$ of dimensions $i\le d$.  The connected components of 
$X_{sm}$ are called strata.  Maps from one partially smooth space to 
another are given by commutative diagrams
$$
\begin{array}{ccc} X_{sm}&\stackrel{\io} \to& X\\
    \downarrow&&\downarrow\\
Y_{sm}&\stackrel{\io} \to& Y,
\end{array}
$$
but for short they are often written $f:X\to Y$.
Also any compact smooth manifold is a partially smooth space in which $P$ is 
given the usual topology and $P_{sm}$ has one stratum.
Hence a partially smooth map $f:X\to P$  is continuous when thought 
of as a map from the Hausdorff space $X$ to the metric space $P$,
and smooth when restricted to each stratum of
$X_{sm}$.
The virtual moduli cycle is a {\bf compact branched partially smooth 
labelled pseudocycle},
or  {\bf (compact) branched pseudocycle} for short.
This means it is a $d$-dimensional 
partially smooth space such that each $d$-dimensional stratum $X_j$ is oriented,
 has a rational label,
and fits together with $(d-1)$-dimensional strata to form a 
branched manifold.  More precisely, the closure $\ov{X_j}$ 
in $\io(X^d\cup X^{d-1})$ of each component $X_j$ of 
$X^d$ can be given the structure of an oriented 
manifold with boundary;  moreover, when one divides
the top dimensional faces that meet each $(d-1)$-dimensional component 
 into two sets according to their 
orientations, the sum of the labels in each of these sets must 
be equal.  If $X$ has such structure then any map $f:X\to P$ represents 
a unique rational $d$-dimensional homology class. Just as with 
pseudocycles, there is an obvious notion of bordism.  

A  compact branched pseudocycle  $X$ is said to have a  free $S^1$-action 
with local slices if each point  $x\in X$ has a 
neighborhood  $U_{sm}\to U$ that is isomorphic to
the product $Y_{sm}\times S^1\to  Y\times S^1$ with action 
$t\cdot (y,s)= (y,s+t).$

The relevance of these definitions to the current problem is clear from 
the following lemma. 

\begin{lemma} \label{le:free} Suppose that the smooth manifold $P$ supports an 
$S^1$-action, that $X$ is a compact branched labelled pseudocycle with free 
$S^1$-action and that $f:X\to P$ is equivariant.  Then $f$ can be 
perturbed to an equivariant map that is disjoint from the closure of 
the image of any $S^1$-invariant pseudocycle $\al:Z\to P$ of 
complementary dimension.  Hence $f\cdot\al = 0$.
\end{lemma}
\begin{proof}\, 
    As before, local equivariant perturbations of $f$  may be 
    constructed by perturbing $f$ on the local slices $Y$.  The 
    perturbation is constructed by induction over the strata $S$, starting 
    with those of lowest dimension.  The perturbations have the form 
    $f|_{Y\cap S}\mapsto \phi\circ f|_{Y\cap S}$ where $\phi$ is a 
    suitable small diffeomorphism of $P$, and hence, even though we 
    have little control over the way the strata in $Y$ fit together, always 
    extend from $Y\cap S$ to $Y$.
    Further details are left to the 
    reader.
\end{proof}

The above lemma  does not extend to pseudocycles with free $S^1$-action,
since the definition  of pseudocycle does not give us enough control
of the boundary.   As an example, consider the standard $S^1$-action 
on $S^2$ and take $f:\C \oursetminus\{0\}\to S^2 = 
\C\cup\{\infty\}$ and $g: \{pt\}\to\{0\}$.   Then $f$ is a 
pseudocycle representing the fundamental class, and it has a free
$S^1$-action in the sense that $f$ is equivariant with respect to
a free $S^1$-action on its domain.
Nevertheless $f\cdot g \ne 0$.  The reason is the 
following.  By definition,
 $f\cdot g$ is calculated by first perturbing $f$ so that its boundary
 is {\it disjoint} from that of $g$ and then counting 
 transverse intersection points: see~\cite[Chapter~6]{MS2}. In this example, the
boundary of $f$ meets $\im\, g$ in an essential way and we cannot 
prevent this by a hypothesis concerning only the $S^1$-action on
the open set
$Z$; we must work with a closed domain.

\subsubsection{Construction of the regularization}

 The basic idea  in the construction of $\oMmnk$ is
 to perturb  the compactification $\oMm: = \oMm_{0,k}(P,J,A)$ of 
 $\Mm_{0,k}(P,J,A)$ to a cycle 
 of the correct dimension.   The analytic input to the
 construction explained 
 below is  the standard gluing result, see~\cite{LIUT,FO} 
 or~\cite[Chapter~10]{MS2} for example; 
 the rest of the construction is purely topological. 
The most important step in the proof of 
Proposition~\ref{prop:invar} is to construct  
the local uniformizers of  Step 1 below
so that they support a free  $S^1$-action.

The regularization process has four steps.\MS

\NI
{\bf Step 1:}\,\,
Denote by $\Tilde\Bb$  the space of 
$k$-pointed stable maps $\ttau=(\Si(\bu),\bu, \bz)$ where $\ttau$, though not 
necessarily $J$-holomorphic, has the property that the group of 
self-maps $\Ga_{\tau} = \{\ga: \bu\circ\ga = \bu\}$ is finite.  
Let 
$\Bb$  be the space of equivalence classes of such $\ttau$.  (More 
details are given in \S\ref{sss:orb} below.)   The elements of
$\Bb$ are organized into strata, depending on the topological types of
their domains, and one can show that $\Bb$  has 
the structure of an  orbifold in the 
partially smooth category. (Objects in this category are spaces 
$B_{sm}\to B$ with two topologies, where the first is Hausdorff and the second
is a finite union of disjoint Banach manifolds.)
Thus each point $\tau\in \Bb$ has a neighborhood $U$ with a uniformizer 
$(\TU,\pi,\Ga)$ where $\pi:\TU\to \TU/\Ga=U$ identifies $U$ with the 
quotient of $\TU$ by the action of the finite group $\Ga:=\Ga_\tau$. 
Since the elements in $\TU$ are stable maps, constructing $\TU$ amounts to 
choosing a consistent set of parametrizations for the elements $\tau\in U$.
More details are given below.
Because $\oMm\subset \Bb$ is compact, it is contained 
in the union $\Ww$ of a finite number $U_1,\dots,U_N$
of such locally uniformized sets $U$, each of which is a neighborhood 
of some point $\tau\in \oMm$.   Throughout the construction 
one decreases the size of each $U_i$ (and hence increases their number) 
as appropriate.  Our notational  
convention is that objects living on the uniformizers $\TU$ are 
designated with tildes.
\MS

\NI
{\bf Step 2:}\,\, We interpret the operator $\pJ$ as a section of
an orbibundle $\Ll\to \Ww$.
For each $U$ there is a locally trivial bundle $\TLl_{\TU}\to \TU$ on 
which the local isotropy group $\Ga$ acts.  
The fiber of $\TLl_{\TU}$ at $\ttau=(\Si(\bu), \bu, \bz)\in \TU$ is the
space
$$
\TLl_{\ttau}: = L^p\bigl(\Si(\bu), \La_J^{0,1}\otimes \bu^*(TP)\bigr)
$$
of  anti-$J$-holomorphic $1$-forms on $\Si(\bu)$  of class $L^p$
with values in $\bu^*(TM)$.  The perturbations used to define 
$\oMmn$ are built from  sections of the local bundles $\TLl_{\TU}\to \TU$.
In order to extend these local sections, we construct 
another object  $\TLl\to \TWw$ from $\Ll\to \Ww$ that is called a {\bf
multibundle.}
Here  $\TLl\to \TWw$ is a 
collection of compatible maps $\TLl_I\to \TV_I$, where $I$ is a 
subset of the indexing set $\{1,\dots,N\}$ for the $U_i$, $V_I$ is a suitable 
subset of  $\cap_{i\in I} U_i$ and $\TV_I$ (resp. $\TLl_I$)
is the fiber product of 
the $\TU_i$ (resp.$ \TLl_i$) over $V_I$.   The details of this 
construction are not important for what follows.  All we need to know 
is that each section $\ts(\nu)$ of $\TLl\to \TWw$ (called a multisection)
consists  of a compatible collection 
$\{\ts(\nu)_I\}$ of multivalued sections of $\TLl_I\to \TV_I$.  It 
turns out that each $\ts(\nu)_I$ is single valued over the top strata,
but may well be multivalued over lower dimensional strata.\MS

\NI
{\bf Step 3:}\,\, We construct a finite dimensional vector
space $R$ and a map $\nu\mapsto \ts(\nu)$ of $R$ into the space of
multisections  of  $\TLl\to \TWw$ with the property that for
generic small $\nu \in R$ the section $\pJ + \ts(\nu)$ is transverse
to the zero section.
 The vector  space $R$ is a sum $\oplus_{i\in I}R_i$, 
 where $\{U_i\}$ is an open covering of
$\Ww$ and for each $i$ $R_i$ is a suitable finite dimensional space of
sections of $\TLl_{\TU_i}\to \TU_i$.  This space $R_i$ is formed from
the local obstruction bundle.  The
essential requirement is that for each
stable map $\ttau= (\Si(\bu), \bu, \bz)\in \TU_i$ the subspace of
$\TLl_{\ttau}$ formed by the values
$\{\nu(\ttau): \nu\in R_i\}$ projects onto the
cokernel of the linearization $D_\bu$ of $\pJ$ at $\bu$.  The fact 
that suitable finite dimensional spaces $R_i$ exist is a consequence 
of the gluing construction and the compactness of $\oMm$.  To see 
this, choose for each 
$\tau\in U\subset\oMm$:\MS

(a)  a lift $\ttau=(\Si(\bu), \bu, \bz) \in \TU$ of $\tau$; and 
\smallskip

(b)  a subspace  $R_{\ttau}\subset \TLl_{\ttau}$ that covers the 
cokernel of   $D_{\bu}$.\MS

Then  extend the elements $\nu\in R_{\ttau}$ 
by parallel translation along small paths in $M$ to sections 
$\ttau'\mapsto \nu(\ttau')$ of 
$\TLl_{\TU}$ defined over some small neighborhood $\Nn(\ttau)$ of 
$\ttau$ in $\TU$.  By the gluing construction,
 the subspace 
 $$
 R_{\ttau}(\ttau') = \{\nu(\ttau'): \nu\in R_{\ttau}\}\subset \TLl_{\ttau'} 
 $$
projects onto $\coker D_{\bu'}$ when $\ttau'$ is sufficiently 
close to $\ttau$.  Moreover, we can choose this subset $\TU_{\ttau}$ of 
$\TU$ to be 
invariant under the stabilizer group $\Ga_{\tau}$ so that 
it has the form $\pi^{-1}(U_{\tau})$ for some neighborhood $U_{\tau}$ 
of $\tau$ in $\oMm$.
Therefore,  by compactness of $\oMm$, there is a finite set 
$\tau_i$  such that the corresponding pairs $(U_i, R_i): = (U_{\tau_i}, 
R_{\ttau_i})$ have the required properties.

This defines the finite set of local pairs $(U_i,R_i), 1\le i \le N$.  
 One shows that each $\nu\in R_i$, 
when multiplied by a suitable cutoff function, 
gives rise to a multisection  $\ts(\nu)$ of $\TLl\to \TWw$. 
The most important point here is that the construction is local in 
$\oMm$, i.e. for each $\nu\in R_i$ the section
$\ts(\nu)_I = 0$ whenever the closure $\ov{U}_i$ is disjoint
from all the sets $\ov{U}_j, j\in I$.  Now set
$R: = \oplus_i R_i$.
It follows from the construction that the local 
multisections  $\pbar_J + \ts(\nu)_I$ are transverse 
to the zero section for generic small
$\nu\in R$.  Hence the local
zero sets $\TZ^{\nu}_I\subset \TV_I$ are submanifolds of the correct 
dimension $d$. 
\MS

\NI
{\bf Step 4:}\,\, We construct from the local zero sets $\TZ^{\nu}_I$
of $\pJ + \ts(\nu)$ a compact branched $d$ dimensional
pseudomanifold $\oMmnk$.  Its bordism class is independent of choices.
There is a natural projection map
$$
{\rm proj\,}: \oMmnk(P, J,A) \to \Ww
$$
such that each element in the image lies in the zero set of 
the multivalued section $\pJ +
\nu$, and the evaluation map factors through this projection.
Moreover,  the  strata in $\oMmnk$ of dimensions $d,d-1$ project to
the top stratum of $\Ww$, i.e. into stable maps whose domain has a
single component. Therefore, when one evaluates the intersection
number of $\ev:  \oMmnk(P, J,A)\to P^k$ with a cycle in $P^k$ one will
be counting rationally weighted curves $u: S^2\to P$ that satisfy a 
perturbed Cauchy--Riemann equation $\pbar_J u + \nu(u)= 0$. \MS

\begin{definition} The Gromov--Witten invariant
$GW_P(a_1,\dots,a_k; A)$ is the
intersection number of the evaluation map $\ev: \oMmnk(P,J,A)\to P^k$ 
with a generic representing pseudocycle $\al: = Z\to P^k$ 
for the class
$a_1\times\cdots\times a_k$:
$$
GW_P(a_1,\dots,a_k; A): = \ev\cdot \al.
$$
It is zero by definition if the dimensional condition
$\dim\,P + 2c_1(A)  + 2k-6 + \sum_i\dim a_i = k\dim\,P$ is not satisfied.
\end{definition}

Lemma~\ref{le:cutdown} claims that
one would get the same answer by  first cutting down the moduli space
to $\oMm(P,J,A;Z)$
and then regularizing each of its components separately.
\MS

\NI
{\bf Proof of Lemma~\ref{le:cutdown}.}  
 Let $\Cc_j, 1\le j\le \ell,$ be the connected components 
of $\oMm_{0,k}(P,J,A; Z)$. By compactness there is $\eps> 0$ so that
the set 
$$
\Nn^{2\eps}: = \{\tau\in \oMm: d(\ev(\tau),{\ov{\al(Z)}})\le 2\eps\},
$$
where $d$ is the metric in $P^k$, has $\ell$ connected components
$\Nn_j^{2\eps}\supset \Cc_j$.   Now choose pairs $(U_i, R_i)$ as in 
Step 3, where the open subsets $U_i\subset \Bb$ 
separate out the components $\Cc_j$ in 
the following sense: if ${\ov U}_i\cap(\cup_j \Nn_j^{\eps}\Cc_j)\ne\emptyset$ 
and ${\ov U}_k\cap(\oMm \oursetminus (\cup_j\Nn_j^{2\eps})\ne\emptyset$ then $U_i$ and
$U_k$ are disjoint.  Then construct a regularization $\oMmn$ as 
described above.  By the definition of $\Nn_j^{\eps}$, 
the image under $\ev$ of
the set $\oMmn \oursetminus 
{\rm proj\,}^{-1}(\cup_j\Nn_j^{\eps})$ has distance at least $\eps$ from 
${\ov {\al(Z)}}$ and so does not contribute to the intersection
$\ev\cdot \al$.  On the other hand because the construction in Step 
3 is local, the structure of $\Nn(\Cc_j^{\nu}): = 
\oMmn\cap{\rm proj\,}^{-1}(\Nn_j^{\eps})$ 
depends only on the choices made for the open sets covering 
$\Nn_j^{2\eps}$.    Hence if we define the intersection number of
$
\ev:\Nn(\Cc_j^\nu)\to P^k$ with $\al$ as the local contribution of $\Cc_j$ 
to the Gromov--Witten invariant $\ev\cdot \al$, 
this invariant is the sum of local and independent
contributions as claimed.
\QED\MS


\begin{cor}  Let $\al:Z\to P^k$  represent the class 
$a_1\times\dots\times a_k$ in $P^k$.
If $\oMm_{0,k}(P,J,A;Z) = \emptyset$ then 
 $GW_P(a_1,\dots,a_k; A) = 0.$
 \end{cor}

\subsubsection{The moduli space of stable maps as an 
orbifold}\label{sss:orb}

As preparation for the proof of Proposition~\ref{prop:invarv},
we describe the orbifold structure on the 
space of stable maps.  

Let $(T,E)$ be a finite tree where $T$ denotes the set of 
vertices and the relation $E\subset T\times T$ 
describes the set of  oriented edges.
A genus zero stable map with $k$ marked points modelled on $T$ is a tuple 
$$
\Bigl(\{u_\al\}_{\al\in T},\,\{z_{\al\be}\}_{\al E\be},\, 
\{z_i,\al_i\}_{1\le i\le k}\Bigr)
$$
where $u_\al:S^2\rightarrow P$ is a map, $z_{\al\be}\in S^2$ 
denotes the point on  the $\al$-th sphere that attaches 
to the $\be$-th 
sphere, and $z_i\in S^2$ is the $i$th marked point lying on the $\al_i$-th 
sphere. 
Thus its domain $\Si(\bu)$ is the 
quotient of $S^2\times T$ in which 
$(z_{\al\be},\al)\sim (z_{\be\al},\be)$ whenever $\al E\be$.
We require that $u_{\al}(z_{\al\be}) = u_{\be}(z_{\be\al})$ whenever 
$\al E\be$ so that the $u_{\al}$ induce a map
$\bu:\Si(\bu)\rightarrow P$.
The special points $\{z_{\al\be}\,|\, \be\in T\}\cup \{z_i\,|\,\al_i= 
\al\}$ on the $\al$-th sphere are assumed distinct.   The stability 
condition states  that every ghost component (i.e. component of
$\Si(\bu) $ on which $\bu$ is constant) has at least $3$ special 
points.

Two such tuples $\bigl(u_\al,z_{\al\be}, (z_i,\al_i)\bigr)$ and
$\bigl(u_{\al'},z_{\al'\be'}, (z_i',\al'_i)\bigr)$ modelled on $T,T'$ are 
equivalent if there is a tree isomorphism $f:T\rightarrow T'$ and a 
collection $\phi_{\al}\in \PSL(2,\C), \al\in T,$ such that 
$$
   u_{f(\al)}\circ\phi_\al = u_\al, \quad \phi_\al(z_{\al\be}) = 
   z_{f(\al)f(\be)},\quad  
   (\phi_{\al_i}(z_i),f(\al_i))= (z_i',\al_i'). 
$$
We shall call such tuples $(\bu, \bz)$ for short.
The elements $\tau= [\bu,\bz]$ of the moduli space of 
stable maps $\Bb$ are equivalence  classes  of such tuples.
Each stratum consists of equivalence  classes of stable maps modelled 
on a fixed tree $T$ and has an obvious smooth topology.  The 
Hausdorff topology on the whole space is discussed below.

We now describe Liu--Tian's construction of the uniformizers
$(\TU_{\tau},\pi,\Ga_{\tau})$, where
 $\tau= [\bu,\bz]$ is modelled on $T$. Note that each uniformizer 
 $\TU_{\tau}$ is a subset of the ambient space $\Tilde\Bb$.
 Let us first suppose that 
 $\Ga_{\tau} = \{1\}$.  Then the problem is to find a consistent way 
 of parametrizing all the stable maps near $\tau$. 
 Choose a parametrization
 $$
\ttau: =  (\bu,\bz): = \bigl(u_\al,\, z_{\al\be}, \,(z_i,\al_i)\bigr)
$$
of $\tau$. Add 
 the minimum number of points $\bw: = (w_1,\al_{k+1}),\dots, 
 (w_{\ell},\al_{k+\ell})$ to the set of labelled points in $\ttau$ to make its 
 domain stable, i.e. so that
 each component has at least $3$ special points.
 Pick out three of them for each component $\al$ and
 denote by $Y$  the resulting subset of the $z_{\al\be}, z_i$.
The set $\bw$  is chosen to be invariant under the action 
of any element in $\Ga_{\tau}$ that permutes the components of 
$\Si(\bu)$, but so that in each component
no two are on the same orbit of the stabilizer of
this component in $\Ga_\tau$. Next choose for each $j=1,\dots,\ell,$
a small open codimension-$2$ disc ${\bf 
H}_j$ in $P$ that is transverse to the image of $\bu$ at the points 
$\bu(w_{j})$. (This is possible because the ghost components are 
already stable and hence never contain any of the added points $w_i$.)
If $\Dd_T$ denotes the stratum in $\Bb$ containing $\tau$ 
we define 
$
\TU_\tau\cap \Dd_T$ to be a neighborhood of $\ttau$ in the slice
$$
\Ss_T: =  \Bigl\{\bigl(u_\al',z_{\al\be}', 
(z_i',\al_i)\bigr)\,|\, u_{\al_{k+j}}'(w_{j}) \in {\bf H}_j, 
z_{\al\be}' = z_{\al\be} \mbox{ if } z_{\al\be}\in Y,
z_i'=z_i \mbox{ if } z_i\in Y\Bigr\}.
$$

The domains of the
stable maps $(\bu',\bz')$ near $\ttau$ are formed from the domain of $\ttau$ by 
gluing its components via the gluing parameters $a_{\al\be}\in  
T_{z_{\al\be}}(S^2_\al)\otimes T_{z_{\be\al}}(S^2_{\be})$. (Here for 
convenience we denote the $\al$th component by $S^2_{\al}$.)  Assuming 
$r=|a_{\al\be}|$ is sufficiently small we glue 
$S^2_{\al}\oursetminus B_r(z_{\al\be})$ to 
$S^2_{\be}\oursetminus B_r(z_{\be\al})$ along their boundaries by a 
rotation determined by $\arg(a_{\al\be})$.
To describe this more precisely,  let us 
suppose for simplicity that the tree $T$ has two vertices $\{0,\infty\}$
and one edge, so that $\Dd_T$ has codimension $2$. We may suppose that
$$
z_{0\infty} = \{0\}\;\in\; \C\cup\{\infty\}\;=\; S^2_0,\quad
z_{\infty 0} = \{\infty\}\;\in \;\C\cup\{\infty\}\;=\; S^2_{\infty}.
$$
There are two special points $y_{\al m}\in Y$ on
each component that, together with  $z_{0\infty}, z_{\infty 0}$,  are fixed 
on the slice $\Ss_T$.  By minimality  the added points $w_{j}$  (if 
there are any)
form a subset of 
the four  points  $y_{01},y_{02}, y_{\infty 1}, y_{\infty 2}$.
Again, for the sake of clarity, let us suppose that there is one added 
point $w_1: = y_{01}$, 
and that $z_1=y_{02}, z_2: =y_{\infty 1}, z_3: = y_{\infty 2}$. 
Consider the tuple $\bigl(a; u', (z_i',\al_i)_{1\le i\le k}\bigr)$ 
where the domain $S_a$ is the sphere
$$
S_a: = \bigl(S^2_0\oursetminus B_r(0)\bigr)\cup
\bigl(S^2_\infty\oursetminus B_r(\infty)\bigr),
$$
$u':S_a\rightarrow P$ is close to $\bu$ in 
the obvious $C^0$-sense and 
$z_i'\in S^2_{\al_i}\oursetminus B_r\subset S_a$ is close to the image 
of $z_i$.
Each such sphere $S_a$ has a unique identification $\psi_a: S_a\rightarrow 
S^2$ with
$S^2$ under which the four marked points  $w_1,z_1',z_2',z_3'$
 are taken to 
$0,1,\infty,c(a)$, where $c(a)$ is their cross ratio. (Here we 
identify $w_1\in S^2_0$ with its image in $S^2_a$ in the obvious way.
Note also that one can define $a$ so that $c(a) =a$.)
Hence the tuple  $(a; u', (z_i',\al_i))$ can  be written uniquely as
a stable map $(u'',\bz'')\in \Mm_{0,k}(S^2,A,J)$ 
where $u'':S^2\rightarrow P$ is the composite 
$u'\circ(\psi_a)^{-1}$, and  $z_i'' = \psi_a(z_i'), i=1,\dots,k.$ 
Conversely, each stable map that is sufficiently close to $\ttau$ 
does correspond to a unique tuple $(a; u', (z_i',\al_i))$ since the gluing 
parameter $a$ is determined by the cross ratio of the four marked 
points $0,z_1'',z_2'',z_3''$.  Therefore
we may extend the slice $\Ss_T$ by setting
$$
\begin{array}{l}
\Ss: = \Bigl\{\bigl(a; u', (z_i',\al_i)\bigr)\,|\, u':S_a\rightarrow P,
u'(w_{k+j})\in {\bf H}_j,\\\qquad\qquad\qquad\qquad\qquad
z_i'=z_i \mbox{ if } z_i= y_{\al m}\mbox{ for some }\al, m\Bigr\}.
\end{array}
$$
Finally we define $\TU_\tau$ to be
a neighborhood of $\ttau$ in $\Ss\cup\Ss_T$.  The projection to $\Bb$ 
is given by dividing by the reparametrization group, i.e. by taking a 
stable map to its equivalence class.
(We have not given a 
satisfactory description of the topology on $\Bb$: for this 
see~\cite{LIUT,Mcv}.)

Now suppose that $\Ga_{\tau}\ne \{\1\}$.
We must extend the action of $\Ga_{\tau}$ to
$\TU_\tau$.  Suppose first that
$\Ga_{\tau}$ is a rotation group of order $n>1$ with generator $\ga\in 
\PSL(2,\C)$ that acts on a single component $\al_0$ of $\Si(\bu)$.
This component can have at most two special points $y_i$. 
  Let us suppose that it has precisely two, say $y_1,y_2$,  and 
therefore one added point
that we will call $w_1$.  We may suppose $w_1$ chosen so that
 the set $u_{\al_0}^{-1}(u_{\al_0}(w_1))$ contains
$n$ distinct points at which $du_{\al_0}\ne 0$.
Choose  disjoint little discs in $S^2_{\al_0}$ about these points
that are permuted by $\ga$. For any element
$(\bu', \bz')$ that is close to $\ttau$ and
in the same stratum,  $(u_{\al_0}')^{-1}({\bf H}_1)$
is a collection of $n$ points, one in each of 
the little discs.  Therefore there is
 unique point  $w'$  in the little
disc containing $\ga(w_1)$ such that
$
\bu'(w')\in {\bf H}_1,
$
 and we define $\psi_{\bu'}^\ga\in
\PSL(2,\C)^{|T|}$ to be the unique element that 
acts as the identity in all components except for the $\al_0$-th and there
fixes $y_1, y_2$ and takes $w_1$ to $w'$.  Then set 
$$ 
\ga\cdot (\bu', \bz') = 
(\bu'\circ \psi_{\bu'}^\ga, \bz')\;\in\; \TU_\tau.
$$
It is not hard to check that  
this does define an action of $\Ga_\tau$  on a neighborhood of
$\ttau$ in $\TU_\tau\cap \Dd_T$. 

It extends over the full neighborhood $\TU_\tau$ by acting  
 on the gluing parameters $a$.  We give a precise 
description in the case with $|T|=2$ considered above.
If $\bigl(a; u', w_1, (z_i',\al_i)\bigr)\in \Ss$ is sufficiently close to 
$\ttau$, then $(u')^{-1}{\bf H}_1\subset S_a$ consists of $n$ points
with precisely one, call it $w'$, in the 
little disc containing $\ga(w_1)$.  Because the map $a\mapsto c(a)$ is a 
diffeomorphism, there is a unique gluing parameter
$\ga(a)$ for which there is a biholomorphic map
$$
\psi_{\ga}: \bigl(S_{\ga(a)};\, y_{01}= w_1, y_{02}, y_{\infty 1}, y_{\infty 
2}\bigr)\to
\bigl(S_a;\, y_{01}= w', y_{02}, y_{\infty 1}, y_{\infty 
2}\bigr).
$$
We define
$$
\ga\cdot \bigl(a;u', (z_i',\al_i)\bigr): =
\Bigl(\ga(a); u'\circ \psi_{\ga}, \bigl(\psi_{\ga}^{-1}(z_i'),\al_i\bigr)\Bigr)\in 
\Ss.
$$
 Alternatively, if we write the elements of $\Ss$ in the 
form  $(u'', z_i'')$ where $u'':S^2\to P$ then 
$$
\ga\cdot (u'', z_i'') = (u''\circ h, h^{-1}(z_i'')),\quad 
h: = \psi_a\circ\psi_\ga\circ(\psi_{\ga(a)})^{-1}.
$$
Again, one can check that this gives a well defined action of 
$\Ga_{\tau}$ on $\TU_{\tau}$.
Its continuity 
(which is somewhat tricky) is proved in~\cite[\S4.2]{Mcq}; see 
also~\cite{Mcv}.   

The construction for other groups $\Ga_{\tau}$ is 
similar. 
We need to consider the case when $\Ga_{\tau}$  acts in a single component 
with one  or no special points; then consider products of such actions; 
and finally consider an action that also permutes the components. 
These extensions are described in~\cite{LIUT}. 

\MS

\NI{\bf Proof of Proposition~\ref{prop:invarv}}.
 Let $\Cc$ be a component of $\oMm_{0,k}(P,J,A;Z)$ on which the 
 induced action of $S^1$ is locally free. 
 We will show that its regularization $\Nn(\Cc^{\nu})$ can be 
 constructed so as to support a free $S^1$-action.  The result 
 then follows from Lemma~\ref{le:free}.
 
 We show below that $\Cc$ can be covered by $S^1$-invariant
 sets $W_j$ such that their 
 uniformizers $\TW_j$, as well as the uniformizers of 
 all sets they meet, support a free $S^1$-action with local slices 
 $\TY_j$.
 Granted this, the construction of the $(U_i,R_i)$ in Step 3
can be made so that the sections in $R_i$ are 
$S^1$-invariant. To see this,  choose for each $\tau\in W_j$ with lift $\ttau$
 a suitable finite 
dimensional space $R_{\ttau}$ of the fiber $\TLl_{\ttau}$, 
extend its elements to a neighborhood ${\TY}_{\ttau}$ of $\ttau$ in the 
slice ${\TY}_{j}$ by parallel 
translation in $P$ and then extend over the product
$\TU_{\ttau}: = {\TY}_{\ttau}\times S^1$ using the 
$S^1$-action on  $P$.    Since this action 
preserves $J$  the transversality conditions continue to 
hold over the $S^1$-orbit.  

Hence the local zero sets $\TZ^\nu_I$ all carry a free $S^1$-action
with  local slices.  The local virtual cycle 
$\Nn(\Cc^{\nu})$  is made from these zero sets using partitions of unity,
and one can check that its construction can carried out
 in a way that respects
the $S^1$-action.  Moreover the induced 
$S^1$-action is free because each point
in $\Nn(\Cc^{\nu})$ projects to one of the sets $\TV_I$ and hence to
 $\TU_i, i\in I,$ where the action is free by construction: for 
 details see Proposition~4.13 in~\cite{Mcv}.

 Hence it remains to construct the $W_j$.
 To do this, we construct a different set of local uniformizers 
 $(\TW_\tau, \pi,\TGa_\tau)$ for 
a neighborhood $\Nn(\Cc)$ of $\Cc$ in $\Bb$ whose  
stabilizer subgroups $\TGa_\tau$ 
incorporate not only the automorphism groups $\Ga_\tau$ of the stable maps
$\tau\in V_i$  but also the (finite) stabilizer subgroups 
$\Stab(\tau)\subset S^1$ of the locally 
free $S^1$-action on $\Nn(\Cc)$.
This amounts to defining an orbifold structure on 
the quotient $\Nn(\Cc)_S: = 
\Nn(\Cc)/S^1$ whose elements are equivalence classes 
$[\Si(\bu),\bu,\bz]_S$, 
where  the equivalence relation $\sim_S$ is generated by the previous 
relation $\sim$ coming from the action of the reparametrization group 
together with the equivalence 
$$
(\Si(\bu),\bu,\bz) \sim_S (\Si(\bu),\phi_t\circ \bu,\bz),\qquad t\in 
S^1,
$$
where $\phi_t:P\to P$ denotes the action of $t\in S^1$.

The first task is to 
 define the local group $\TGa_\tau$ at $\tau\in \Cc$.
Choose a parametrization $\ttau= (\bu,\bz)$.  
Let $\Ga_{\tau}: = \{\ga\in \Aut(\Si(\bu)): 
\bu\circ \ga=\bu\}$ denote its automorphism group,
and denote by $N$ the order of the stabilizer subgroup 
$\Stab(\tau)$  of $\tau$ in $S^1$.  Then define
$$
    \TGa_{\tau}: = \{(\ga,k)\in \Aut(\Si(\bu))\times \Stab(\tau)
    \,|\, \bu\circ\ga = \phi_{k/N}\circ\bu\}.
$$
There are exact sequences
$$
 \Stab(\tau)'\;\hookrightarrow\; \Stab(\tau)\;\twoheadrightarrow\; 
 \Stab(\tau)'',\quad
 \Ga_\tau\times \Stab(\tau)' \;\hookrightarrow \;\TGa_{\tau} 
\;\twoheadrightarrow\; \Stab(\tau)''.
$$
where $\Stab(\tau)'=\{t\in \R/\Z\,|\,\phi_t\circ\bu = \bu\}$
is the stabilizer subgroup of the image of $\bu$ in $P$.

We must show that every $\tau\in \Cc$ has a neighborhood $W_\tau$
with a uniformizer $(\TW_\tau,\pi,\TGa_{\tau})$ such that $
(\TW_\tau,\TGa_\tau)$ is equivariantly isomorphic  to a product
$\TY_{\tau}\times S^1$ with action induced by
$$
(\ga,k)\cdot (\bu',t): =  (\bu'\circ\ga,t-k/N).
$$
Then the projection $\pi$ given by $\pi(\bu',\bz',t): = \phi_t\circ 
(\bu',\bz')$ 
is well defined and $S^1$-equivariant, and does quotient out by the 
action of $\TGa_\tau$.  For these formulas to make sense
$\TY_{\tau}$ must be invariant under the action of 
$\Ga_\tau\times \Stab(\tau)'$.  To find such a slice $\TY_{\tau}$
we will use the fact that, by hypothesis, $\Stab(\tau)$ is finite.

Suppose first that  $\Ga_{\tau} = \{\1\}$.  Choose the added points
$w_j$ to be generic, i.e. so that  $d\bu(w_j)\ne 0$ and
the stabilizers  $\Stab(\bu(w_j))$
of the points $\bu(w_j)\in P$ are 
as small as possible, and then choose the slices ${\bf H}_j\subset P$ to be 
$\Stab(\bu(w_j))$-invariant.  Note that 
$\Stab(\tau)'\subseteq \Stab(\bu(w_j))$ for all $j$.
Suppose in addition that it is possible 
to choose one of the added points, say  $w_2$, so that 
$\Stab(\bu(w_2))$ is finite.  Then there is a  
$\Stab(\bu(w_2))$-invariant codimension $1$ disc $X$ through $\bu(w_2)$ 
that is transverse both to the $S^1$-action and to ${\bf H}_2$, and we set  
${\bf H}'_2: = {\bf H}_2\cap X$.  Then
$$
\TY_{\tau}: = \{\ttau'\in \TU_{\tau} \,|\, \bu'(w_2)\in {\bf H}'_2\}
$$
is a slice for the induced local $S^1$-action on $\TU_{\tau}$; in 
particular it is $\Stab(\tau)'$-invariant.  Hence 
we may take  
$$
\TW_{\tau}: = \TY_{\tau}\times S^1,\qquad W_{\tau}: = \{\phi_t\cdot\tau'\,|\, \tau'\in \pi(\TY_{\tau}),\, t\in 
S^1\}.
$$
The projection $\TW_{\tau}\to W_{\tau}$ is given by
$(\ttau',t)\mapsto \phi_t\cdot\tau'$. Note that the uniformizer $\TW_\tau$
is no longer  a subset of $\Tilde\Bb$, but is defined so that it 
supports a free $S^1$-action.

Suppose now that we cannot choose $w_2$ as above.  (For example, 
there may be no need to add any $w_j$ or the unstable components may 
all map into the fixed set.)  Then, we choose any point
$w_0\in \Si(\bu)$ so that $\Stab(\bu(w_0))$ is finite. (This exists
since $\Stab(\tau)$ is finite.)
We choose the slice $X$ through $\bu(w_0)$ 
as before and define $\TY_{\tau}\subset \TU_{\tau}$ by the 
condition $\bu'(w_0)\in X$.  It is obvious what this means when 
$[\bu',\bz']$ is in the same stratum at $\tau$.
One extends to a neighboring strata as before.
Note that in this case
 $w_0$ lies on a component with at
least $3$ special points.

Finally suppose that $\Ga_{\tau} \ne \{\1\}$.  As before we treat the 
case when $\Ga_\tau$ is cyclic and acting on one component of 
$\Si(\bu)$. 
Because this component contains at most $2$ special points, $w_0$ (if 
it has been defined) always lies on some other component. Thus the only 
case that needs special consideration is when $w_2$ 
lies on the component on which $\Ga_{\tau}$ acts
and so equals the point previously called $w_1$. 
But then we may simply repeat the previous 
construction for the action of $\Ga_{\tau}$, 
replacing ${\bf H}_1$ by ${\bf H}_1': = {\bf H}_1\cap X$.  This defines an
action of $\Ga_{\tau}$ on $\TY_{\tau}$ and hence completes the 
construction.\QED\MS

\section{Applications and Examples}\label{sec:appl}

In the first section,  we describe the small quantum cohomology 
of toric manifolds.  
Next, we work out $\Ss(\La)$ in specific cases to illustrate 
what may happen when the hypotheses of the main theorems do not 
hold.

\subsection{The small quantum homology of smooth toric varieties}
\labell{ss:toric}

This section describes the general form of a set of 
generators and relations for the 
small quantum cohomology ring $\QH^*(M)$ of a toric manifold: see 
Proposition~\ref{prop:qhtoric}.
In the case of a Fano variety the description is 
completely explicit; it is determined by
a simple algorithm from
the moment polytope $\Delta$
and agrees with Batyrev's  presentation~\cite{Ba}.  
In the NEF case we show that
$\QH^*(M)$ is determined by a simple algorithm involving
its moment polytope $\Delta$ together 
with the  Seidel elements of the circle actions corresponding to the 
primitive outward normals  $\eta_1,\dots,\eta_N$ to the facets of $\Delta$.
 (Of course, calculating the Seidel elements is a very nontrivial problem that we do not attempt.)
In the general case, the relations correspond to certain products of the Seidel elements but are not immediately determined by them.
Our result 
elaborates on a very small  part of Givental's work on the 
mirror conjecture: see
Cox--Katz~\cite[Examples~8.1.2.2,~11.2.5.2]{CK}.   
Throughout we work with quantum cohomology with the Novikov ring 
coefficients defined in \S\ref{ss:qh}
though one can extend the result to the full Novikov ring:
see Remark~\ref{rmk:Novik}. 
Batyrev used complex 
coefficients; for a discussion of the relation of these 
coefficient systems see~\cite[8.1.3]{CK}.

Before beginning our computation,
let us review a few facts about quantum cohomology.
First, as in~(\ref{eq:vcheck}),
define a valuation $\check v$ on $\Q[x_1,\ldots,x_N] \otimes \check\Lambda$ 
by 
$$
\check v \Bigl(\sum_{d,\ka} a_{d,\ka} \otimes q^d t^{\ka}\Bigr) = 
\min \{\ka\;|\;\exists \, d :  a_{d,\ka}\ne 0\}.
$$

\begin{lemma}\labell{le:qhtor}
Let $(M,\omega)$ be a symplectic manifold.  Fix $x_1,\ldots,x_N
\in H^*(M)$, and consider the natural homomorphisms of rings
$$
\theta: \Q[x_1,\ldots,x_N] \to H^*(M),  \quad \mbox{and }
$$
$$
\Theta: \Q[x_1,\ldots,x_N] \otimes \check\Lambda \to \QH^*(M).
$$
\smallskip

\NI {\rm  (i)}
If $\theta$ is surjective, then $\Theta$ is also surjective.
Further, given $z \in \QH^*(M)$, there 
is
$\tilde{z} \in \Q[x_1,\ldots,x_N]\otimes~\check\Lambda$ so that
$\Theta(\tilde{z}) = z$ and so that 
$\check v(\tilde{z}) \geq \check v(z)$.
\smallskip

\NI {\rm  (ii)}
Let  $p_1,\ldots,p_m \in \Q[x_1,\ldots,x_N]$  generate the
kernel of $\theta$,  
and suppose $q_1,\dots,q_m $ $\in \Q[x_1,\ldots,x_N]\otimes \check{\La}$ are such
that $\Theta(q_i) = 0$ and $\check{v}(p_i-q_i) > 0$ for all $i$.
Then $q_1,\ldots,q_m$ generate
the kernel of $\Theta$.
\end{lemma}

\begin{proof}\,  
Fix $\hbar> 0$ such that $\hbar$ is less than
the energy $\om(B)$ of every class $B\ne 0$ 
that contributes to the quantum multiplication, i.e. for which 
there is a nonzero three point Gromov--Witten invariant.
Then $\check v(\al*\be - \al\cup \be)\ge \hbar$ for all $\al,\be\in H^*(M)$, and 
hence
\begin{equation}\label{eq:vmult}
    \check v\bigl(\Theta(\tilde z) -\theta(\tilde z)\bigr)\ge \hbar,\quad 
    \forall \tilde z\in \Q[x_1,\ldots,x_N].
\end{equation}
By possibly shrinking $\hbar$, we can also assume that $v^*(p_i - q_i) > \hbar$
for all $i = 1,\ldots,m$.

Fix $z \in \QH^*(M)$. To prove (i) 
it is enough to  find $\tilde{z}  
\in \Q[x_1,\ldots,x_N] \otimes \check\Lambda$ 
so that 
$$
\check v(z - \Theta(\tilde{z})) \geq \check v(z) + \hbar,
\quad\mbox{ and }\;\; 
\check v(\tilde{z}) \geq \check v(z),
$$ 
since then the argument can be completed by induction.
Write
$$
z = \sum_{i=1}^k~z_i~\otimes~q^{d_i}~t^{\kappa_i} + r,
$$
where $\check v(r) \geq \check v(z) + \hbar$, 
$z_i \in H^*(M)$,
$d_i \in \Z$, and $\kappa_i \geq \check v(z)$. 
Since $\theta$ is surjective,  there exists $\tilde{z_i} \in 
\Q[x_1,\ldots,x_N]$
so that $\theta(\tilde{z_i}) = z_i$. 
Then 
$\check v\bigl(z_i - \Theta(\tilde{z_i})\bigr) \geq \hbar$ 
by~(\ref{eq:vmult}); so
let $\tilde{z}  = \sum_{i=1}^k \tilde{z_i} \otimes q^{d_i} t^{\kappa_i}$.

Now fix $\tilde y \in \ker \Theta$.  
To prove (ii),
it is enough to find $\tilde{z} \in \langle q_1,\ldots,q_m \rangle$
so that 
$$
\check v(\tilde{z}-\tilde y) \geq \check v(\tilde y) + \hbar,
\quad \mbox { and }\;\; 
\check v(\tilde{z}) \geq \check v(\tilde y),
$$ 
since then this argument can  also be completed by induction.
Write
$$
\tilde y = \sum_{i=1}^k \tilde y_i \otimes q^{d_i} t^{\kappa_i} + \tilde r,
$$
where $\check v( \tilde r) \geq \check v( \tilde y) + \hbar$, 
$\tilde{y_i} \in \Q[x_1,\ldots,x_N]$, $d_i \in \Z$, and 
$\check v(\tilde y) + \hbar >  \kappa_i \geq \check v(\tilde y)$.
We may also assume that $(d_i,\kappa_i) \neq (d_j,\kappa_j)$ if $i \neq j$.
Note that by~(\ref{eq:vmult})
$$
0 = \Theta(\tilde y) = \sum \Theta(\tilde y_i \otimes q^{d_i} t^{\kappa_i}) + 
\Theta(\tilde r)
= \sum \theta(\tilde y_i) \otimes q^{d_1} t^{\kappa_i}  + \tilde r',
$$
where $\check v(\tilde r') \geq \check v(\tilde y) + \hbar$.
Therefore, for all $i$, $\theta(\tilde y_i) = 0$,
and hence $\tilde y_i$ lies in the ideal generated by $p_1,\ldots,p_m$.
Hence, there exists $\tilde{z_i} \in \langle q_1,\ldots,q_m\rangle$,
so that $\check v(\tilde y_i - \tilde{z}_i) \geq \hbar$.
Let $\tilde{z} = \sum \tilde{z}_i \otimes q^{d_i} t^{\kappa_i}.$
\end{proof}

We will now give a brief review of toric geometry.
Good basic references are Cox--Katz~\cite[Ch~3]{CK} and 
Batyrev~\cite{Ba}.   

Consider a torus $T$ with Lie algebra $\ft$ and lattice $\ell$.
Let $(M,\omega)$ be a smooth toric variety with moment
map $\Phi: M \to \ft^*$, chosen so that each of its components is mean 
normalized.
Let $\Delta \subset \ft^*$ be the image of the moment map.
Let $D_1,\ldots, D_N$ be the facets of $\Delta$ (the
codimension one faces), and
let $\eta_1, \ldots \eta_N \in \ell$ denote the outward
primitive integral normal vectors.\footnote
{
Choosing the $\eta_i\in \ft$ to be the outward rather than the inward normal
is more natural in our context.  For then the corresponding circle 
action has $\Phi\inv(D_i)$ as its maximal fixed point component, and 
it is this, rather than the minimal fixed point component, that 
is seen by the Seidel element: cf. Theorem~\ref{thm:max}.  However, 
the authors of~\cite{Ba,CK} make the other choice, defining the 
polytope $\Delta$  by equations of the form $\{v\in 
\ft^*:\langle \eta_i, v\rangle \ge -a_i\}$.  If we take the inward  
normals then in the definition of $SR_Y$ in~(\ref{eq:SRom}) 
$\be_I$ should be replaced by $-\be_I$.}
Let $\ell^*\subset \ft^*$ denote the lattice dual to $\ell$.
Let $\Sigma$ be the set of subsets
$I = \{i_1,\ldots,i_k\} \subseteq \{1,\ldots,N\}$
so that   $D_{i_1} \cap \cdots \cap D_{i_k} \neq \emptyset$.
Define two ideals in $\Q[x_1,\ldots,x_N]$:
$$
P(\Delta) = \left< \left. \sum (\xi,\eta_i)\, x_i \right| \xi \in \ell^* \right>,
\quad \mbox{and} \quad
SR(\Delta) = \langle x_{i_1} \cdots x_{i_k}  \mid 
\{i_1,\ldots,i_k\} \not\in \Sigma \rangle.
$$
A subset $I \subseteq \{1,\ldots,N\}$ is called {\bf primitive}
if $I$ is not in $\Sigma$ but every proper subset is.
Clearly,
$$
SR(\Delta) = \left< x_{i_1} \cdots x_{i_k} \mid
\{i_1,\ldots,i_k\} \subseteq \{1,\ldots,N\} \mbox{ is primitive} \right>.
$$

The map which sends $x_i$ to the Poincar\'e dual of $\Phi\inv(D_i)$
(which we shall also denote by $x_i\in H^2(M)$)
induces an isomorphism
\begin{eqnarray}\label{eq:H2}
\Q[x_1,\ldots,x_N]/ (P(\Delta) + SR(\Delta)) \cong H^*(M,\Q).
\end{eqnarray}
Moreover, there is a natural isomorphism between
$H_2(M;\Z)$
and the set of tuples $(a_1,\ldots,a_N) \in \Z^N$ 
such that $\sum a_i \eta_i = 0$, under which
the pairing between such an element of $H_2(M,\Z)$ and
$x_i$ is $a_i$.
The linear functional $\eta_i$ is constant on $D_i$;
let $\eta_i(D_i)$ denote its value.
Under the isomorphism of~(\ref{eq:H2}) (extended to real coefficients)
\begin{eqnarray}\label{eq:H3}
    [\om] = \sum_i \eta_i(D_i) x_i,\quad \mbox{and} \quad c_1(M)=\sum_i x_i.
\end{eqnarray}

We are now ready to examine the quantum cohomology of
a toric variety.
The Seidel representation  in cohomology is the homomorphism
$$
\Ss^*: \pi_1(\Ham(M, \om))\to \QH_\ev(M;\La)^\times,\quad \La\mapsto \PD(\Ss(\La)),
$$
where $\QH_\ev(M;\La)^\times$ is the group of even units in $\QH^*(M)$
and $\Ss$ is the representation in homology.
For each $\eta_i$ define  $\Phi^{\eta_i}: M \to \R$
to be the composite of the moment map $\Phi:M\to \ft^*$ with the 
linear functional  $\eta_i \in \ft=\Hom(\ft^*,\R)$.  Thus
$\Phi^{\eta_i}$ is  the  moment map for the circle action $\La_i$ 
with tangent vector $\eta_i\in \ft$ and with $F_{\max} = \Phi\inv(D_i)$.  
Denote:
$$
\Ss^*(\La_i) = y_i\otimes q^{-1} t^{-\eta_i(D_i)}\in \QH_\ev(M;\La)^\times.
$$
By Theorem~\ref{thm:max} and the
formula given for Poincar\'e duality in \S\ref{ss:qh},
$y_i = x_i + $ higher order terms, where the terms are ordered by $\check v$.   

Given any face of $\Delta$, let $D_{j_1}, \ldots, D_{j_\ell}$
be the facets that intersect to form this face.
The {\bf dual cone} is  the set of elements in $\ft$ which can
be written as a positive linear combination of 
$\eta_{j_1},\ldots,\eta_{j_\ell}$.
Every vector in $\ft$ lies in the dual cone of a unique face of $\Delta$.
Therefore, given any subset
$I = \{i_1,\ldots,i_k\} \subseteq \{1,\ldots,N\}$
there is a unique face of $\Delta$ so
that  $\eta_{i_1} + \cdots + \eta_{i_k}$ lies in its dual cone.
Let $D_{j_1},\ldots,D_{j_\ell}$ be the  facets that intersect
to form this unique face.  Then there exist unique positive integers 
$c_1,\ldots,c_\ell$ so that 
$$
\eta_{i_1} + \cdots + \eta_{i_k} - c_1 \eta_{j_1} -
\cdots - c_\ell \eta_{j_\ell} = 0.
$$
Batyrev showed that if $I$ is primitive the sets $I$ and $J = 
\{j_1,\dots,j_{\ell}\}$ are disjoint.
Let $\beta_I \in H_2(M,\Z)$ be the class corresponding to the above relation.
By~(\ref{eq:H3}), we see that
\begin{eqnarray*}
c_1(\beta_I) & = &  k - c_1 - \cdots - c_\ell, \quad \mbox{and}\\ 
 \omega(\beta_I) &= &\eta_{i_1}(D_{i_1}) + \cdots + \eta_{i_k}(D_{i_k})
- c_1 \eta_{j_1}(D_{j_1}) - \cdots - c_\ell \eta_{j_\ell}(D_{j_\ell}).
\end{eqnarray*}
Since
$\eta_{i_1} + \cdots + \eta_{i_k} = c_1 \eta_{j_1} + \cdots + c_\ell 
\eta_{j_\ell}$,
the corresponding circle actions are also equal.
Using the fact that the Seidel representation is in fact a homomorphism,
we have
\begin{eqnarray*}
&& y_{i_1} *\cdots* y_{i_k}\otimes q^{-k} 
t^{- \eta_{i_1}(D_{i_1}) - \cdots -\eta_{i_k}(D_{i_k})}
=\\
 &&\qquad\quad \qquad
{y_{j_1}}^{c_1} * \cdots * {y_{j_\ell}}^{c_\ell}\otimes q^{-c_1 -  \ldots - 
c_\ell}\;
t^{-c_1 \eta_{j_1}(D_{j_1}) - \cdots - c_k \eta_{j_\ell}(D_{j_\ell})}.
\end{eqnarray*}
Therefore
$$
y_{i_1} *\cdots* y_{i_k} - {y_{j_1}}^{c_1} * \cdots * 
{y_{j_\ell}}^{c_\ell}\otimes
q^{c_1(\beta_I)}
t^{\omega(\beta_I)} = 0.
$$

Since $x_1,\ldots,x_N$ generate $H^*(M)$,
by Lemma~\ref{le:qhtor}  the natural  homomorphism
$$
\Theta : \Q[x_1,\ldots,x_N]\otimes \check\La \to \QH^*(M)
$$
which takes $x_i$ to the Poincar\'e dual of $\Phi\inv(D_i)$ is surjective. 
 To compute $\QH^*(M)$, we need
to find the kernel of $\Theta$. 
By Lemma~\ref{le:qhtor}, there exists
\begin{equation}\label{eq:Yi}
Y_i = x_i + \mbox{ higher order terms }
\end{equation}
such that $\Theta(Y_i) = y_i$.
Define an ideal $SR_Y(\Delta) \subset \Q[x_1,\ldots,x_N]\otimes 
\check\La$ by
\begin{eqnarray}\label{eq:SRom}
SR_Y(\Delta) & = & \left< \left. Y_{i_1} \cdots Y_{i_l} - 
{Y_{j_1}}^{c_1}  \cdots {Y_{j_k}}^{c_k} \otimes 
q^{c_1(\beta_I)} t^{\omega(\beta_I)}\,\right|\right. \\\notag
&&\qquad \qquad \qquad \qquad I = \{i_1,\ldots,i_l\} \mbox{ is primitive}\} \left.\right>,
\end{eqnarray}
where the $Y_i$ are as in~(\ref{eq:Yi}).  
Note that $SR_Y(\Delta)$ depends on  the 
$Y_i$.  
Additionally, even if 
the Seidel element 
$y_i$ is known, it is not
in general possible to describe 
its lift
$Y_i$ without prior knowledge of the ring
structure on $\QH^*(M)$.
On the other hand,  $SR_Y$ is clearly contained in the kernel  of  $\Theta$.
Moreover, Batyrev shows that $\omega(\beta_I) > 0$ for all
primitive $I$.
Hence, applying Lemma~\ref{le:qhtor}, we obtain the
following proposition:

\begin{proposition}\label{prop:qhtoric}
Let $\QH^*(M)$ denote the small quantum cohomology of the toric 
manifold $(M, \om)$.
The map which sends $x_i$ to the Poincar\'e dual of $\Phi\inv(D_i)$ induces
an isomorphism
$$
\Q[x_1,\ldots,x_N]\otimes \check\La/ (P(\Delta) + SR_Y(\Delta)) \cong \QH^*(M).
$$
\end{proposition}

This is especially simple in the Fano case.  

\begin{example}[Fano toric varieties]\label{ex:torifan}\rm
Assume that $M$ is Fano, i.e. that $c_1(B) > 0$ for every
class $B\in H_2(M)$ with a  holomorphic representative.
In this case the higher order terms in $\Ss(\La_i)$ vanish 
by part (iii) of Theorem~\ref{thm:max}.  Therefore  $y_i = x_i$ for all $i$, so 
that
we may set $Y_i = x_i$.  Hence
$$
SR_Y(\Delta) = \left< x_{i_1} *\cdots* x_{i_l} -
{x_{j_1}}^{c_1} * \cdots *{x_{j_k}}^{c_k}\otimes e^{\beta_I}  \mid
I = \{i_1,\ldots,i_l\} \mbox{ is primitive}\} \right>.
$$
This gives exactly the formula for the small quantum cohomology
of a Fano toric variety given by Batyrev and proved by Givental.
\end{example}

\begin{example}[NEF toric varieties]\label{ex:torinef}\rm
Now assume that $M$ is NEF, i.e. that $c_1(B) \ge 0$ for every
class $B\in H_2(M)$ with a  holomorphic representative.
Now there may be   
higher order terms in the Seidel 
elements $y_i$.  However, part (ii) of Theorem~\ref{thm:max}
implies that the higher order terms in $\Ss^*(\La_i)$
have the form 
$$
\alpha_B\otimes q^{-1 + c_1(B)} t^{\eta_i(D_i) + \omega(B)}
$$
where $B \in H_2(M)$ satisfies $c_1(B) = 0$ or $1$.
Since $\Ss^*(\Lambda_i)$ is homogeneous of degree $0$, every nonzero
$\alpha_B$ must have degree $0$ or $2$.
Therefore $\al_B$ either lifts to the unit $\1$ in 
$\Q[x_1,\dots,x_N]\otimes \check\La$ or to some linear 
combination of the $x_i$ that is unique modulo the additive relations 
$P(\Delta)$.  Hence 
the lifts $Y_i$ of the Seidel elements $y_i$ are determined by the linear relations $P(\Delta)$.   
The other information needed to determine 
the multiplicative structure of $\QH^*(M)$ is the set of primitive classes $I$.
Thus, in the NEF case,
once one knows the Seidel elements $\Ss(\La_i), i = 1,\dots,N$,
there is an easy formula based on the combinatorics of 
its moment polytope $\Delta$ for
the multiplicative relations in
the quantum cohomology ring. This 
 substitution of the $Y_i$ for
the $x_i$ in the Stanley--Reisner ring $SR_Y$ is one way of 
looking at  Givental's change of 
variable formulae as discussed in~\cite[11.2.5.2]{CK}.
\QED\end{example}

\begin{rmk}\label{rmk:Novik}\rm
    Often one wants to consider quantum cohomology with coefficients 
    in a completion of the group ring of $H_2^S(M)$ rather than 
    of a quotient of $H_2^S(M)$. Our methods give similar results in
    this case, but one must use a slightly different 
    version of the Seidel representation.  For more details see
    McDuff--Salamon~\cite[Chapter~11.4]{MS2}.
    \end{rmk}

\subsection{The Seidel representation: examples}\label{ss:ex}

The examples in this section show that even in the case of the 
simplest manifolds, namely rational ruled symplectic 
$4$-manifolds, the Seidel element can be quite 
complicated.
The first example is Fano.   We show how lower order terms may appear in the
formula for $\Ss(\La_K)(a)$ and  
discuss a circle action with at most twofold isotropy.
The second 
example illustrates the NEF case, in which, as already noted in 
Seidel~\cite{Sei}, the expression for
$\Ss(\La)$ can have infinitely many nonzero terms. 
We also show what
can happen when the isotropy 
has order greater than two.

\begin{example}[The one point blowup of $\CP^2$]\label{ex:CPB}\rm
Fix $\mu \in (0,1)$.
Identify the one point blow up $M_*$ of $\CP^2$  with the region
$$
\left\{(z_{1}, z_{2})\in \C^{2} \left|\ 
\frac{\mu^{2}}{\pi}\le |z_{1}|^{2} + |z_{2}|^{2}\le \frac 1{\pi}
\right. \right\}
$$
with boundaries collapsed along the Hopf flow, and give it the
corresponding symplectic form $\om_{\mu}$.  
Let $E \in H_2(M_*)$ denote the class of the exceptional divisor, let $L = [\CP^1]$,
and let $B = L - E$ be the fiber class. 
Thus $\om_{\mu}(L) = 1$.
Let $p \in H_0(M)$ denote the homology class of a point,
and let $\1$ be the generator of $H_4(M)$.
The space $M_*$ is a toric variety,
where $T = S^1 \times S^1$  acts on $M_*$ by
$(\alpha_1,\alpha_2) \cdot (z_1,z_2) = (\alpha_1 z_1, \alpha_2 z_2)$.
The standard complex structure $J$ on $M_*$
is $T$-invariant and is compatible with $\omega_{\mu}$.
The  moment map $\Phi: M \to \R^2$ is given  by
$$
\Phi(z_1,z_2) = (|z_1|^2 - \epsilon, |z_2|^2 - \epsilon), \quad
\mbox{where} \quad  \epsilon =  \frac{1-\mu^6}{3(1- \mu^4)}.
$$
The primitive outward normals are 
$$
\eta_1 = (-1,0),\ \eta_2 = (0,-1), \
\eta_3 = (1,1),\  \mbox{and} \ \eta_4 = (-1,-1).
$$

Let $\Lambda_i$ denote the circle action corresponding to $\eta_i$.
Since the moment map for $\La_1$ takes its maximum on the set $z_1=0$,
$\La_1$ is the action $(z_1,z_2)\mapsto (e^{-2\pi i 
t}z_1,z_2)$.  Similar arguments give explicit formula for the 
other $\La_i$.
Since $(M_*,J)$ is Fano, 
part (iii) of Theorem~\ref{thm:max} implies that
$$
\Ss(\Lambda_1) = \Ss(\Lambda_2)  = B \otimes q t^\epsilon, \ \
\Ss(\Lambda_3) = L \otimes q t^{1 - 2 \epsilon}, \ \  \mbox{and}  \ \
\Ss(\Lambda_4) = E \otimes q t^{2 \epsilon - \mu^2}.
$$
There are two primitive subsets, namely $\{3,4\}$ and $\{1,2\}$.
Since $\eta_3 + \eta_4 = 0$,
$$
\1 = \Ss(\Lambda_3)*\Ss(\Lambda_4) = L * E \otimes q^2 t^{1 - \mu^2}.
$$
Since $\eta_1 + \eta_2  = \eta_4$,
$$ 
E \otimes q t^{2 \epsilon - \mu^2} = \Ss(\Lambda_4) = \Ss(\Lambda_1)
*\Ss(\Lambda_2) =  B * B \otimes q^2 t^{2 \epsilon}.
$$
Therefore  
\begin{equation}\labell{product}
B * B = E \otimes q^{-1} t^{-\mu^2} 
\quad \mbox{and} \quad
L * E = \1 \otimes q^{-2} t^{\mu^2 - 1}.
\end{equation}

The circle action $(\Lambda_1)^{-1}$  also has  a semifree maximum,
namely 
the point $[(0,1)] \in M_*$, the inverse image of the vertex $D_2\cap D_3$.
The holomorphic spheres $C$ through
$\Fmax$ all have $c_1(C) \geq 2$.  Hence, again applying
 part (iii) of Theorem~\ref{thm:max},
we conclude 
$$
\Ss(\Lambda_1^{-1}) = p \otimes q^2 t^{1 - \epsilon}.
$$
Since  $-\eta_1 = \eta_3 + \eta_2$, 
$$
p \otimes q^2 t^{1 - \epsilon} = \Ss((\Lambda_1)^{-1}) =
\Ss(\Lambda_3) \Ss(\Lambda_2) = B * L \otimes q^2 t^{1 - \epsilon}.
$$
Therefore,
\begin{equation}\labell{product2}
B * L = p.
\end{equation}
Note that equation (\ref{product}) determines $\QH_*(M_*)$
as a ring, but does not determine the product above.  Together,
equations (\ref{product}) and (\ref{product2}) determine all possible 
products in $\QH_*(M_*).$ 
In particular, using associativity, we find
$$
p*p = L\otimes q^{-3}t^{-1},\quad  E * p  = B \otimes q^{-2}t^{\mu^2 - 1},
\quad \mbox{and}
\quad p*B = \1\otimes q^{-3}t^{-1}.
$$
These products 
may also be derived directly 
from the $3$-point Gromov--Witten invariants:
it is not hard to check  
that the only nonzero invariants
involving the classes $p,B,$ and $E$ are
\begin{eqnarray*}
\GW_{L,3}(p, p, B) = 1;&& \GW_{B,3}(p, E, E) = 1; 
\quad \mbox{and} 
\\
\GW_{E,3}(A_{1}, A_{2}, A_{3}) = \pm 1&&\mbox{where }\;\; A_{i} = E\mbox{
or }B.
\end{eqnarray*}

The natural action of $U(2)$ on $\C^2$ induces an action on $M_*$;
this action contains the torus $T$.
Since $\pi_1(U(2)) = \Z$, this shows that, as elements of
$\pi_1(\Symp(M_*,\omega))$, $\Lambda_1 = \Lambda_2$. 
Hence 
$$
\Lambda_3 = (\Lambda_4)^{-1} = {\Lambda_1}^{-2} = {\Lambda_2}^{-2}
$$
It is a worthwhile exercise to check that 
$\Ss(\Lambda_3) = \Ss(\Lambda_4)^{-1} = \Ss(\Lambda_1)^{-2}= 
\Ss(\Lambda_2)^{-2}.$

Since $\Lambda_4$ is semifree, we 
can also apply Theorem~\ref{thm:semifree} to this action.
It has $F_{\max} = E$, the exceptional divisor.
Let $r \in H_0(E)$  be the homology
class of a point. Then the downwards extension
$r^- = B \in H_2(M)$, and the upwards extension $r^{+} = p \in H_0(M)$.
Then 
$$ 
\Ss(\Lambda_4)(r^-) = (E \otimes qt^{2\epsilon - \mu^2}) * B
= r^+  \otimes q t^{2 \epsilon - \mu^2} - 
E \otimes  (q t^{2 \epsilon - \mu^2})(q^{-1} t^{-\mu^2}).
$$
This agrees with Theorem~\ref{thm:semifree}, but also
shows that  lower order terms can appear, even in
this simple example. 
This lower order term comes from an invariant chain consisting of 
the sphere  $F_{\max}$ (in class $E$) together with a section $\si_z$ for $z\in 
F_{\max}$.

Now consider the  circle action $\Lambda'$ corresponding to
$\eta_1 +\eta_4 = (-2,-1)$.  
The corresponding moment map has a semifree maximum, namely
the point  $[(0,\mu)] \in M_*$ that maps down to $D_1\cap D_4\in\Phi(M)$.
Hence  part (i) of Theorem~\ref{thm:max} applies, but
part (iii) does not
because there is a holomorphic sphere $E$ through the maximum 
with $2c_1(E) = 2 \leq \codim \Fmax = 4$.  
Therefore our results do not rule out the presence of lower order 
terms in $\Ss(\Lambda')$ and indeed these exist:
since $(-2,-1) = \eta_1 + \eta_4$,  
$$
\Ss(\Lambda') = \Ss(\Lambda_1)*\Ss(\Lambda_4) =
B * E \otimes q^2 t^{3 \epsilon - \mu^2} = 
p \otimes q^2 t^{3 \epsilon - \mu^2}  - 
E \otimes q t^{3\eps-2\mu^2}.
$$
Observe also that $\Lambda'$ has at most twofold isotropy, with
isotropy submanifold $(M_*)^{\Z/(2)}$  equal to $\Phi^{-1}(D_2)$.
One can check this by writing $(-2,-1) = 2\eta_4  - \eta_2$: as 
explained in the proof of Proposition~\ref{prop:torsim}
in \S\ref{ss:Srep}  the coefficient of $-\eta_4$ in this 
expression equals the weight on the transverse edge $D_2$.
Therefore Theorem~\ref{thm:pfasimplemax} applies to the fixed 
components $F_{13}: = \Phi^{-1}(D_1\cap D_3)$ and
$F_{24}: = \Phi^{-1}(D_2\cap D_4)$, which are both isolated points.
Since $F_{13}$ is semifree 
the Euler class $e(F_{13})$ is nonzero.  On the other hand,
 $e(F_{24})= 0$.
Further, if $c_{ij}\in H_0(F_{ij})$ denotes a generator, we find
$(c_{13})^{-} = L$, while $(c_{24})^{-} = E$.
Therefore Theorem~\ref{thm:pfasimplemax} implies that
$
\Ss(\La')(L)$ has a nontrivial summand 
$c_{0,0}\otimes t^{K'(F_{13})}$ where
$c_{0,0}\cdot L = 1$ and $K'$ denotes the moment map 
$\Phi^{\eta_1+\eta_4}$ of $\La'$.  (This is the contribution to $
\Ss(\La')(L)$ of the constant section at the homologically visible
point $F_{13}$.)  On the other hand, because $F_{24}$ is not
homologically visible,
the constant section at $F_{24}$ makes no contribution to 
$\Ss(\La')(E)$ and so
the coefficient of $q t^{K'(F_{24})}$ in $
\Ss(\La')(E)$ vanishes.  This can be checked by direct calculation.
For example $K'(F_{13}) = 3\eps -1$ and
$
\Ss(\La')(L) = \Ss(\La')(E+B)$ contains one nonzero term of the form 
$a\otimes  t^{\ka}$, namely $B\otimes t^{3\eps-1}$.

The manifold $M_*$ has many other toric structures; correspondingly 
there are many other elements of $\pi_1(\Ham(M_*, \om_{\mu}))$
that are represented by semifree
circle actions.  Indeed, whenever $\mu^2 > k/(k+1)$, there is
an $\om_\mu$-compatible  complex structure $J_k$ on $M_*$ such that
the underlying complex manifold $(M_*, J_k)$ can be identified with the
 projectivization
${\mathbb P}(L_k \oplus \C)$, where $L_k$ is the holomorphic line bundle
over $\CP^1$ with Chern class $2k+1$.  The loop that rotates
the fibers of $\C$ by $e^{2\pi it}$ is semifree and represents the
class $(4k+2)\al$, where $\al = [-\La_1]\in \pi_1(\Ham(M_*, \om_{\mu}))$.\footnote
{
The formula in~\cite{AM} Lemma 2.11(i) is  slightly incorrect.}
The  classes $(2k+1)\al$ are also represented by
circle actions that preserve $J_k$ and
rotate the base of the ruled surface $(M_*, J_k)$.
When $k=0$, $J_0$ is the standard complex structure discussed above, and 
the representative for $2\al$ is $\La_3$ while the 
representative for $\al$ is $\La_1^{-1}$.
When $k>0$ explicit formulas for these actions can be derived from the
description of  $(M_*, J_k)$ as a toric manifold given in~\cite{AM}~\S2.3.  
In this case these actions
have $4$ isolated fixed points, two each in the fibers lying above the
fixed points of the base rotation.  However, these fixed points are
not semifree except when $k=1$, in which case the fixed points in the
fiber containing the overall minimum are semifree.
In the next example we shall discuss a similar action on 
$S^2\times S^2$ in detail.
\end{example}

\begin{example}[Circle actions on $S^2 \times S^2$.]\label{ex:S2S2}\rm
Consider $M = \CP^1 \times \CP^1$ with the symplectic form
$\omega_\mu = \mu\pi_1^*(\si) + \pi_2^*(\si)$, where $\pi_i$
is projection onto the $i$'th factor, $\si$ is the standard
symplectic form on $\CP^1$ with total area $1$.
Assume that  $\mu \ge 1$.
Define $A$ and $B$ in $H_2(M)$ by $A = [\CP^1 \times \{q\}]$  
and $B = [ \{q\} \times \CP^1]$, where $q \in \CP^1.$ 
Note that $\omega_\mu(A) = \mu$ and $\omega_\mu(B) = 1$.
Let $p \in H_0(M)$ denote the homology class of a point,
and let $\1$ denote the generator of $H_4(M)$.

The standard action of the torus $T = S^1 \times S^1$ on $(M,\om_\mu)$
in which each $S^1$-factor rotates the corresponding sphere
has moment map
$\Phi: M \to \R^2$ given by 
$$
\Phi([x_1:x_2],[y_1:y_2]) = 
\left(\mu\frac{|x_1|^2 - |x_2|^2}{|x_1|^2 +|x_2|^2}, 
\frac{|y_1|^2 - |y_2|^2}{|y_1|^2 + |y_2|^2}\right).
$$
The primitive outward normals to the moment
image $\Delta = \Phi(M)$ are 
$$
\eta_1 = (1,0),\ \eta_2 = (-1,0), \
\eta_3 = (0,1), \ \mbox{and} \ \eta_4 = (0,-1).
$$
Let $\Lambda_i$ be the circle action associated to $\eta_i$.

Since the standard complex structure on $\CP^1 \times \CP^1$
is Fano and $T$-invariant, and since $\Lambda_i$ acts semifreely
for all $i$, by Theorem~\ref{thm:max}
$$
\Ss(\Lambda_1) = \Ss(\Lambda_2) =  B \otimes q t^{\frac{\mu}{2}},
\quad \mbox{and} \quad
\Ss(\Lambda_3) = \Ss(\Lambda_4) = A \otimes q t^{\frac{1}{2}}.
$$
Since $\eta_1 + \eta_2 = 0$ and $\eta_3 + \eta_4 = 0$,
$\Ss(\Lambda_1) *\Ss(\Lambda_2) = \1$ and 
$\Ss(\Lambda_3) * \Ss(\Lambda_4) = \1$.
This implies that  
$$
B * B = \1 \otimes  q^{-2} t^{-\mu} \quad \mbox{and} \quad 
A * A = \1 \otimes q^{-2} t^{-1}.
$$

Let $\Lambda' \subset S^1 \times S^1$ be the circle  associated to 
$\eta_1 + \eta_3: = (1,1)$.
Then $\Lambda'$ acts by the diagonal action and so is semifree.
Since $c_1(C) \geq 2$ for every holomorphic sphere $C$,
$$
\Ss(\Lambda') = p \otimes q^{2} t^{\frac{1 + \mu}{2}}.
$$
Because $\Ss(\Lambda') = \Ss(\Lambda_1)* \Ss(\Lambda_3)$ we find
$
A * B = p.
$
As before, these products determine all the products in $\QH_*(M)$.
In particular
$$
p*A = B\otimes q^{-2}t^{-1}\quad \mbox{and} \quad p*B = A\otimes q^{-2}t^{-\mu}.
$$

We now describe a second toric structure on $M$.
Let $L_n$ denote the holomorphic bundle over $\CP^1$ with Chern class $n$.
Let $M'$ be the projectivization of the bundle $L_2 \oplus \C$.
Two commuting circles act naturally on $M'$.
First, the standard circle action on $\CP^1$ lifts
naturally to an action on $T^*(\CP^1) = L_2$, and hence to $M'$.
Denote this circle action by $\Gamma'$.
Another circle, 
say $\Ga''$,
acts by  rotating each fiber.
The standard complex structure $J_2$ on $M'$ is invariant under
the resulting $S^1 \times S^1$-action.
Moreover, 
if we assume that $\mu> 1$,
there exists a $J_2$-compatible invariant symplectic form $\omega$
on $M'$ so that $M'$ is symplectomorphic to $(\CP^1 \times 
\CP^1,\omega_\mu)$,
which we consider to be fibered over $\CP^1$ via projection to the 
first factor.
In fact, we may assume that this symplectomorphism
lifts the identity map on the base $\CP^1$ and
is equivariant with respect to the action of $\Lambda'$ on $M$
and $\Gamma'$ on $M'$:
see for example~\cite{AM}.
Hence, we immediately conclude
$$
\Ss(\Gamma') = \Ss(\Lambda') = p \otimes q^2 t^{\frac{1 + \mu}{2}}.
$$
Here, and elsewhere, we 
identify $p$, $A$, $B$, and $\1$ with their image in $H_*(M')$.

As described above, $M'$ is a smooth toric  variety with
moment map $\Phi'$.
The moment image $\Delta' = \Phi'(M')$ 
is a quadrilateral with
outward normals
$$
\gamma_1 = (0,1), \ \gamma_2 =  (0,-1),\  \gamma_3 = (1,-1), \
\mbox{and} \ \gamma_4 = (-1,-1).
$$
Further $(\Phi')^{-1}(D_1)$ is the diagonal in $M'\equiv \CP^1\times \CP^1$
and so contains points  that we will call $v_{ss}$ and $v_{nn}$, 
where $v_{ss}$ corresponds to $([0:1],[0:1])$, the pair (south pole, 
south pole),
in $\CP^1\times \CP^1$  and $v_{nn}$ corresponds to (north pole, north 
pole).   Similarly $(\Phi')^{-1}(D_2)$ is the antidiagonal and 
contains $v_{sn}, v_{ns}$.  Indeed
$$
\begin{array}{cc}
(\Phi')^{-1}(D_1\cap D_3) = v_{nn},& (\Phi')^{-1}(D_3\cap D_2) = 
v_{ns},\\ (\Phi')^{-1}(D_2\cap D_4) = v_{sn},&(\Phi')^{-1}(D_4\cap D_1) = v_{ss}.\end{array}
$$
The moment image itself is
$$ 
\Delta' = \{ \alpha \in \R^2  \mid (\alpha,\gamma_i) \leq c_i \},
\quad \mbox{where} 
$$
$$ 
c_1 = \frac{1}{2}  + \frac{\mu}{2} - \epsilon,\;\;\;  c_2 = \epsilon + \frac{1}{2}
- \frac{\mu}{2}, \;\;\; c_3 = c_4 = \epsilon,\;\;\;
\mbox{and} \quad
\epsilon = \frac{\mu}{2} + \frac{1}{6 \mu}.
$$
Let $\Gamma_i \subset S^1 \times S^1$ be the circle  associated
to $\gamma_i$.
In our previous notation, $\Ga' = \Ga_1 + \Ga_3$ and $\Ga'' = \Ga_1$.

The circle  $\Gamma_1$ acts semifreely.
Since every holomorphic sphere $C$ which
intersects $\Fmax$ has $c_1(C) \geq 2$,
it follows from part (iii) of Theorem~\ref{thm:max} that
there are no lower order terms in $\Ss(\Gamma_1)$. 
Since $[\Fmax] =
[(\Phi')^{-1}(D_1)]= 
A + B$, 
$$
\Ss(\Gamma_1) = (A + B) \otimes q 
t^{\frac{1}{2} + \frac{\mu}{2} - \epsilon}.
$$

The circle $\Gamma_2$ also acts semifreely.
In this case, $\Fmax$ itself is a holomorphic sphere in class $A - B$, 
so $c_1(\Fmax) = 0$.
Therefore, part (iii) of Theorem 1.10 does not exclude lower
order terms.  On the other hand, every holomorphic
sphere $C$ with $c_1(C) \leq 1$ lies entirely in $\Fmax$,
so every  term which contributes comes
from a $C$ which lies in   $\Fmax$.
Indeed, 
since $\gamma_1  = -\gamma_2$,
$$
\Ss(\Gamma_2) = \Ss(\Gamma_1)^{-1} 
=  (A - B) \otimes 
\frac{q t^{\frac{1}{2} - \frac{\mu}{2} + \epsilon}}{1 - t^{1-\mu}} = 
(A - B) \otimes  q t^{\frac{1}{2} - \frac{\mu}{2} + \epsilon}
\left( 1 + t^{1 - \mu} +  t^{2(1 - \mu)} + \cdots \right).
$$  
This calculation  also appears in Remark~11.5 of~\cite{Sei}.

Now consider $\Gamma_3$.
Once again,  Theorem 1.10 does not rule out lower order terms.
Since $\gamma_3 = \gamma_2 + (1,0)$,
$$
\Ss(\Gamma_3) = \Ss(\Gamma_2)* \Ss(\Gamma')
=   B \otimes q t^{\epsilon}
- (A - B) \otimes q t^{\epsilon} \frac{t^{1-\mu}}{1 - t^{1-\mu}}. 
$$
A similar argument applies to $\Gamma_4$.

Let's now pause for a moment to compare these
results with the previous section.
As above, let
$D_i$ denote the facet  that corresponds to $\gamma_i$;
let $x_i$ denote the Poincare dual of $\Phi^{-1}(D_i)$;
note that $[(\Phi')^{-1}(D_3)] = [(\Phi')^{-1}(D_4)] = A$,
$[(\Phi')^{-1}(D_1)] = A + B$, and $[(\Phi')^{-1}](D_2) = A - B$.
Converting the equations above into cohomology,
and using this notation, we find:
\begin{eqnarray*}
\Ss^*(\Gamma_1) & = & x_1 \otimes q^{-1} 
t^{- \frac{1}{2} - \frac{\mu}{2} + \epsilon},\\
\Ss^*(\Gamma_2) 
& = & x_2 \otimes \frac{q^{-1} t^{-\frac{1}{2} + \frac{\mu}{2} - \epsilon}}
{1 - t^{\mu-1}}, \\
\Ss^*(\Gamma_3)  
& = &  \left( x_3 -  x_2 \otimes \frac{t^{\mu - 1}}{1-t^{\mu-1}} \right) 
q^{-1} t^{-\epsilon}
 ,\quad \mbox{and}\\
\Ss^*(\Gamma_4)  
& = & \left( x_4 -  x_2 \otimes \frac{t^{\mu - 1}}{1-t^{\mu-1}} \right) 
q^{-1} t^{-\epsilon}.
\end{eqnarray*}
Thus in equation~(\ref{eq:Yi}) we may take
$$
Y_1 =  x_1, \quad  Y_2 = x_2 \otimes \frac{1}{1 - t^{\mu - 1}}, \quad
\mbox{and} \quad  Y_4 = Y_3 = x_3  - x_2 \otimes \frac{ t^{\mu - 1}}{1 - t^{\mu - 1}}.
$$
 
We now look at $\Ss(\Tilde\Ga)$ for the circle action $\Tilde\Ga$ 
 associated with $\ga=(1,2)$, which has 
threefold isotropy.
In notation introduced earlier, we can describe the fixed set of 
$\Tilde{\Gamma}$ as consisting of
the points $v_{nn}$ (the maximum), the saddle points
$v_{ss}, v_{ns}$ and the minimum $v_{sn}$.
The maximum is not semifree; in fact, 
because $(1,2) = \eta_1 + 
3\eta_3$, the diagonal $(\Phi')^{-1}D_1$ is
stablized 
by $\Z/(3)$.
Since the action does not have at most twofold
isotropy, the arguments of 
Theorems~\ref{thm:semifree} and~\ref{thm:pfasimplemax}
do not apply.  We show that the conclusions
of these theorems also fail. 
Since $(1,2) = 2 \gamma_1 + (1,0)$,
$$
\Ss(\tilde{\Gamma}) = \Ss(\Gamma_1)^2 * \Ss(\Gamma') = 
(p + p \otimes t^{1 - \mu} + 2 q^{-2} t^{-\mu}) \otimes 
q^2 t^{ \frac{1}{2} + \frac{3\mu}{2} - 2\epsilon}.
$$

First consider the minimum $v_{sn}$ which is semifree.
Then, in the notation of Theorem~\ref{thm:pfasimplemax}, $(v_{sn})^- = p$
and $(v_{sn})^+ = \1$ and so one might expect the
leading order term of $\Ss(\Tilde{\Gamma})(p)$ to come from the
section $\sigma_{sn}$ and so have the form $\1 \otimes q^d t^\kappa$.
But
$$
\Ss(\Tilde{\Gamma})(p) = 
(\1 \otimes q^{-4}t^{-1-\mu}  + \1 \otimes q^{-4} t^{ - 2\mu} + 2p \otimes q^{-2} t^{-\mu}) \otimes 
q^2 t^{ \frac{1}{2} + \frac{3\mu}{2} - 2\epsilon}
$$
has the leading order term $p \otimes t^{\frac{1}{2} + \frac{\mu}{2} 
- 2\epsilon}$.
It is not hard 
to check that this term comes from the invariant chain
$$
x =v_{sn} \stackrel{A - B}\to v_{ns} \stackrel{2B}\to v_{nn}
\stackrel{\si_{nn}}\to 
y = v_{nn},
$$
where $\si_{nn}$ is the constant section at $v_{nn}$.
This lies in class $A + B + [\si_{nn}] = A+B+ [\si_{nn}] = 
[\si_{sn}] -B$ since $[\si_{sn}] - [\si_{nn}] = A+2B$.

Next consider the semifree saddle point $v_{ss}$.  Then $(v_{ss})^- = B$ and 
$(v_{ss})^+ = A+B$. Therefore, from Theorem~\ref{thm:semifree} one 
would expect the leading order term in $\Ss(2\gamma + 
\tau_1+\tau_2)(B)$ to be $(A+B)\otimes t^{K(v_{ss})}$,
while in fact it is 
$(A+2B)\otimes t^{K(v_{ss})}$.  Since $(v_{ns})^+ = B$, one can get this extra 
term from an invariant chain going from $x\in (v_{ss})^-$ to $y\in 
(v_{ns})^-$  that 
lies in class $[\si_{ss}]$.  Since 
$[\si_{nn}] = [\si_{ss}] -(A+B)$ such a chain is given by
$$
x = v_{sn}\in (v_{ss})^- \stackrel E\to v_{ns}\stackrel B\to v_{nn}
\stackrel {\si_{nn}}\to 
v_{nn} \stackrel B\to y= v_{ns}\in (v_{ns})^-.
$$
\end{example}


\begin{thebibliography}{999999999}

\bibitem{AM}  M. Abreu and D. McDuff, Topology of symplectomorphism
groups of rational ruled surfaces,  SG/9910057,
 {\it Journ. of Amer. Math. Soc.}  {\bf 13}, (2000) 971--1009.

 \bibitem{AB} D. Austin and P. Braam, Morse--Bott theory and 
 equivariant cohomology, in {\it The Floer Memorial Volume},
 Progress in Mathematics {\bf 133}, Birkh\"auser (1995).

  \bibitem{Ba} V. Batyrev, Quantum cohomology rings of toric manifolds,  {\it Ast\'erisque} {\bf 218} (1993), 9--34.
  
 \bibitem{CK} D. Cox and S. Katz: {\it Mirror Symmetry and Algebraic 
 Geometry}, Math. Surveys and Monographs vol 68, AMS, Providence (1999).


\bibitem{FO}  K. Fukaya and K. Ono, Arnold conjecture and
Gromov--Witten invariants, {\it Topology}  {\bf 38} (1999), 933--1048.

\bibitem{Gonz} E. Gonzalez, Quantum cohomology and $S^1$-actions with
         isolated fixed points, SG/0310114.

\bibitem{GRO}    M. Gromov,  Pseudo holomorphic
curves in symplectic manifolds,
       {\it  Inventiones Mathematicae\/}, {\bf 82} (1985), 307--47.

\bibitem{Lalf} F. Lalonde, A field theory for symplectic fibrations 
      over surfaces, {\it Geometry and Topology} {\bf 8} (2004)
      1189--1226.


\bibitem{LMP}  F. Lalonde, D. McDuff and L. Polterovich,
Topological rigidity of Hamiltonian
loops and quantum homology, {\it Invent. Math} {\bf 135}, 369--385
(1999)


\bibitem{LIUT} Gang Liu and Gang Tian, Floer homology and Arnold
         conjecture, {\it Journ. Diff. Geom.}, {\bf 49} (1998), 1--74.


\bibitem{LIUT2} Gang Liu and Gang Tian, On the equivalence of multiplicative
        structures  in Floer Homology and Quantum Homology,
        {\it Acta Math. Sinica} {\bf 15} (1999).

	      
\bibitem{Mcv}  D. McDuff, The virtual moduli cycle, 
        {\it Amer. Math. Soc. Transl.} (2) {\bf 196} (1999), 73--102


\bibitem{Mcq}  D. McDuff, Quantum homology of Fibrations over $S^2$,
       {\it International
        Journal of Mathematics}, {\bf 11}, (2000), 665--721.

	

 \bibitem{MS1}  D. McDuff and D. Salamon, {\it Introduction to
       Symplectic Topology}, 2nd edition (1998) OUP, Oxford, UK


\bibitem{MS2}  D. McDuff and D. Salamon, {\it $J$-holomorphic curves 
      and Symplectic Topology}, Amer. Math. Soc. Colloq Publications,
      to appear (2004)


\bibitem{Mch}  D. McDuff,
       Geometric variants of the Hofer norm, 
       SG/0103089,   {\it Journal of Symplectic Geometry},
       {\bf 1} (2002), 197--252.



\bibitem{MSlim} D. McDuff and J. Slimowitz, Hofer--Zehnder capacity
      and length minimizing Hamiltonian paths, SG/0101085, {\it Geom. 
      Topol.} {\bf 5} (2001), 799--830.

\bibitem{McTLie} D. McDuff and S. Tolman, On circle actions with
        semifree fixed points, preprint (2005).
        
\bibitem{McTtor} D. McDuff and S. Tolman,  Polytopes with 
       mass-linear functions, preprint (2005).

\bibitem{SAL}
      D. Salamon, Morse theory, the Conley index and Floer
      homology, {\it Bulletin of the London Mathematical Society\/}, {\bf
      22} (1990), 113--40.


\bibitem{SZ}
       D. Salamon, and E. Zehnder,  Morse theory for periodic
       solutions of Hamiltonian systems and the Maslov index. {\it
       Communications in Pure and Applied Mathematics\/}, {\bf 45} (1992),
       1303--60.



\bibitem{Sch} M. Schwarz, Equivalences for Morse homology, in {\it
       Geometry and Topology in Dynamics} ed M. Barge, K. Kuperberg,
       Contemporary Mathematics {\bf  246}, Amer. Math. Soc. (1999), 197--216.


\bibitem{Sei}
       P. Seidel, $\pi_1$  of symplectic automorphism groups
       and invertibles in quantum cohomology rings, {\it Geom. and     Funct.  Anal.} {\bf 7} (1997), 1046 -1095.


\bibitem{TW}  S. Tolman and J. Weitsman, The cohomology rings of 
      abelian symplectic quotients, {\it Communications in Analysis and Geometry}.
to appear.

\bibitem{Wei} A. Weinstein, Cohomology of symplectomorphism groups 
      and critical values of Hamiltonians, {\it Math Z.} {\bf 201} 
      (1989), 75--82.


\end{thebibliography}
\end{document}